# A new approach to stochastic McKean–Vlasov limits with low-regularity coefficients

*By* **Robert Alexander Crowell** *of ETH Zürich*

**Abstract:** The empirical measure flow of a McKean–Vlasov $n$-particle system with common noise is a measure-valued process whose law solves an associated martingale problem. We obtain a stability result for the sequence of martingale problems: all narrow cluster points of the sequence of laws solve the formally limiting martingale problem. Through the solution of the limiting problem, we are able to characterize the dynamics of limits of the empirical measure flows.
A major new aspect of our result is that it requires rather weak regularity assumptions for the coefficients, for instance a form of local continuity of the drift in the measure argument and ellipticity of the diffusion coefficient for the idiosyncratic noise. In fact, the formally limiting martingale problem may fail to have any solution if there are discontinuities in the measure argument, so that the stability property is in general not true under low regularity assumptions.
Our novel approach leverages an emergence of regularity property for cluster points of the empirical measure flow. This provides us with a priori analytic regularity estimates which we use to compensate for the low regularity of the drift.

**Robert Alexander Crowell**

rc@math.ethz.ch

Department of Mathematics
ETH Zürich
Rämistrasse 101
CH-8092 Zürich

**Acknowledgements:** This work is part of my PhD thesis written at ETH Zürich. I am deeply indebted to my advisor Martin Schweizer who gave me the freedom to pursue this line of research, and at the same time guided me with his immense experience throughout. I benefited greatly from his detailed feedback, and the arguments and the exposition in this manuscript improved markedly from his careful comments.
I would also like to thank participants at workshops held at ETH Zürich, Hammamet, CIRM Luminy, Oxford University and Imperial College London for their feedback.
Numerous discussions with my colleagues, in particular with my office-mate Florian Krach and with Andrew Allan, should not go unmentioned. They sharpened my reasoning.

## CONTENTS

## INTRODUCTION

This paper is the second in a series studying the limiting behavior of particle systems of McKean–Vlasov type with common noise when the coefficients have only low regularity. We characterize the dynamics satisfied by cluster points of the sequence $(\mu^n)_{n\in\mathbb{N}}$ of *empirical measure flows*, where $\mu^n = (\mu^n_t)_{t\in[0,T]}$ is associated with the $\mathbb{R}^d$-valued *n-particle system* given for $i \in \{1,\dots,n\}$ and $t \in [0,T]$ by

$$\begin{cases} X^{i,n}_t = x_0 + \int_0^t b_s(X^{i,n}_s, \mu^n_s)\,\mathrm{d}s + \int_0^t \sigma_s(X^{i,n}_s, \mu^n_s)\,\mathrm{d}B^{i,n}_s + \int_0^t \bar\sigma_s(X^{i,n}_s, \mu^n_s)\,\mathrm{d}Z^n_s\,, \\ \mu^n_t = \dfrac{1}{n}\sum_{i=1}^n \delta_{X^{i,n}_t}\,. \end{cases} \tag{1}$$

Above we call $B^{i,n}$ the *idiosyncratic noise* and $Z^n$ the *common noise*. The system (1) possesses a unique weak solution under very general assumptions on the coefficients. In addition, the process $\mu^n$ has continuous trajectories in the space $\mathcal{S}'(\mathbb{R}^d)$ of tempered distributions and is valued in the subset $\mathbf{M}^+_1(\mathbb{R}^d) \subseteq \mathcal{S}'(\mathbb{R}^d)$ of Borel probability measures.

Instead of the sequence of processes $(\mu^n)_{n\in\mathbb{N}}$, we study their laws $(\mathbf{P}^n_{x_0})_{n\in\mathbb{N}}$ on the space $\Omega := \mathbf{C}([0,T];\,\mathcal{S}')$ with the canonical coordinate process $\Lambda = (\Lambda_t)_{t\in[0,T]}$. For each $n \in \mathbb{N}$, the measure $\mathbf{P}^n_{x_0}$ solves the *canonical martingale problem* $\mathrm{MP}(\delta_{x_0}, \mathscr{A}, \mathscr{Q} + \frac{1}{n}\mathscr{C})$, where $\mathscr{A}$, $\mathscr{Q}$, $\mathscr{C}$ are certain operators on $\mathbf{M}^+_1(\mathbb{R}^d) \subseteq \mathcal{S}'(\mathbb{R}^d)$ specified in terms of the coefficients $b$, $\sigma$, $\bar\sigma$ in (1). In particular, the process $\Lambda$ can be written as

$$\Lambda_t = \delta_{x_0} + \int_0^t \mathscr{A}_s(\Lambda_s)\,\mathrm{d}s + M_t\,, \qquad \text{for } t \in [0,T] \text{ under } \mathbf{P}^n_{x_0}, \tag{2}$$

where $M = (M_t)_{t\in[0,T]}$ is an $\mathcal{S}'$-valued martingale with continuous trajectories under $\mathbf{P}^n_{x_0}$, whose quadratic variation is given by

$$\langle\!\langle M \rangle\!\rangle_t = \int_0^t \Big(\mathscr{Q}_s(\Lambda_s) + \frac{1}{n}\mathscr{C}_s(\Lambda_s)\Big)\,\mathrm{d}s\,, \qquad \text{for } t \in [0,T] \text{ under } \mathbf{P}^n_{x_0}. \tag{3}$$

Thus (2) can be viewed as an $\mathcal{S}'$-valued semimartingale decomposition of $\Lambda$, with characteristics given by $\mathscr{A}$ and $\mathscr{Q} + \frac{1}{n}\mathscr{C}$. It is also well known that the sequence $(\mathbf{P}^n_{x_0})_{n\in\mathbb{N}}$ is tight and therefore possesses cluster points.

Our main result in Theorem 1.6 proves under rather weak regularity assumptions on $b$, $\sigma$, $\bar\sigma$ that *any cluster point* $\mathbf{P}^\infty_{x_0}$ *of* $(\mathbf{P}^n_{x_0})_{n\in\mathbb{N}}$ *solves the formally limiting martingale problem* $\mathrm{MP}(\delta_{x_0}, \mathscr{A}, \mathscr{Q})$. In other words, the decomposition in (2) is stable under narrow convergence and thus remains valid under $\mathbf{P}^\infty_{x_0}$, with (3) replaced by its formal limit $\langle\!\langle M \rangle\!\rangle_t = \int_0^t \mathscr{Q}_s(\Lambda_s)\,\mathrm{d}s$. This result requires ellipticity of $\sigma$ and a type of local continuity in the measure argument $\mu$ of the drift $b$; see Assumptions 1.2 and 1.5 below. These assumptions are indispensable for solutions of $\mathrm{MP}(\delta_{x_0}, \mathscr{A}, \mathscr{Q})$ to exist, and at the same time generous enough to cover interesting examples that appear prominently in the literature; see Section 5 for examples and Section 6 for a discussion.

In (2), $M$ is a martingale under $\mathbf{P}^n_{x_0}$. If the coefficients $b$, $\sigma$, $\bar\sigma$ are jointly continuous in $(x,\mu)$, it is straightforward to show that $M$ remains a martingale under $\mathbf{P}^\infty_{x_0}$, and this readily yields a result as in Theorem 1.6. This works because the operators $\mathscr{A}$, $\mathscr{Q}$, $\mathscr{C}$ inherit continuity properties from the coefficients $b$, $\sigma$, $\bar\sigma$. However, the situation changes significantly if for example $b_t$ is not continuous in $x$ because then the operator $\mathscr{A}_t$ is in general not continuous either. This discontinuity presents a serious obstacle, and the literature has no strategy to deal with it. As a consequence, results are limited to specific cases; see Section 6 for a discussion.

The systematic approach we develop here is novel and provides a general strategy to address this discontinuity issue. It rests on an *emergence of regularity* property, namely a quantitative a priori regularity estimate for $\Lambda$ under cluster points $\mathbf{P}^{\infty}_{x_0}$ of $(\mathbf{P}^{n}_{x_0})_{n\in\mathbb{N}}$ which we obtained in Crowell [12]. With this estimate and an additional uniform regularity estimate for the family of $n$-particle systems in (1), we are able to compensate for the low regularity of $b_t$ in $(x, \mu)$ and deduce a certain form of *local continuity* of $\mathscr{A}_t$. This provides a firm basis to use narrow convergence arguments and obtain the stability of the decomposition in (2) under cluster points of $(\mathbf{P}^{n}_{x_0})_{n\in\mathbb{N}}$. We demonstrate the virtues of our approach by implementing it under Assumptions 1.2 and 1.5. However, we emphasize that its conceptual aspects carry over to more general assumptions and also to other settings. A case in point is the forthcoming work Crowell [13], where the stability result in Theorem 1.6 below serves as a stepping stone to obtain an existence result for solutions of $\mathrm{MP}(\delta_{x_0}, \mathscr{A}, \mathscr{Q})$ under more general assumptions. We preview this result in Section 6.

Solutions of $\mathrm{MP}(\delta_{x_0}, \mathscr{A}, \mathscr{Q})$ yield solutions to the associated stochastic Fokker–Planck equation, a certain SPDE with coefficients given in terms of $b$, $\sigma$, $\bar{\sigma}$. Via superposition theory, solutions of this SPDE are related to solutions of the associated common noise McKean–Vlasov SDE, linking the existence theory of both. The stability question, i.e., the view that $\mathrm{MP}(\delta_{x_0}, \mathscr{A}, \mathscr{Q})$ can be interpreted as the limit of $(\mathrm{MP}(\delta_{x_0}, \mathscr{A}, \mathscr{Q} + \frac{1}{n}\mathscr{C}))_{n\in\mathbb{N}}$, is connected to the classical propagation of chaos property. If available, this property justifies viewing a McKean–Vlasov SDE as a limiting model of the particle systems in (1) when $n \to \infty$. Both the solution theory for McKean–Vlasov SDEs and the propagation of chaos property have received considerable attention in the literature, mostly for regular coefficients. But an understanding of $\mathrm{MP}(\delta_{x_0}, \mathscr{A}, \mathscr{Q})$ as a limiting problem lies at the heart of obtaining existence and propagation of chaos results for McKean–Vlasov SDEs also with a low-regularity drift. Important examples that come to mind are instances when $b$ is given in terms of an irregular or singular kernel—a case which arises often in applications. We give an example of an irregular kernel in Section 5. In Section 6, we provide an outlook on how extensions of our approach can be used to treat more general singular kernels, e.g. unbounded integrable kernels or distributional kernels.

Another motivation for martingale problems with low regularity comes from mean-field games and control in their weak formulations. Consider the case of drift control as an example; and more specifically the case in which the dynamics (2) are controlled via $b$ which appears in $\mathscr{A}$. Since the regularity of optimal policies cannot be taken for granted and these policies feed back into $b$, the operator $\mathscr{A}$ may inherit poor regularity. Moreover, optimal policies may depend on $n$, and so we should expect added difficulties due to an explicit dependence on $n$ in (2), leading to a sequence $(\mathrm{MP}(\delta_{x_0}, \mathscr{A}^n, \mathscr{Q} + \frac{1}{n}\mathscr{C}))_{n\in\mathbb{N}}$ of controlled martingale problems. By the same reasoning, we have for the hypothetical case $n = \infty$ a controlled problem $\mathrm{MP}(\delta_{x_0}, \mathscr{A}^{\infty}, \mathscr{Q})$. This no longer has heuristic meaning as a formal limiting problem because it is unclear how $\mathscr{A}^{\infty}$ should be understood as a limit of $(\mathscr{A}^n)_{n\in\mathbb{N}}$. All this not only emphasizes the stability question but, also raises the more fundamental existence question for $\mathrm{MP}(\delta_{x_0}, \mathscr{A}^{\infty}, \mathscr{Q})$; see Section 6 for a discussion.

### Organization

This paper is structured as follows. The probabilistic setup, the main result in Theorem 1.6 and the outline of its proof are contained in Section 1. Section 2 develops the basic technical terminology and studies the martingale problem $\mathrm{MP}(\delta_{x_0}, \mathscr{A}, \mathscr{Q} + \frac{1}{n}\mathscr{C})$. Section 3 develops the necessary regularity results. In particular, we recall the emergence of regularity. Section 4 combines all these ingredients into the proof of the main result in Theorem 1.6. Finally, Sections 5 and 6 concludes this work with examples, and a discussion and an outlook, respectively. The Appendices A and B contain useful results from functional analysis and probability, as well as proofs omitted from the main body.

### Notation

We collect the essential notation. Appendix A contains more background and references.

- Throughout $c(a_1, \dots, a_k)$ denotes a generic finite constant which depends on the quantities $a_1, \dots, a_k$, and which may change from line to line.

- The space of *Schwartz functions* is denoted by $\mathcal{S} = \mathcal{S}(\mathbb{R}^d)$ and endowed with its usual Fréchet topology. We write $\mathcal{S}' = \mathcal{S}'(\mathbb{R}^d)$ for the topological vector space of *(tempered) distributions*, and $\lambda[\phi] := \langle \lambda; \phi \rangle_{\mathcal{S}' \times \mathcal{S}}$ for the duality pairing of $\lambda \in \mathcal{S}'$ and $\phi \in \mathcal{S}$. The Fourier transform $\mathscr{F}$ is a linear automorphism in $\mathcal{S}'$ with inverse $\mathscr{F}^{-1}$.

- For a measure-space $(\mathrm{X}, \mathcal{X}, \chi)$ and a Banach space $(B, \|\cdot\|_B)$, we denote the associated *Bochner space* by $\mathbb{L}^p(\chi; B) = \mathbb{L}^p((\mathrm{X}, \mathcal{X}, \chi); B)$. If $B = \mathbb{R}$, we write $\mathbb{L}^p(\chi) = \mathbb{L}^p(\mathrm{X}, \mathcal{X}, \chi)$ for the usual *Lebesgue space*.

- The *Bessel potential space* with integrability $r \in (1, \infty)$ and regularity $s \in \mathbb{R}$ is defined by $\mathsf{H}^s_r = \mathsf{H}^s_r(\mathbb{R}^d) := \{f \in \mathcal{S}'(\mathbb{R}^d) : \|f\|_{\mathsf{H}^s_r(\mathbb{R}^d)} < \infty\}$, where $\|f\|_{\mathsf{H}^s_r(\mathbb{R}^d)} := \|J^s f\|_{\mathbb{L}^r(\mathrm{d}x)}$ with $J^s f = \mathscr{F}^{-1}(h^s(\mathscr{F} f))$ for $f \in \mathcal{S}'$ and $h(\xi) := \langle \xi \rangle := (1 + |\xi|^2)^{1/2}$.

- The *Hermite–Fourier space* with regularity $p \in \mathbb{R}$ is the Hilbert space obtained by endowing $\mathscr{H}_p = \mathscr{H}_p(\mathbb{R}^d) := \{f \in \mathcal{S}'(\mathbb{R}^d) : \|f\|_{\mathscr{H}_p(\mathbb{R}^d)} < \infty\}$ with the (semi-)norm $\|f\|_{\mathscr{H}_p(\mathbb{R}^d)} := (\sum_{k \in \mathbb{N}^d} (\langle k \rangle^p f_k^{\#})^2)^{1/2}$, where $f_k^{\#}$ is the $k$'th Hermite coefficient of $f \in \mathcal{S}'$.

- If $(\mathrm{X}, \mathcal{X})$ is a measurable space then the set of *probability measures* is denoted by $\mathbf{M}_1^+(\mathrm{X}) = \mathbf{M}_1^+(\mathrm{X}, \mathcal{X})$. For $\mu, \nu \in \mathbf{M}_1^+(\mathrm{X})$ with $\int_{\mathrm{X}} |x| \mathrm{d}\mu, \int_{\mathrm{X}} |x| \,\mathrm{d}\nu < \infty$, we define the *Kantorovich–Rubinstein* distance by $\mathsf{d}_{\mathbb{W}_1}(\mu, \nu) := \sup\{\int f \,\mathrm{d}(\mu - \nu) : \|f\|_{\mathsf{Lip}} \leq 1\}$. We write $\mathcal{P}_{\mathrm{wk}^*} = \mathcal{P}_{\mathrm{wk}^*}(\mathbb{R}^d)$ for the set $\mathbf{M}_1^+ = \mathbf{M}_1^+(\mathbb{R}^d)$ with the *narrow topology*.

- Note that $\mathbf{M}_1^+ \subseteq \mathcal{S}'$ and that $\mathcal{P}_{\mathrm{wk}^*} \hookrightarrow \mathcal{S}'$, which is to say that the embedding is continuous. In particular, if $\mu \in \mathbf{M}_1^+$ and $\phi \in \mathcal{S}$, then $\mu[\phi] = \langle \mu; \phi \rangle_{\mathcal{S}' \times \mathcal{S}} = \int_{\mathbb{R}^d} \phi(x) \,\mu(\mathrm{d}x)$.

# 1 SETUP AND MAIN RESULT

## 1.1 Probabilistic setup for the particle systems

Fix $T < \infty$ and $d \in \mathbb{N}$. The basis for our study is the family of *n-particle systems* of McKean–Vlasov type formalized for each $n \in \mathbb{N}$, $i = 1, \ldots, n$ and $t \in [0, T]$, by

$$X_t^{i,n} = x_0 + \int_0^t b_s(X_s^{i,n}, \mu_s^n)\,\mathrm{d}s + \int_0^t \sigma_s(X_s^{i,n}, \mu_s^n)\,\mathrm{d}B_s^{i,n} + \int_0^t \bar{\sigma}_s(X_s^{i,n}, \mu_s^n)\,\mathrm{d}Z_s^n\,, \tag{1.1}$$

$$\mu_t^n = \frac{1}{n}\sum_{i=1}^n \delta_{X_t^{i,n}}\,. \tag{1.2}$$

The $n$-particle system (1.1) consists of $n$ coupled stochastic differential equations (SDEs), each with values in $\mathbb{R}^d$. We refer to $X^{i,n}$, the $i$'th component of the system, as a *particle*, to $B^{1,n}, \ldots, B^{n,n}$ as the *idiosyncratic noise* and to $Z^n$ as the *common noise*. To state the assumptions on the coefficients, we need on $[0, T] \times \mathbb{R}^d \times \mathbf{M}_1^+(\mathbb{R}^d)$ the product $\sigma$-algebra $\mathcal{E}$ generated by the Lebesgue-measurable sets of $[0, T] \times \mathbb{R}^d$ and the Borel-measurable sets of $\mathcal{P}_{\mathrm{wk}^*}(\mathbb{R}^d)$.

**Assumption 1.1** | The drift coefficient $b : [0, T] \times \mathbb{R}^d \times \mathbf{M}_1^+(\mathbb{R}^d) \to \mathbb{R}^d$ is bounded and measurable relative to $\mathcal{E}$.

**Assumption 1.2** | The diffusion coefficients $\sigma : [0, T] \times \mathbb{R}^d \times \mathbf{M}_1^+(\mathbb{R}^d) \to \mathbb{R}^{d\times d}$ as well as $\bar{\sigma} : [0, T] \times \mathbb{R}^d \times \mathbf{M}_1^+(\mathbb{R}^d) \to \mathbb{R}^{d\times m}$ are bounded and measurable relative to $\mathcal{E}$. In addition,

(E) there is $\kappa > 0$ such that the uniform ellipticity condition

$$z^\top\Big((\sigma_t\sigma_t^\top)(x, \mu)\Big)z \geq \kappa\,|z|^2$$

holds for all $z \in \mathbb{R}^d$, $t \in [0, T]$ and $(x, \mu) \in \mathbb{R}^d \times \mathbf{M}_1^+(\mathbb{R}^d)$;

(H) there exist $\beta > 0$ and a constant $C < \infty$ such that the Hölder-regularity conditions

$$|\sigma_t(x, \mu) - \sigma_t(x', \mu')| \leq C\Big(|x - x'|^\beta + \mathsf{d}_{\mathbb{W}_1}(\mu, \mu')^\beta\Big)\,,$$
$$|\bar{\sigma}_t(x, \mu) - \bar{\sigma}_{t'}(x', \mu')| \leq C\Big(|t - t'|^\beta + |x - x'|^\beta + \mathsf{d}_{\mathbb{W}_1}(\mu, \mu')^\beta\Big)\,,$$

hold for all $t, t' \in [0, T]$ and $(x, \mu), (x', \mu') \in \mathbb{R}^d \times \mathbf{M}_1^+(\mathbb{R}^d)$.

The $n$-particle system (1.1) is solved in the usual weak sense of SDEs; cf. Stroock and Varadhan [63, Ch. 6]. To construct such a solution, we write

$$\Omega_{\mathrm{par}}^n := \mathbf{C}([0, T];\, \mathbb{R}^{nd} \times \mathbb{R}^m)$$

with the coordinate process representing the particles and the common noise, i.e.,

$$(t, \omega) \mapsto \omega_t =: (X_t^n, Z_t^n)(\omega) = (X_t^{1,n}, \ldots, X_t^{n,n}, Z_t^n)(\omega) \qquad \text{for } (t, \omega) \in [0, T] \times \Omega_{\mathrm{par}}^n\,.$$

We endow $\Omega_{\mathrm{par}}^n$ with the filtration $\mathbb{G}^n := (\mathcal{G}_t^n)_{0\leq t\leq T}$ given by the right-continuous version of the canonical raw filtration, i.e.,

$$\mathcal{G}_t^n := \bigcap_{0<\varepsilon<T-t} \sigma\Big((X_s^n, Z_s^n)\,:\, 0 \leq s < t + \varepsilon\Big) \qquad \text{for } t \in [0, T),$$
$$\mathcal{G}^n := \mathcal{G}_T^n := \sigma\Big((X_s^n, Z_s^n)\,:\, 0 \leq s \leq T\Big)\,.$$

The following technical result is well-known and can be found e.g. in Crowell [12, Lem. 1.3].

**Lemma 1.3** | *Under Assumptions 1.1 and 1.2, there exist a unique probability measure* $\mathbb{P}^n_{x_0}$ *on* $(\Omega^n_{\text{par}}, \mathcal{G}^n)$ *and independent standard d-dimensional* $(\mathbb{G}^n, \mathbb{P}^n_{x_0})$*-Brownian motions* $B^{1,n}, \dots, B^{n,n}$ *such that* $Z^n$ *is a standard m-dimensional* $(\mathbb{G}^n, \mathbb{P}^n_{x_0})$*-Brownian motion that is independent of* $B^n := (B^{1,n}, \dots, B^{n,n})$*, and* $(X^n, B^n, Z^n)$ *solves* (1.1)*,* (1.2)*. Moreover, the law of* $X^n$ *under* $\mathbb{P}^n_{x_0}$ *is exchangeable and we have*

$$\mathbf{E}^{\mathbb{P}^n_{x_0}}\left[\sup_{t\in[0,T]} |X^{i,n}_t|^q\right] =: c_q < \infty \qquad \textit{for any } q > 1, \tag{1.3}$$

*where the constant* $c_q$ *depends only on* $q$*,* $T$ *and the bounds for* $b$*,* $\sigma$*,* $\bar{\sigma}$*, but not on* $n$*.*

By exchangeability, we mean that for any permutation $\pi$ of $[n] := \{1, \dots, n\}$, we have

$$\mathrm{Law}_{\mathbb{P}^n_{x_0}}(X^{1,n}, \dots, X^{n,n}) = \mathrm{Law}_{\mathbb{P}^n_{x_0}}(X^{\pi(1),n}, \dots, X^{\pi(n),n}).$$

In the sequel, we denote the weak solution from Lemma 1.3 by

$$\boldsymbol{X}^n = \Big(\Omega^n_{\text{par}}, \mathcal{G}^n, \mathbb{G}^n, \mathbb{P}^n_{x_0}, (X^n_t)_{t\in[0,T]}, (B^n_t)_{t\in[0,T]}, (Z^n_t)_{t\in[0,T]}\Big). \tag{1.4}$$

Note that Assumption 1.2(E) also allows $\bar{\sigma} = 0$, in which case there is no common noise $Z^n$ in (1.1) and the above simplifies in the obvious way.

## 1.2 Probabilistic setup for the measure flow

From the $n$-particle system (1.1) we have the *empirical measure flow* $\mu^n := (\mu^n_t)_{t\in[0,T]}$. The sequence $(\mu^n)_{n\in\mathbb{N}}$ of these processes is a key object in this work. To focus on it, we transfer it to a canonical setup. For this, let

$$\Omega := \mathbf{C}\Big([0,T];\, \mathcal{S}'(\mathbb{R}^d)\Big)$$

denote the *canonical trajectory space*. Its topological vector space structure is induced by convergence uniformly on $[0,T]$ for bounded sets in $\mathcal{S}(\mathbb{R}^d)$; see e.g. Kallianpur and Xiong [36, Ch. 2.4] or Appendix A. We denote the canonical process on $\Omega$ by

$$(t,\omega) \mapsto \omega_t =: \Lambda_t(\omega) = \Lambda_t \qquad \text{for } (t,\omega) \in [0,T] \times \Omega,$$

We endow $\Omega$ with the right-continuous version $\mathbb{F} := (\mathcal{F}_t)_{0\le t\le T}$ of the the canonical filtration generated by $\Lambda$, i.e.,

$$\mathcal{F}_t := \bigcap_{0<\varepsilon<T-t} \sigma(\Lambda_s \,:\, 0 \le s < t+\varepsilon) \qquad \text{for } t \in [0,T),$$
$$\mathcal{F} := \mathcal{F}_T := \sigma(\Lambda_s \,:\, 0 \le s \le T).$$

The set $\mathbf{M}^+_1(\mathbb{R}^d)$ is not closed in the topology of $\mathcal{S}'(\mathbb{R}^d)$. We therefore need to establish part 1) of the following technical result.

**Lemma 1.4** | *Under Assumptions 1.1 and 1.2, let $\mathbf{X}^n$ be the weak solution from Lemma 1.3 of the $n$-particle system* (1.1), (1.2)*. Then there exists a version of the process $\mu^n = (\mu^n_t)_{t\in[0,T]}$ which has continuous trajectories in the strong topology of $\mathcal{S}'(\mathbb{R}^d)$. In particular,*

$$\mathbf{P}^n_{x_0} := \mathbb{P}^n_{x_0} \circ (\mu^n)^{-1} \tag{1.5}$$

*is a well-defined probability measure on $(\Omega, \mathcal{F})$. Moreover:*

1) *The sequence $(\mathbf{P}^n_{x_0})_{n\in\mathbb{N}}$ is tight and any narrow cluster point $\mathbf{P}^\infty_{x_0}$ of $(\mathbf{P}^n_{x_0})_{n\in\mathbb{N}}$ satisfies*

$$\mathbf{P}^\infty_{x_0}\Big[\,\mathbf{C}\Big([0,T];\, \mathcal{S}'(\mathbb{R}^d) \cap \mathbf{M}^+_1(\mathbb{R}^d)\Big)\,\Big] = 1\,, \tag{1.6}$$

*that is, $\mathbf{P}^\infty_{x_0}$ concentrates its mass on the set of probability-measure-valued processes.*

2) *For any $\varepsilon, K > 0$ and $q > 1$, we have for all $n \in \mathbb{N}\cup\{\infty\}$ the concentration bound*

$$\mathbf{P}^n_{x_0}\Big[\omega \in \Omega : \Lambda_t\Big[[-K,K]^d\Big] < 1 - \varepsilon \text{ for some } t \in [0,T]\Big] \le \frac{c_q}{\varepsilon K^q}\,, \tag{1.7}$$

*where $c_q$ is the constant from* (1.3)*, which depends on $q$ but not on $K$, $\varepsilon$ nor $\mathbf{P}^n_{x_0}$.*

***Proof*** See [12, Lem. 1.4]. □

Lemma 1.4 lets us consider, for each $n \in \mathbb{N}$, the $\mathcal{S}'$-valued canonical process $\Lambda$ under $\mathbf{P}^n_{x_0}$. Note the distinction between $\mathbb{P}^n_{x_0}$ and $\mathbf{P}^n_{x_0}$: the former denotes the joint law of $n$ particles and the common Brownian motion, while the latter models the flow of the empirical measure of the particles. Accordingly, we term $\mathbf{P}^n_{x_0}$ the *empirical measure law*. Under this empirical measure law $\mathbf{P}^n_{x_0}$, defined via (1.5), the process $\Lambda$ represents the dynamics of the empirical measure of the $n$-particle system in the sense that we have $\mathrm{Law}_{\mathbb{P}^n_{x_0}}(\mu^n) = \mathrm{Law}_{\mathbf{P}^n_{x_0}}(\Lambda)$ for all $n \in \mathbb{N}$.

## 1.3 Main result

For each $n \in \mathbb{N}$, the dynamics of $\Lambda$ under $\mathbf{P}^n_{x_0}$ can be described explicitly. Indeed, in Section 2.3 below, we recall that $\Lambda$ under $\mathbf{P}^n_{x_0}$ is an $\mathcal{S}'$-valued semimartingale with continuous trajectories, which means that there is a unique decomposition $\Lambda = A + M$, where $A$ is an $\mathcal{S}'$-valued process of finite variation and $M$ is an $\mathcal{S}'$-valued local martingale, both with continuous trajectories in the strong topology of $\mathcal{S}'$; see Section 2.1 for a brief recap of these terms.

The processes $A$ and $M$ in the semimartingale decomposition can be described efficiently using a *martingale problem*. We develop the necessary terminology and notation for this in Sections 2.1 and 2.2, below. This description ultimately allows us to *represent the limiting dynamics* of $\Lambda$ under cluster points $\mathbf{P}^\infty_{x_0}$ of $(\mathbf{P}^n_{x_0})_{n\in\mathbb{N}}$.

To state our main result and discuss the strategy of proof, we preview the basic recurrent notation that we require. The data of the martingale problem is given in terms of maps defined for each $t \in [0,T]$ by

$$\mathscr{A}_t : \mathcal{S}' \to \mathcal{S}' \cup \{\infty\}\,, \tag{1.8}$$

with $\lambda \mapsto \mathscr{A}_t(\lambda)$, and

$$\mathscr{Q}_t : \mathcal{S}' \to \mathcal{L}^+_{\mathrm{bl}}(\mathcal{S}) , \tag{1.9}$$

$$\mathscr{C}_t : \mathcal{S}' \to \mathcal{L}^+_{\mathrm{bl}}(\mathcal{S}) , \tag{1.10}$$

with $\lambda \mapsto \mathscr{Q}_t(\lambda)$ and $\lambda \mapsto \mathscr{C}_t(\lambda)$, respectively, where $\mathcal{L}^+_{\mathrm{bl}}(\mathcal{S}')$ is a space of bilinear maps $\mathcal{S} \times \mathcal{S} \to \mathbb{R}$. In (2.26)–(2.28), below, we make $\mathscr{A}$, $\mathscr{Q}$, $\mathscr{C}$ explicit terms of $b$, $\sigma$, $\bar{\sigma}$.

As we recall in Proposition 2.6, the measure $\mathbf{P}^n_{x_0}$ solves for each $n \in \mathbb{N}$ the canonical $\mathcal{S}'$-valued martingale problem $\mathrm{MP}(\delta_{x_0}, \mathscr{A}, \mathscr{Q} + \frac{1}{n}\mathscr{C})$; see Definition 2.3 for what precisely this entails. As a consequence, the $\mathcal{S}'$-valued semimartingale decomposition of $\Lambda$ under $\mathbf{P}^n_{x_0}$ can be written explicitly as

$$\Lambda_t = \delta_{x_0} + \int_0^t \mathscr{A}_s(\Lambda_s)\,\mathrm{d}s + M_t\,, \qquad t \in [0,T]\,, \tag{1.11}$$

where $M$ is an $\mathcal{S}'$-valued square-integrable martingale with continuous trajectories, i.e., $M \in \mathscr{M}^c_2(\mathbb{F}, \mathbf{P}^n_{x_0}; \mathcal{S}')$, with (tensor) quadratic variation

$$\langle\!\langle M \rangle\!\rangle_t = \int_0^t \mathscr{Q}_s(\Lambda_s)\,\mathrm{d}s + \frac{1}{n}\int_0^t \mathscr{C}_s(\Lambda_s)\,\mathrm{d}s\,, \qquad t \in [0,T]\,. \tag{1.12}$$

By Lemma 1.4, we know that the sequence $(\mathbf{P}^n_{x_0})_{n\in\mathbb{N}}$ possesses at least one narrow cluster point $\mathbf{P}^\infty_{x_0}$. Moreover, $\Lambda$ under $\mathbf{P}^\infty_{x_0}$ remains $\mathbf{M}^+_1(\mathbb{R}^d)$-valued. If we strengthen Assumption 1.1 to the following, we are able to characterize the dynamics of $\Lambda$ under $\mathbf{P}^\infty_{x_0}$.

**Assumption 1.5** | The function $b$ is bounded and measurable relative to $\mathcal{E}$. In addition, there exists $r_0 > 1$ such that $b$ satisfies:

(C$_\mathrm{s}$) Whenever $\mu \in \mathbf{M}^+_1$ is absolutely continuous and $(\mu_k)_{k\in\mathbb{N}}$ is a sequence in $\mathbf{M}^+_1$ with $\mu_k \to \mu$ narrowly, then

$$\lim_{k\to\infty} \sup_{x\in K} |b_t(x,\mu_k) - b_t(x,\mu)| = 0$$

for any compact subset $K \subseteq \mathbb{R}^d$ and $t \in [0,T]$.

Our main result in Theorem 1.6, below, establishes that under Assumptions 1.2 and 1.5, any cluster point $\mathbf{P}^\infty_{x_0}$ solves the canonical $\mathcal{S}'$-valued martingale problem $\mathrm{MP}(\delta_{x_0}, \mathscr{A}, \mathscr{Q})$, the *formally limiting problem* of the sequence $(\mathrm{MP}(\delta_{x_0}, \mathscr{A}, \mathscr{Q} + \frac{1}{n}\mathscr{C}))_{n\in\mathbb{N}}$. As a consequence, we obtain the stability of the decomposition in (1.11) and thus of the $\mathcal{S}'$-valued semimartingale property of $\Lambda$ under narrow convergence.

**Theorem 1.6** | *Under Assumptions 1.2 and 1.5, any narrow cluster point $\mathbf{P}^\infty_{x_0}$ of the sequence $(\mathbf{P}^n_{x_0})_{n\in\mathbb{N}}$ solves the martingale problem $\mathrm{MP}(\delta_{x_0}, \mathscr{A}, \mathscr{Q})$. In particular, the coordinate process $\Lambda$ on $(\Omega, \mathcal{F}, \mathbf{P}^\infty_{x_0})$ possesses the $\mathcal{S}'$-valued semimartingale decomposition*

$$\Lambda_t = \delta_{x_0} + \int_0^t \mathscr{A}_s(\Lambda_s)\,\mathrm{d}s + M_t\,, \qquad t \in [0,T]\,, \tag{1.13}$$

*where $M = (M_t)_{t\in[0,T]}$ is in $\mathscr{M}^c_2(\mathbb{F}, \mathbf{P}^\infty_{x_0}; \mathcal{S}')$ with quadratic variation*

$$\langle\!\langle M \rangle\!\rangle_t = \int_0^t \mathscr{Q}_s(\Lambda_s)\,\mathrm{d}s\,, \qquad t \in [0,T]\,. \tag{1.14}$$

*Finally, the* empirical propagation of chaos property *holds, that is,*

$$\mathrm{Law}_{\mathbb{P}^{n_k}_{x_0}}(\mu^{n_k}) \longrightarrow \mathrm{Law}_{\mathbf{P}^{\infty}_{x_0}}(\Lambda) \text{ narrowly in } \mathbf{M}^+_1(\Omega) \tag{1.15}$$

*as* $k \to \infty$ *for any subsequence* $(n_k)_{k\in\mathbb{N}}$ *such that* $\mathbf{P}^{n_k}_{x_0} \to \mathbf{P}^{\infty}_{x_0}$.

Solutions of $\mathrm{MP}(\delta_{x_0}, \mathscr{A}, \mathscr{Q})$ are related to probabilistically and analytically weak solutions of the *stochastic Fokker–Planck equation*. This is a certain SPDE that, under appropriate assumptions, can be interpreted as modeling the limiting empirical measure flow of convergent subsequences of $(\mu^n)_{n\in\mathbb{N}}$. Weak solutions of the stochastic Fokker–Planck equation are in turn closely related to weak solutions of *McKean–Vlasov SDEs with common noise*. This relation is established using superposition theory; see Figalli [21] and Lacker et al. [44]. In this way, a weak solution theory for McKean–Vlasov SDEs can be obtained, and the empirical propagation of chaos (1.15) translates to a version of the propagation of chaos for weak solutions of McKean–Vlasov systems. We make these comments precise in forthcoming work that is based on Crowell [14].

With this in mind, our contribution in Theorem 1.6 is twofold. Firstly, the existence of solutions of $\mathrm{MP}(\delta_{x_0}, \mathscr{A}, \mathscr{Q})$ was obtained previously only for regular coefficients or in special cases; see Section 6 for a discussion and references. In particular, the stability property for the formally limiting martingale problem, i.e., the validity of (1.15), is by no means obvious under Assumptions 1.2 and 1.5. We make this point in the discussion of the strategy of proof below. Still, our assumptions are general enough to cover interesting examples appearing in the literature; see Section 5. Secondly, Theorem 1.6 is key to obtaining an existence result for $\mathrm{MP}(\delta_{x_0}, \mathscr{A}, \mathscr{Q})$ under more general assumptions than those in Theorem 1.6, yet at the expense of the convergence property (1.15); see Section 6.

### 1.4 Strategy of proof

The main difficulty in the proof of Theorem 1.6 is to establish the $(\mathbf{P}^{\infty}_{x_0}, \mathbb{F})$-martingale property of $M$. We now discuss why, and how we manage to obtain it.

**Preliminary thoughts — Steps 1 and 2** The starting point of our study is the measure $\mathbf{P}^n_{x_0}$ obtained in (1.5); this constitutes the first step. As the second step, we make the sequence $(\mathbf{P}^n_{x_0})_{n\in\mathbb{N}}$ akin to a systematic analysis that ultimately leads to Theorem 1.6. An efficient setting is provided by the framework of $\mathcal{S}'$-valued processes and $\mathcal{S}'$-valued martingale problems. We recall the basic definitions for this setting in Section 2.1. In Section 2.2, we derive using Itô's formula and the particle system (1.1), (1.2) the basic objects used to formulate the martingale problems we study, specifically, maps $\mathscr{A}$, $\mathscr{Q}$, and $\mathscr{C}$. Then, in Section 2.3, we show in Proposition 2.6 that for each $n \in \mathbb{N}$ the empirical measure flow law $\mathbf{P}^n_{x_0}$ solves the *martingale problem* $\mathrm{MP}(\delta_{x_0}, \mathscr{A}, \mathscr{Q} + \frac{1}{n}\mathscr{C})$. This implies in particular that the process $M$ in (1.11) is an $(\mathbf{P}^n_{x_0}, \mathbb{F})$-martingale, so that

$$\mathbf{E}^{\mathbf{P}^n_{x_0}}[(M_t - M_s)g] = 0 \qquad \text{in } \mathcal{S}' \tag{1.16}$$

for all $0 \le s \le t \le T$ and all continuous bounded $\mathcal{F}_s$-measurable functions $g : \Omega \to \mathbb{R}$. The expectation in (1.16) is understood as a Bochner integral. Standard linearity properties of this integral let us equivalently state (1.16) as $\mathbf{E}^{\mathbf{P}^n_{x_0}}[(M_t[\phi] - M_s[\phi])g] = 0$ for all $\phi \in \mathcal{S}$ and $0 \le s \le t \le T$, and $g$ as above. Our natural urge in (1.16) is to pass to narrow limit points $\mathbf{P}^{\infty}_{x_0}$ of $(\mathbf{P}^{n}_{x_0})_{n\in\mathbb{N}}$ and conclude that

$$\mathbf{E}^{\mathbf{P}^{\infty}_{x_0}}[(M_t - M_s)g] = 0 \qquad \text{in } \mathcal{S}' \tag{1.17}$$

for all $0 \le s \le t \le T$ and for $g$ as above. In fact, if we can establish (1.17), then by standard monotone class arguments, $M$ is a $(\mathbf{P}^{\infty}_{x_0}, \mathbb{F})$-local martingale. Establishing this martingale property is what we aim for, as it is needed to argue that $\mathbf{P}^{\infty}_{x_0}$ solves the martingale problem $\mathrm{MP}(\delta_{x_0}, \mathscr{A}, \mathscr{Q})$, which is the *formal limit* of the sequence $(\mathrm{MP}(\delta_{x_0}, \mathscr{A}, \mathscr{Q} + \frac{1}{n}\mathscr{C}))_{n\in\mathbb{N}}$.

**The difficulty** In passing (1.16) to the limit, we face an obstacle. To see this, use (1.11) under $\mathbf{P}^{n}_{x_0}$ to write

$$M_t - M_s = \Lambda_t - \Lambda_s - \int_s^t \mathscr{A}_u(\Lambda_u)\,\mathrm{d}u. \tag{1.18}$$

This identifies $M_t - M_s$ as a progressive functional on $(\Omega, \mathbb{F})$, i.e., as a measurable map of $\Lambda_{\cdot\wedge t}$ (or rather $\omega_{\cdot\wedge t}$). However, the function

$$\mathcal{S}' \cap \mathbf{M}^{+}_{1}(\mathbb{R}^d) \ni \lambda \mapsto \mathscr{A}_u(\lambda) \in \mathcal{S}',$$

is *not globally continuous* when $\mathbf{M}^{+}_{1}(\mathbb{R}^d)$ is equipped with either the narrow topology or the topology inherited from $\mathcal{S}'$. To see this, let us preview the definition of $\mathscr{A}_u$ given in (2.8) and (2.26). For each $\lambda \in \mathbf{M}^{+}_{1}$, $\mathscr{A}_u(\lambda)$ is given by

$$\mathcal{S} \ni \phi \mapsto \mathscr{A}_u(\lambda)[\phi] = \lambda[\mathscr{L}_u(\lambda)\phi] = \int_{\mathbb{R}^d} \Big(\mathscr{L}_t(\lambda)\phi\Big)(x)\,\lambda(\mathrm{d}x) \in \mathbb{R}, \tag{1.19}$$

and $\mathscr{A}_u(\lambda) = \infty$ whenever $\lambda \in \mathcal{S}' \setminus \mathbf{M}^{+}_{1}$, where

$$\Big(\mathscr{L}_u(\lambda)\phi\Big)(x) := b_u(x,\lambda)\cdot\nabla\phi(x) + \frac{1}{2}(\sigma_u\sigma_u^{\mathsf{T}} + \bar\sigma_u\bar\sigma_u^{\mathsf{T}})(x,\lambda) : \nabla^2\phi(x). \tag{1.20}$$

Now consider, for instance, $\lambda_m := \delta_{(\frac{1}{m},\dots,\frac{1}{m})}$ with $m \in \mathbb{N}$, which converges to $\delta_0$ in both the narrow topology and the strong topology of $\mathcal{S}'$, and let $b_t(x,\lambda) = B(x) := \mathbb{1}_{\{x_i \le 0, i=1,\dots,d\}}$, $\sigma = \bar\sigma \equiv 1$. Then the definition (1.19) gives $\mathscr{A}_u(\lambda_m) = B(1/m) = 2$, while by the same definition $\mathscr{A}_u(\delta_0) = B(0) = 3$; so $\lim_{m\to\infty}\mathscr{A}_u(\lambda_m) \neq \mathscr{A}_u(\delta_0)$. If $\mathscr{A}_u$ is not continuous, then by (1.18) neither is $M_t - M_s$ as a functional on $(\Omega, \mathbb{F})$. Consequently, passing in (1.16) to narrow cluster points to deduce (1.17) is generally not possible. As a matter of fact, continuity of the data of a martingale problem plays a central role in standard proofs of existence results; see e.g., Karatzas and Shreve [38, Proof of Thm. 5.4.22]. Without continuity, we need to devise a more careful argument.

Fortunately, the above obstacle is not insurmountable. To identify a viable strategy, we use the Skorokhod representation theorem for any subsequence $(\mathbf{P}^{n_k}_{x_0})_{k\in\mathbb{N}}$ with $\mathbf{P}^{n_k}_{x_0} \to \mathbf{P}^{\infty}_{x_0}$ narrowly as $k \to \infty$, to construct $\mathbf{C}([0,T];\mathcal{S}')$-valued random variables $\tilde\Lambda^k$ for $k \in \mathbb{N}\cup\{\infty\}$ on an abstract probability space $(\tilde\Omega, \tilde{\mathcal{F}}, \tilde{\mathbf{P}})$, satisfying that

(a) for all $k \in \mathbb{N}\cup\{\infty\}$, we have $\mathrm{Law}_{\mathbf{P}^{n_k}_{x_0}}(\Lambda) = \mathrm{Law}_{\tilde{\mathbf{P}}}(\tilde\Lambda^k)$;

(b) $\tilde{\mathbf{P}}$-a.s., $\tilde\Lambda^k \to \tilde\Lambda^\infty$ in $\mathbf{C}([0,T];\mathcal{S}')$ as $k \to \infty$.

With the random variables $(\tilde{\Lambda}^k)_{k\in\mathbb{N}\cup\{\infty\}}$, (1.18) may be recast as

$$\tilde{M}_t^k - \tilde{M}_s^k = \tilde{\Lambda}_t^k - \tilde{\Lambda}_s^k - \int_s^t \mathscr{A}_u(\tilde{\Lambda}_u^k)\,\mathrm{d}u \qquad \text{for } k \in \mathbb{N}\cup\{\infty\}\,, \tag{1.21}$$

and (1.16) as $\mathbf{E}^{\tilde{\mathbf{P}}}[(\tilde{M}_t^k - \tilde{M}_s^k)g] = 0$. By (1.21) and property (b), we see that in this setting, (1.17) follows from (1.16) once we establish that

$$\lim_{k\to\infty}\left|\mathbf{E}^{\tilde{\mathbf{P}}}\Big[\Big(\int_s^t \Big(\mathscr{A}_u(\tilde{\Lambda}_u^k) - \mathscr{A}_u(\tilde{\Lambda}_u^\infty)\Big)\,\mathrm{d}u\Big)g\Big]\right| = 0\,. \tag{1.22}$$

To prove (1.22), we show that while under Assumptions 1.2 and 1.5, and for $\phi \in \mathcal{S}$ the map $\mathcal{S}' \cap \mathbf{M}_1^+ \ni \lambda \mapsto \mathscr{A}_u(\lambda)[\phi] \in \mathbb{R}$ in (1.19) is not continuous in general, it is actually $(\mathrm{d}u \otimes \mathrm{d}\tilde{\mathbf{P}})$-a.s. *locally continuous along the sequence* $(\tilde{\Lambda}^k)_{k\in\mathbb{N}\cup\{\infty\}}$, in the sense that $(\mathrm{d}u \otimes \mathrm{d}\tilde{\mathbf{P}})$-a.s. we have $\mathscr{A}_u(\tilde{\Lambda}_u^k)[\phi] \longrightarrow \mathscr{A}_u(\tilde{\Lambda}_u^\infty)[\phi]$ as $k\to\infty$, and standard linearity properties of the Bochner integral allow us to deduce (1.22).

**The role of regularity — Steps 3 and 4** To understand the origin of the local continuity of the map in (1.19), we use that $\mathscr{A}_u(\lambda)[\phi] = \lambda[\mathscr{L}_u(\lambda)\phi]$ takes a specific form, namely that *the measure* $\lambda$ *integrates the function* $(\mathscr{L}_u(\lambda)\phi)(x)$. If $\lambda$ is a measure with a density possessing favorable analytic regularity, this feature allows us to compensate for a lack of regularity in the function $\mathscr{L}_t(\lambda)\phi$, and thus via the definition in (1.20), for a lack of regularity in the drift $b$. This observation is key. To exploit it fruitfully in the context of the random convergent sequence $\tilde{\Lambda}^k \to \tilde{\Lambda}^\infty$, we establish in the third step towards the proof of Theorem 1.6 the necessary regularity estimates; see Section 3. These regularity estimates come in two forms: an estimate for the limit, i.e., for $k=\infty$, and a uniform estimate for the pre-limit, i.e., for $k\in\mathbb{N}$.

For the limit $k=\infty$, the key ingredient is an *emergence of regularity* property that we established in Crowell [12] and recall in Section 3.1, below. The emergence of regularity shows for the limit $\tilde{\Lambda}^\infty$ that for each $t\in(0,T]$, we have $\tilde{\mathbf{P}}$-a.s. that $\tilde{\Lambda}_t^\infty \ll \mathrm{d}x$. In fact, for any $q\geq 1$, we can find $r>1$ and $0<\gamma<1$ such that the function $t\mapsto(\tilde{\omega}\mapsto \mathrm{d}\tilde{\Lambda}_t^\infty(\tilde{\omega})/\mathrm{d}x)$ defines a measurable map

$$\tilde{p} : (0,T] \to \mathbb{L}^q\Big((\tilde{\Omega},\tilde{\mathcal{F}},\tilde{\mathbf{P}});\, \mathbb{L}^r(\mathbb{R}^d)\Big)$$

which satisfies for some constant $c<\infty$ the bound

$$\Big\|\|\tilde{p}(t)\|_{\mathbb{L}^r(\mathbb{R}^d)}\Big\|_{\mathbb{L}^q(\tilde{\mathbf{P}})} \leq c(1\wedge t)^{-\gamma} \quad \text{for all } t\in(0,T]. \tag{1.23}$$

The quantitative estimate in (1.23) is the type of favorable analytic regularity of the density that we use to exploit the above observation.

The emergence of regularity is a limiting phenomenon. For $k\in\mathbb{N}$, the measure $\tilde{\Lambda}_t^k$ is $\tilde{\mathbf{P}}$-a.s. purely atomic for all $t\in[0,T]$ and thus does not possess comparable regularity properties. However, for $t\in(0,T]$, the time-marginal laws of each particle $X^{i,n_k}$, i.e., $\mathrm{Law}_{\mathbb{P}_{x_0}^{n_k}}(X_t^{i,n_k})$, satisfy a regularity estimate of the same type as in (1.23), and this uniformly for $n\in\mathbb{N}$. We establish this in Section 3.2, below.

In the fourth and final step, we deduce the property in (1.22) and prove Theorem 1.6. This uses a careful approximation argument for the drift $b$ that is facilitated by the regularity estimates from step 3. The estimate in (1.23) obtained from the emergence of regularity can be used directly, exploiting the uniform regularity estimate for the time-marginals $\mathrm{Law}_{\mathbb{P}^{n_k}_{x_0}}(X^{i,n_k}_t)$ requires an intricate compactness argument. This uses the $n_k$-particle systems and their exchangeability. We explain how this is done in detail before carrying out the relevant estimates in Lemma 4.6, below.

## 2 EMPIRICAL MEASURE FLOW MARTINGALE PROBLEM

This section serves a preparatory purpose. In Section 2.1, we recall some of the basic terminology from $\mathcal{S}'$-valued stochastic analysis. In particular, the notion of an *canonical $\mathcal{S}'$-valued martingale problem*, which we recall in Definition 2.3 below, is used repeatedly in this paper. In Section 2.2, starting with the empirical measure flow $\mu^n$ from (1.2), we develop the basic notation used to define for each $n \in \mathbb{N}$ the canonical $\mathcal{S}'$-valued martingale problem $\mathrm{MP}(\delta_{x_0}, \mathscr{A}, \mathscr{Q} + \frac{1}{n}\mathscr{C})$. Finally, in Section 2.3, we show that for each $n \in \mathbb{N}$, the empirical measure flow law $\mathbf{P}^n_{x_0}$ from (1.5) solves $\mathrm{MP}(\delta_{x_0}, \mathscr{A}, \mathscr{Q} + \frac{1}{n}\mathscr{C})$.

### 2.1 Background on $\mathcal{S}'$-valued semimartingales and martingale problems

Martingale problems form a backbone of the weak solution theory for $\mathbb{R}^d$-valued SDEs. This formulation is well known and goes back to a series of works by Stroock and Varadhan [61, 62]. We refer to the classic treatment in the book by the same authors [63, Ch. 6], or the more modern account in Karatzas and Shreve [38, Ch. 5.3].

In this work, we require the martingale problem formulation for processes valued in $\mathcal{S}'$ which we recall below; see e.g. Mikulevicius and Rozovskii [9, Ch. 6.3] for more details. To efficiently develop this theory, we recall first some of the basic structural properties of $\mathcal{S}'$ and basic notions of $\mathcal{S}'$-valued semimartingales. A foundational treatment of $\mathcal{S}'$-valued stochastic calculus is found in Itô [35]. More modern expositions are available in e.g. Kallianpur and Xiong [36] or Dalang and Sanz-Solé [15].

**Structure of $\mathcal{S}'$** Let $(\,\cdot\,;\,\cdot\,)_0$ be the usual inner product on $\mathbb{L}^2(\mathrm{d}x)$. The map $f \mapsto (\,\cdot\,; f)_0$ defines a continuous embedding $\mathcal{S} \hookrightarrow \mathcal{S}'$, which the allows us to identify $\mathcal{S}$ and $\mathbb{L}^2(\mathrm{d}x)$ with subspaces of $\mathcal{S}'$. The Fréchet topology of $\mathcal{S}$ and the strong topology of $\mathcal{S}'$ can be described by a family of Hilbertian seminorms, or equivalently, in terms of a family of Hilbert spaces; see e.g. Simon [58] or Appendix A. In concrete terms, we can consider the family $(\mathscr{H}_p)_{p\in\mathbb{R}}$ of *Fourier–Hermite spaces*, defined as follows. Let $(h_k)_{k\in\mathbb{N}^d}$ be the Hermite functions on $\mathbb{R}^d$. Note that $h_k \in \mathcal{S}$ for each $k \in \mathbb{N}^d$. Hence for $\lambda \in \mathcal{S}'$, the $k$'th Hermite–Fourier coefficient $\lambda^{\#}_k := \lambda[h_k]$ is well defined. For $k \in \mathbb{N}$, let $\langle k\rangle := (1+|k|^2)^{1/2}$. For each $\lambda \in \mathcal{S}'$, there exists $p \in \mathbb{R}$ such that

$$\|\lambda\|_{\mathscr{H}_p} := \left(\sum_{k\in\mathbb{N}} (\langle k\rangle^{p/d}\lambda^{\#}_k)^2\right)^{1/2} < \infty\,.$$

In fact, if $\lambda \in \mathcal{S} \subsetneq \mathcal{S}'$, then the above sum is finite for any $p \in \mathbb{R}$. With this in mind, we can define a Hilbert-subspace of $\mathcal{S}'$ by considering

$$\mathscr{H}_p := \{\lambda \in \mathcal{S}' : \|\lambda\|_{\mathscr{H}_p} < \infty\}$$

with the norm $\|\cdot\|_{\mathscr{H}_p}$. For $p \geq 0$, $\mathscr{H}_p$ can be identified with a subset of $\mathbb{L}^2(\mathrm{d}x)$ and thus with a set of functions. Finally, the $\mathbb{L}^2$ duality pairing $f \mapsto (f;\,\cdot\,)_0$ extends for any $p \in \mathbb{R}$ to an isometric isomorphism $J : \mathscr{H}_{-p} \to \mathscr{H}_p'$ so that $\mathscr{H}_{-p} \simeq \mathscr{H}_p'$. In sum, the above shows that $\mathcal{S} \hookrightarrow \mathscr{H}_p \hookrightarrow \mathbb{L}^2 \hookrightarrow \mathscr{H}_{-p} \hookrightarrow \mathcal{S}'$ for all $p > 0$. It is a classical fact that $\mathcal{S} = \bigcap_{p\in\mathbb{R}} \mathscr{H}_p$ and $\mathcal{S}' = \bigcup_{p\in\mathbb{R}} \mathscr{H}_p$, where the right-hand sides are respectively endowed with the limit and colimit tolology. This makes precise the statement that the topology of $\mathcal{S}$ and $\mathcal{S}'$ can be described by a family of Hilbert spaces.

**Semimartingales valued in** $\mathcal{S}'$ An $\mathcal{S}'$-valued process $M = (M_t)_{t\in[0,T]}$ is called a *square-integrable, continuous $\mathcal{S}'$-valued martingale* relative to $(\mathbb{F}, \mathbf{P})$ if for each $\phi \in \mathcal{S}$, there is a version of $M[\phi] = (M_t[\phi])_{t\in[0,T]}$ in $\mathscr{M}_2^{\mathrm{c}}(\mathbb{F}, \mathbf{P}; \mathbb{R})$; we write $M \in \mathscr{M}_2^{\mathrm{c}}(\mathbb{F}, \mathbf{P}; \mathcal{S}')$. The process $M$ is called a *continuous $\mathcal{S}'$-valued local martingale* if for each $\phi \in \mathcal{S}$, there is a version of $M[\phi]$ in $\mathscr{M}_{\mathrm{loc}}^{\mathrm{c}}(\mathbb{F}, \mathbf{P}; \mathbb{R})$; we write $M \in \mathscr{M}_{\mathrm{loc}}^{\mathrm{c}}(\mathbb{F}, \mathbf{P}; \mathcal{S}')$.

We call $M = (M_t)_{t\in[0,T]}$ a *square-integrable, continuous $\mathscr{H}_{-p}$-valued martingale* if it is $\mathscr{H}_{-p}$-valued and a martingale with continuous trajectories relative to the $\mathscr{H}_{-p}$-norm, $\mathbf{P}$-a.s., and such that

$$\|M\|_{\mathscr{M}_2^{\mathrm{c}}(\mathbb{F},\mathbf{P};\mathscr{H}_{-p})} \coloneqq \mathbf{E}[\|M_T\|^2_{\mathscr{H}_{-p}}]^{\frac{1}{2}} < \infty\,. \tag{2.1}$$

In this case, we write $M \in \mathscr{M}_2^{\mathrm{c}}(\mathbb{F}, \mathbf{P}; \mathscr{H}_{-p})$. We call $M$ a *continuous $H$-valued local martingale* if there exists a localizing sequence $(\tau^n)_{n\in\mathbb{N}}$ such that $M_{\cdot\wedge\tau^n} \in \mathscr{M}_2^{\mathrm{c}}(\mathbb{F}, \mathbf{P}; \mathscr{H}_{-p})$. In this case we write $M \in \mathscr{M}_{\mathrm{loc}}^{\mathrm{c}}(\mathbb{F}, \mathbf{P}, \mathscr{H}_{-p})$.

We should emphasize the asymmetry in the requirements placed on the localizing sequence. In the $\mathcal{S}'$-valued case, this sequence may depend on $\phi$, in the $\mathscr{H}_{-p}$-valued case, it must not.

**Definition 2.1** | Let $\Lambda$ be a continuous $\mathcal{S}'$-valued process. We call $\Lambda$ a *continuous $\mathcal{S}'$-valued $(\mathbb{F}, \mathbf{P})$-semimartingale* if there are a random variable $\Lambda_0$ and processes $A = (A_t)_{t\in[0,T]}$, $M = (M_t)_{t\in[0,T]}$, all $\mathcal{S}'$-valued, such that

$$\Lambda_t = \Lambda_0 + A_t + M_t \quad \text{in } \mathcal{S}' \text{ for all } t \in [0, T]$$

and the following conditions are satisfied:

*1)* $\Lambda_0$ is $\mathcal{F}_0$-measurable.

*2)* $A$ is $\mathbb{F}$-adapted, $A_0 = 0$ and for all $\phi \in \mathcal{S}$, the process given by $A[\phi] = (A_t[\phi])_{t\in[0,T]}$ has continuous paths of finite variation.

*3)* $M \in \mathscr{M}_{\mathrm{loc}}^{\mathrm{c}}(\mathbb{F}, \mathbf{P}; \mathcal{S}')$ and $M_0 = 0$. Equivalently, $M[\phi] = (M_t[\phi])_{t\in[0,T]} \in \mathscr{M}_{\mathrm{loc}}^{\mathrm{c}}(\mathbb{F}, \mathbf{P}; \mathbb{R})$ for all $\phi \in \mathcal{S}$ and $M_0 = 0$.

Let $\Lambda$ be an $\mathscr{H}_{-p}$-valued process with continuous trajectories, $\mathbf{P}$-a.s. We call $\Lambda$ a *continuous $\mathscr{H}_{-p}$-valued $(\mathbb{F}, \mathbf{P})$-semimartingale* if there are an $\mathscr{H}_{-p}$-valued random variable $\Lambda_0$ and $\mathscr{H}_{-p}$-valued processes $A$ and $M$ such that

$$\Lambda_t = \Lambda_0 + A_t + M_t \quad \text{in } \mathscr{H}_{-p} \text{ for all } t \in [0, T]\,,$$

and the following conditions are satisfied: (a) $\Lambda_0$ is $\mathcal{F}_0$-measurable, (b) $A$ is $\mathbb{F}$-adapted, $A_0 = 0$ and $A$ has continuous paths of finite variation relative to the $\mathscr{H}_{-p}$-norm, (c) $M \in \mathscr{M}_{\mathrm{loc}}^{\mathrm{c}}(\mathbb{F}, \mathbf{P}; \mathscr{H}_{-p})$ and $M_0 = 0$.

**Lemma 2.2** | *Let $\tilde\Lambda = (\tilde\Lambda_t)_{t\in[0,T]}$ be a cylindrical process on $\mathcal{S}$, i.e. $\tilde\Lambda_t : \mathcal{S} \to \mathbb{L}^0(\Omega, \mathcal{F}, \mathbf{P})$ is a linear map for each $t \in [0, T]$. Suppose that $\tilde\Lambda_t : \mathcal{S} \to \mathbb{L}^0(\Omega, \mathcal{F}, \mathbb{P})$ is continuous for*

*each $t \in [0,T]$ and that for each $\phi \in \mathcal{S}$, the process $\tilde{\Lambda}[\phi] = (\tilde{\Lambda}_t[\phi])_{t\in[0,T]}$ is a real-valued $(\mathbb{F}, \mathbf{P})$-semimartingale with continuous trajectories and canonical decomposition*

$$\tilde{\Lambda}[\phi] = \Lambda_0^\phi + A^\phi + M^\phi\,,$$

*where $\Lambda_0^\phi$ is an $\mathcal{F}_0$-measurable real-valued random variable, $A^\phi$ is a real-valued $\mathbb{F}$-adapted and continuous process of finite variation, and $M^\phi \in \mathscr{M}^{\mathrm{c}}_{\mathrm{loc}}(\mathbb{F}, \mathbf{P}; \mathbb{R})$. Then:*

1) *There is a continuous $\mathcal{S}'$-valued $(\mathbb{F}, \mathbf{P})$-semimartingale $\Lambda$ satisfying $\Lambda_0[\phi] = \Lambda_0^\phi$, $A[\phi] = A^\phi$, $M[\phi] = M^\phi$ for all $\phi \in \mathcal{S}$, $\mathbf{P}$-a.s.*

2) *If for each $\phi \in \mathcal{S}$, the real-valued semimartingale $\Lambda[\phi] = \Lambda_0[\phi] + A[\phi] + M[\phi]$ satisfies in addition that*

$$\Big\| \int_0^T |\mathrm{d}A_t[\phi]| \Big\|_{\mathbb{L}^2(\Omega,\mathbb{F},\mathbf{P})} + \|M[\phi]\|_{\mathscr{M}_2^{\mathrm{c}}(\mathbb{F},\mathbf{P};\mathbb{R})} < \infty\,,$$

*then there exists $p_0 > 0$ such that for all $p \geq p_0$, we have $\mathbf{P}$-a.s. that $\Lambda, A, M \in \mathbf{C}([0,T];\ \mathscr{H}_{-p})$ and $M \in \mathscr{M}_2^{\mathrm{c}}(\mathbb{F}, \mathbf{P}; \mathscr{H}_{-p})$.*

***Proof*** By [46, Thm. 1],we find an $\mathcal{S}'$-valued continuous version $\Lambda$ of $\tilde{\Lambda}$, and part (i) of the lemma then follows from Pérez-Abreu [49, Thm. 1]. Part (ii) follows from e.g. Fonseca-Mora [22, Cor. 3.8]. □

**Tensor quadratic variation** We call a map $Q : \mathcal{S} \times \mathcal{S} \to \mathbb{R}$ *symmetric* if $Q[\phi, \psi] = Q[\psi, \phi]$ for all $\phi, \psi \in \mathcal{S}$, and *nonnegative* if $Q[\phi, \phi] \geq 0$ for all $\phi \in \mathcal{S}$. Denote by $\mathcal{L}^+_{\mathrm{bl}}(\mathcal{S})$ the space of bilinear maps $\mathcal{S} \times \mathcal{S} \to \mathbb{R}$ which are symmetric and nonnegative. This space carries the topology of uniform convergence on products of bounded sets in $\mathcal{S}$; see Huang and Yan [32, Sec. 3.3] for details.

A process $Q : [0,T] \times \Omega \to \mathcal{L}^+_{\mathrm{bl}}(\mathcal{S})$ is called *weakly $\mathbb{F}$-progressive* if for all $\phi, \psi \in \mathcal{S}$, the $\mathbb{R}$-valued process $Q[\phi, \psi]$ is $\mathbb{F}$-progressive, and *strongly $\mathbb{F}$-progressive* if the process $Q$ is $\mathbb{F}$-progressive when $\mathcal{L}^+_{\mathrm{bl}}(\mathcal{S}')$ is equipped with its Borel-$\sigma$-algebra. This is analogous to the weak and strong measurability in known from infinite-dimensional analysis because an element of $\mathcal{L}^+_{\mathrm{bl}}(\mathcal{S})$ can equivalently be seen as a continuous linear map $\mathcal{S} \to \mathcal{S}'$; see Huang and Yan [32, Thm. 3.17] or Simon [58, Sec. 4].

We say that $M \in \mathscr{M}_2^{\mathrm{c}}(\mathbb{F}, \mathbf{P}; \mathcal{S}')$ has the *(tensor) quadratic variation* process $\int_0^t Q_s\,\mathrm{d}s$ with *(tensor) covariance operator* $Q : [0,T] \times \Omega \to \mathcal{L}^+_{\mathrm{bl}}(\mathcal{S})$ if $Q$ is weakly $\mathbb{F}$-progressive and

$$\Big( M_t[\phi]M_t[\psi] - \int_0^t Q_s[\phi,\psi]\,\mathrm{d}s \Big)_{t\in[0,T]} \tag{2.2}$$

is in $\mathscr{M}^{\mathrm{c}}(\mathbb{F}, \mathbf{P}; \mathbb{R})$ for all $\phi, \psi \in \mathcal{S}$. It is well known that if $M \in \mathscr{M}_2^{\mathrm{c}}(\mathbb{F}, \mathbf{P}; \mathcal{S}')$, then the tensor quadratic variation exists and is uniquely characterized by (2.2), and we have

$$\Big( \int_0^t Q_s\,\mathrm{d}s \Big)[\phi,\psi] = \int_0^t Q_s[\phi,\psi]\,\mathrm{d}s \tag{2.3}$$

for all $t \in [0,T]$ and $\phi, \psi \in \mathcal{S}$, $\mathbf{P}$-a.s. Moreover, $Q$ can be taken to be strongly $\mathbb{F}$-progressive so that the integral $\int Q_t\,\mathrm{d}t$ above is understood in the Bochner sense. Indeed, Lemma 2.2, we may assume that there is a version (still denoted by $M$) in $\mathscr{M}_2^{\mathrm{c}}(\mathscr{H}_{-p})$. But in the setting of Hilbert spaces, the result is classical, see Gawarecki and Mandrekar [24, Lem. 2.1], and $(Q_t)_{t\in[0,T]}$ is Bochner-integrable.

We generally denote the *tensor variation process* by $\langle\langle M\rangle\rangle$, that is, $\langle\langle M\rangle\rangle$ is the unique $\mathcal{L}^+_{\mathrm{bl}}(\mathcal{S})$-valued, progressive process such that

$$M[\phi]M[\psi] - \langle\langle M\rangle\rangle[\phi,\psi] \in \mathscr{M}^{\mathrm{c}}(\mathbb{R}) \qquad \text{for all } \phi,\psi \in \mathcal{S}\,.$$

**Canonical $\mathcal{S}'$-valued martingale problems** We take as pre-specified the measurable space $(\Omega,\mathcal{F})$ together with the right-continuous version of the canonical filtration $\mathbb{F} = (\mathcal{F}_t)_{t\in[0,T]}$ and the $\mathcal{S}'$-valued canonical process $\Lambda = (\Lambda_t)_{t\in[0,T]}$.

Assume we are then given a triplet $(\nu_0, A, Q)$, consisting of

– a probability measure $\nu_0$ on $\mathcal{B}(\mathbb{R}^d)$;

– an $\mathbb{F}$-progressive process $A : [0,T]\times\Omega \to \mathcal{S}'\cup\{+\infty\}$, i.e. $A[\phi]$ is $\mathbb{F}$-progressive for each $\phi\in\mathcal{S}$, and $A_0 = 0$;

– an $\mathbb{F}$-progressive process $Q : [0,T]\times\Omega \to \mathcal{L}^+_{\mathrm{bl}}(\mathcal{S})$, i.e. $Q[\phi,\psi]$ is $\mathbb{F}$-progressive for each $\phi,\psi\in\mathcal{S}$.

We say that $(\nu_0, A, Q)$ defines the *data of a martingale problem*, and write $\mathrm{MP}(\nu_0, A, Q)$.

**Definition 2.3** | A probability measure $\mathbf{P}$ on $(\Omega,\mathcal{F})$ is said to be a *solution of* $\mathrm{MP}(\nu_0, A, Q)$ if the following hold:

*1)* We have that $\mathrm{Law}_{\mathbf{P}}(\Lambda_0) = \nu_0$.

*2)* For each $\phi\in\mathcal{S}$, we have that $\mathbf{P}[\int_0^T |A_s[\phi]|\,\mathrm{d}s < \infty]$.

*3)* If $M = (M_t)_{t\in[0,T]}$ is defined by

$$M_t \coloneqq \Lambda_t - \nu_0 - \int_0^t A_s\,\mathrm{d}s \quad \text{for } t\in[0,T],$$

then $M \in \mathscr{M}^{\mathrm{c}}_{\mathrm{loc}}(\mathbb{F},\mathbf{P};\mathcal{S}')$.

*4)* $\langle M[\phi], M[\psi]\rangle = \int Q_t[\phi,\psi]\,\mathrm{d}t$ $\mathbf{P}$-a.s. whenever $\phi,\psi\in\mathcal{S}$.

The reason we allow $A$ to take the value $+\infty$ is because this lets us construct martingale problems whose solutions, provided they exist, are prohibited of putting mass on certain measurable subsets of $\mathcal{S}'$ via condition 1) above; see Mikulevicius and Rozovskii [9, Def. 6.3.3] for an alternative but equivalent formulation.

**Remark 2.4** | If $\mathbf{P}$ solves $\mathrm{MP}(\nu_0, A, Q)$ and if $M \in \mathscr{M}^{\mathrm{c}}_2(\mathbb{F},\mathbf{P};\mathcal{S}')$, then part *3)* in Definition 2.3 becomes

*4')* $\langle\langle M\rangle\rangle = \int Q_t\,\mathrm{d}t$ $\mathbf{P}$-a.s.

**Remark 2.5** | There is an equivalent way to formulate martingale problems in terms of a second order differential operator acting on Fréchet differentiable functions $\mathcal{S}'\to\mathbb{R}$. We refer to e.g. Kallianpur et al. [37] for this approach. We avoid it here because it would require additional notation.

## 2.2 Data for the martingale problem

Consider again the particle system (1.1) and (1.2). A central object in this section is the *empirical measure flow* from (1.2), i.e.

$$\mu_t^n = \frac{1}{n}\sum_{i=1}^{n} \delta_{X_t^{i,n}} \quad \text{with } t \in [0,T]. \tag{2.4}$$

By Lemma 1.4, the canonical process $\Lambda$ on $(\Omega, \mathcal{F}, \mathbf{P}_{x_0}^n)$ represents the dynamics of the empirical measure in the sense that for all $n \in \mathbb{N}$, $\mathrm{Law}_{\mathbb{P}_{x_0}^n}(\mu^n) = \mathrm{Law}_{\mathbf{P}_{x_0}^n}(\Lambda)$. As a steppingstone to describe the dynamics of $\Lambda$ under $\mathbf{P}_{x_0}^n$, we now develop the basic notation based on the dynamics of $n$-particle systems. This notation is used repeatedly. We introduce the maps

$$\mathscr{L}_t : \mathbf{M}_1^+(\mathbb{R}^d) \times \mathcal{S}(\mathbb{R}^d) \to \mathbb{L}^1(\mathbb{R}^d;\, \mathbb{R}) \cap \mathbb{L}^\infty(\mathbb{R}^d; \mathbb{R})\,, \tag{2.5}$$

$$\mathscr{R}_t : \mathbf{M}_1^+(\mathbb{R}^d) \times \mathcal{S}(\mathbb{R}^d) \to \mathbb{L}^1(\mathbb{R}^d;\, \mathbb{R}^d) \cap \mathbf{C}_b(\mathbb{R}^d; \mathbb{R}^d)\,, \tag{2.6}$$

$$\mathscr{U}_t : \mathbf{M}_1^+(\mathbb{R}^d) \times \mathcal{S}(\mathbb{R}^d) \to \mathbb{L}^1(\mathbb{R}^d;\, \mathbb{R}^d) \cap \mathbf{C}_b(\mathbb{R}^d; \mathbb{R}^d)\,, \tag{2.7}$$

defined as follows. We set

$$(\lambda, \phi) \mapsto \mathscr{L}_t(\lambda)\phi\,, \qquad (\lambda, \phi) \mapsto \mathscr{R}_t(\lambda)\phi\,, \qquad (\lambda, \phi) \mapsto \mathscr{U}_t(\lambda)\phi\,,$$

where

$$\Big(\mathscr{L}_t(\lambda)\phi\Big)(x) := b_t(x,\lambda) \cdot \nabla\phi(x) + \frac{1}{2} a_t(x,\lambda) : \nabla^2\phi(x)\,, \tag{2.8}$$

$$\Big(\mathscr{R}_t(\lambda)\phi\Big)(x) := \bar{\sigma}_t^{\mathsf{T}}(x,\lambda)\nabla\phi(x)\,, \tag{2.9}$$

$$\Big(\mathscr{U}_t(\lambda)\phi\Big)(x) := \sigma_t^{\mathsf{T}}(x,\lambda)\nabla\phi(x)\,. \tag{2.10}$$

The notation specifying the function (2.8) is explicitly given by

$$b_t(x,\lambda) \cdot \nabla\phi(x) = \sum_{j=1}^{d} b_t^j(x,\lambda)\partial_j\phi(x)\,,$$

$$a_t(x,\lambda) : \nabla^2\phi(x) = \sum_{i,j=1}^{d} a_t^{ij}(x,\lambda)\partial_{ij}\phi(x)\,,$$

with

$$a_t^{ij}(x,\lambda) \coloneqq (\sigma_t\sigma_t^{\mathsf{T}} + \bar{\sigma}_t\bar{\sigma}_t^{\mathsf{T}})^{ij}(x,\lambda) = \sum_{k=1}^{d}(\sigma_t^{ik}\sigma_t^{kj} + \bar{\sigma}_t^{ik}\bar{\sigma}_t^{kj})(x,\lambda)\,,$$

the notation specifying (2.9) by

$$\bar{\sigma}_t^{\mathsf{T}}(x,\lambda)\nabla\phi(x) = \sum_{j=1}^{d} \bar{\sigma}_t^{j\cdot}(x,\lambda)\partial_j\phi(x)\,,$$

and the notation specifying (2.10) by

$$\sigma_t^{\mathsf{T}}(x,\lambda)\nabla\phi(x) = \sum_{j=1}^{d} \sigma_t^{j\cdot}(x,\lambda)\partial_j\phi(x)\,.$$

The fact that for each $(t, \lambda, \phi) \in [0, T] \times \mathbf{M}_1^+ \times \mathcal{S}$, the function $x \mapsto (\mathscr{L}_t(\lambda)\phi)(x)$ is an element of $\mathbb{L}^1(\mathbb{R}^d; \mathbb{R}) \cap \mathbb{L}^\infty(\mathbb{R}^d; \mathbb{R})$ follows from (2.8) together with observing that $b_t(\,\cdot\,, \lambda)$, $\sigma_t(\,\cdot\,, \lambda)$ and $\bar{\sigma}_t(\,\cdot\,, \lambda)$ are bounded, and that $\partial_i\phi$ and $\partial_{ij}\phi$ are for any $i, j \in [n]$ in $\mathcal{S}$, and thus also elements of $\mathbb{L}^s(\mathbb{R}^d; \mathbb{R})$ for any $s \in [1, \infty]$. Similarly, the fact that $x \mapsto (\mathscr{R}_t(\lambda)\phi)(x)$ and $x \mapsto (\mathscr{U}_t(\lambda)\phi)(x)$ are in $\mathbb{L}^1(\mathbb{R}^d; \mathbb{R}^d) \cap \mathbf{C}_b(\mathbb{R}^d; \mathbb{R}^d)$ follows from (2.9) and (2.10) along with Assumption 1.2 by which $x \mapsto \sigma_t(x, \lambda)$ and $x \mapsto \bar{\sigma}_t(x, \lambda)$ are continuous and bounded, together with the fact that $\partial_i\phi$ is for any $i \in [n]$ in $\mathcal{S}$, and thus also in $\mathbf{C}_b(\mathbb{R}^d; \mathbb{R})$ and in $\mathbb{L}^1(\mathbb{R}^d; \mathbb{R})$.

We often find it convenient to simplify the notation and write

$$\mathscr{L}_t(\lambda)\phi(x) := \Big(\mathscr{L}_t(\lambda)\phi\Big)(x)\,, \tag{2.11}$$

$$\mathscr{R}_t(\lambda)\phi(x) := \Big(\mathscr{R}_t(\lambda)\phi\Big)(x)\,, \tag{2.12}$$

$$\mathscr{U}_t(\lambda)\phi(x) := \Big(\mathscr{U}_t(\lambda)\phi\Big)(x)\,. \tag{2.13}$$

This notation emphasizes that for any $(t, \lambda) \in [0, T] \times \mathbf{M}_1^+$, $\mathscr{L}_t(\lambda)$, $\mathscr{R}_t(\lambda)$ and $\mathscr{U}_t(\lambda)$ are *differential operators* acting on functions in $\mathcal{S}$.

We point out that the maps $\mathscr{L}_t$, $\mathscr{R}_t$ and $\mathscr{U}_t$ arise naturally via Itô's formula. Indeed, for each $n \in \mathbb{N}$ and $i \in [n]$, consider the weak solution $\boldsymbol{X}^n$ of (1.1), (1.2) from Lemma 1.3. If $\phi \in \mathcal{S}$, Itô's formula implies that

$$\phi(X_t^{i,n}) = \phi(x_0) + \int_0^t \mathscr{L}_s(\mu_s^n)\phi(X_s^{i,n})\,\mathrm{d}s + (M^\phi)_t^{i,n} + (\bar{M}^\phi)_t^{i,n} \tag{2.14}$$

for $i \in [n], t \in [0, T]$, where $(M^\phi)^{i,n}$ and $(\bar{M}^\phi)^{i,n}$ are the real-valued $(\mathbb{G}^n, \mathbb{P}_{x_0}^n)$-local martingales

$$(M^\phi)_t^{i,n} := \int_0^t \mathscr{U}_s(\mu_s^n)\phi(X_s^{i,n})\,\mathrm{d}B_s^{i,n}\,, \tag{2.15}$$

$$(\bar{M}^\phi)_t^{i,n} := \int_0^t \mathscr{R}_s(\mu_s^n)\phi(X_s^{i,n})\,\mathrm{d}Z_s^n\,, \tag{2.16}$$

Noting that (2.4) gives

$$\mu_t^n[\phi] = \frac{1}{n}\sum_{i=1}^n \phi(X_t^{i,n})\,, \tag{2.17}$$

we observe upon summing (2.14) over $i$ that

$$\mu_t^n[\phi] = \phi(x_0) + \frac{1}{n}\sum_{i=1}^n \int_0^t \mathscr{L}_s(\mu_s^n)(X_s^{i,n})\,\mathrm{d}s + \frac{1}{n}\sum_{i=1}^n (M^\phi)_t^{i,n} + \frac{1}{n}\sum_{i=1}^n (\bar{M}^\phi)_t^{i,n}\,.$$

Appealing to (2.4) once more shows that

$$\frac{1}{n}\sum_{i=1}^n \mathscr{L}_s(\mu_s^n)\phi(X_s^{i,n}) = \mu_s^n[\mathscr{L}_s(\mu_s^n)\phi]\,,$$

the left-hand side of which is uniformly bounded in view of (2.5). Then linearity of the Lebesgue integral and the previous two displays yield

$$\mu_t^n[\phi] = \phi(x_0) + \int_0^t \mu_s^n[\mathscr{L}_s(\mu_s^n)\phi]\,\mathrm{d}s + \frac{1}{n}\sum_{i=1}^n (M^\phi)_t^{i,n} + \frac{1}{n}\sum_{i=1}^n (\bar{M}^\phi)_t^{i,n}\,. \tag{2.18}$$

This shows that $\mu^n[\phi] = (\mu_t[\phi])_{t\in[0,T]}$ is a real-valued semimartingale whose finite-variation part involves $\mathscr{L}_s$ defined in (2.8). By combining (2.16) with the definition (2.4) of $\mu^n$, we also obtain

$$\frac{1}{n}\sum_{i=1}^{n}(\bar{M}^\phi)_t^{i,n} = \int_0^t \mu_t^n[\mathscr{R}_s(\mu_s^n)\phi]\,\mathrm{d}Z_s^n\,. \tag{2.19}$$

However,

$$\frac{1}{n}\sum_{i=1}^{n}(M^\phi)_t^{i,n} = \frac{1}{n}\sum_{i=1}^{n}\int_0^t \mathscr{U}_s(\mu_s^n)\phi(X_s^{i,n})\,\mathrm{d}B_s^{i,n} \tag{2.20}$$

cannot be written in an analogous way because $i$ appears in $B^{i,n}$.

We next calculate the cross-variation of $\mu^n[\phi]$ with $\mu^n[\psi]$ for any pair of functions $\phi, \psi \in \mathcal{S}$. For this, we note from (2.18) that

$$\langle \mu^n[\phi], \mu^n[\psi]\rangle_t = \frac{1}{n^2}\Big\langle \sum_{i=1}^{n}\Big((M^\phi)^{i,n} + (\bar{M}^\phi)^{i,n}\Big), \sum_{i=1}^{n}\Big((M^\psi)^{i,n} + (\bar{M}^\psi)^{i,n}\Big)\Big\rangle_t\,. \tag{2.21}$$

Now $B^{i,n}$ and $B^{j,n}$ are Brownian motions and independent for $i \neq j$. Hence (2.15) gives

$$\langle (M^\phi)^{i,n}, (M^\psi)^{j,n}\rangle_t = 0 \quad \text{for } i \neq j \tag{2.22}$$

and

$$\langle (M^\phi)^{i,n}, (M^\psi)^{i,n}\rangle_t = \int_0^t \Big(\mathscr{U}_s(\mu_s^n)\phi(X_s^{i,n})\Big)\cdot\Big(\mathscr{U}_s(\mu_s^n)\psi(X_s^{i,n})\Big)\,\mathrm{d}s \quad \text{for } i \in [n]\,, \tag{2.23}$$

where

$$\begin{aligned}\Big(\mathscr{U}_s(\mu_s^n)\phi(X_s^{i,n})\Big)\cdot\Big(\mathscr{U}_s(\mu_s^n)\psi(X_s^{i,n})\Big) &= \Big(\mathscr{U}_s(\mu_s^n)\phi(X_s^{i,n})\Big)^{\mathsf{T}}\Big(\mathscr{U}_s(\mu_s^n)\psi(X_s^{i,n})\Big)\\ &= \sum_{\ell=1}^{d}\Big(\mathscr{U}_s(\mu_s^n)\phi(X_s^{i,n})\Big)_\ell\Big(\mathscr{U}_s(\mu_s^n)\psi(X_s^{i,n})\Big)_\ell\,.\end{aligned}$$

Next, $Z^n$ is a Brownian motion so that (2.16) gives analogously

$$\langle (\bar{M}^\phi)^{i,n}, (\bar{M}^\psi)^{j,n}\rangle_t = \int_0^t \mathscr{R}_s(\mu_s^n)\phi(X_s^{i,n})\cdot\mathscr{R}_s(\mu_s^n)\psi(X_s^{j,n})\,\mathrm{d}s \quad \text{for } i,j \in [n]\,. \tag{2.24}$$

Finally, as $B^{i,n}$ and $Z^n$ are independent Brownian motions, the local martingales $(M^\phi)^{i,n}$ and $(\bar{M}^\psi)^{j,n}$ as well as $(M^\psi)^{j,n}$ and $(\bar{M}^\phi)^{i,n}$ are orthogonal. Plugging (2.22)–(2.24) into (2.21) therefore gives upon using the definition (2.4) of $\mu^n$ that

$$\begin{aligned}\langle \mu^n[\phi], \mu^n[\psi]\rangle_t &= \frac{1}{n^2}\sum_{i=1}^{n}\int_0^t \Big(\mathscr{U}_s(\mu_s^n)\phi(X_s^{i,n})\Big)\cdot\Big(\mathscr{U}_s(\mu_s^n)\psi(X_s^{i,n})\Big)\,\mathrm{d}s\\ &\quad+ \frac{1}{n^2}\sum_{i,j=1}^{n}\int_0^t \Big(\mathscr{R}_s(\mu_s^n)\phi(X_s^{i,n})\Big)\cdot\Big(\mathscr{R}_s(\mu_s^n)\psi(X_s^{j,n})\Big)\,\mathrm{d}s\\ &= \frac{1}{n}\int_0^t \mu_s^n\Big[\Big(\mathscr{U}_s(\mu_s^n)\phi\Big)\cdot\Big(\mathscr{U}_s(\mu_s^n)\psi\Big)\Big]\,\mathrm{d}s + \int_0^t \mu_s^n[\mathscr{R}_s(\mu_s^n)\phi]\cdot\mu_s^n[\mathscr{R}_s(\mu_s^n)\psi]\,\mathrm{d}s\,,\end{aligned} \tag{2.25}$$

where integration of the $\mathbb{R}^d$-valued function $x \mapsto \mathscr{R}_s(\mu_s^n)\phi(x)$ against the measure $\mu_s^n$ is understood as a Bochner integral and thus to be componentwise, that is,

$$\mu_s^n[\mathscr{R}_s(\mu_s^n)\phi] \coloneqq \Big(\mu_s^n\Big[\Big(\mathscr{R}_s(\mu_s^n)\phi\Big)_1\Big], \ldots, \mu_s^n\Big[\Big(\mathscr{R}_s(\mu_s^n)\phi\Big)_d\Big]\Big)\,.$$

This shows that the quadratic variation of $\mu^n[\phi]$ involves the maps $\mathscr{U}_t$ and $\mathscr{R}_t$ defined in (2.10) and (2.9), respectively, and concludes the motivation of the origins of $\mathscr{L}_t$, $\mathscr{R}_t$, $\mathscr{U}_t$.

### 2.3 Martingale problem for the $n$-particle system

The notation introduced in Section 2.2 lets us formulate for each $n \in \mathbb{N}$ a martingale problem in the sense of Section 2.1.

We find it convenient to introduce more succinct notation for the key objects which arose in our calculation of (2.18) and (2.25). These are the *characteristics* of the martingale problem we formulate. For this, we define for each $t \in [0, T]$ the maps

$$\mathscr{A}_t : \mathcal{S}' \to \mathcal{S}' \cup \{\infty\},$$

with $\lambda \mapsto \mathscr{A}_t(\lambda)$, and

$$\mathscr{Q}_t : \mathcal{S}' \to \mathcal{L}^+_{\mathrm{bl}}(\mathcal{S}),$$
$$\mathscr{C}_t : \mathcal{S}' \to \mathcal{L}^+_{\mathrm{bl}}(\mathcal{S}),$$

with $\lambda \mapsto \mathscr{Q}_t(\lambda)$ and $\lambda \mapsto \mathscr{C}_t(\lambda)$, respectively, where we recall that $\mathcal{L}^+_{\mathrm{bl}}(\mathcal{S})$ denotes the space of continuous, nonnegative definite and symmetric bilinear forms $\mathcal{S} \times \mathcal{S} \to \mathbb{R}$; see Section 2.1 for the definition. These maps are explicitly given as follows. We first recall the notation $\lambda[f] = \int f(x)\,\lambda(\mathrm{d}x)$, where the integral is understood to be componentwise if $f$ is an $\mathbb{R}^d$-valued bounded function. For each $\lambda \in \mathbf{M}^+_1$ and $\phi, \psi \in \mathcal{S}$, we set

$$\mathscr{A}_t(\lambda)[\phi] := \lambda[\mathscr{L}_t(\lambda)\phi], \tag{2.26}$$

and $\mathscr{A}_t(\lambda) = \infty$ whenever $\lambda \in \mathcal{S}' \setminus \mathbf{M}^+_1$. Moreover, if $\lambda \in \mathbf{M}^+_1$, we set

$$\mathscr{Q}_t(\lambda)[\phi, \psi] := \lambda[\mathscr{R}_t(\lambda)\phi] \cdot \lambda[\mathscr{R}_t(\lambda)\psi], \tag{2.27}$$
$$\mathscr{C}_t(\lambda)[\phi, \psi] := \lambda\Big[\Big(\mathscr{U}_t(\lambda)\phi\Big) \cdot \Big(\mathscr{U}_t(\lambda)\phi\Big)\Big], \tag{2.28}$$

and $\mathscr{Q}_t(\lambda)[\phi, \psi] = 0$ and $\mathscr{C}_t(\lambda)[\phi, \psi] = 0$ if $\lambda \in \mathcal{S}' \setminus \mathbf{M}^+_1$. On the space $(\Omega, \mathcal{F})$ with the filtration $\mathbb{F} = (\mathcal{F}_t)_{t\in[0,T]}$ and the canonical process $\Lambda = (\Lambda_t)_{t\in[0,T]}$ introduced in Section 1.2, this defines for each $n \in \mathbb{N}$ a martingale problem $\mathrm{MP}(\delta_{x_0}, \mathscr{A}, \mathscr{Q} + \frac{1}{n}\mathscr{C})$ in the sense of Section 2.1 with $\nu_0 = \delta_{x_0}$, $A = \mathscr{A}$ and $Q = \mathscr{Q} + \frac{1}{n}\mathscr{C}$.

**Proposition 2.6** | *Under Assumptions 1.1 and 1.2, the measure $\mathbf{P}^n_{x_0}$ defined by* (1.5) *solves* $\mathrm{MP}(\delta_{x_0}, \mathscr{A}, \mathscr{Q} + \frac{1}{n}\mathscr{C})$ *in the sense of Definition 2.3. In particular, $\Lambda$ on $(\Omega, \mathcal{F}, \mathbf{P}^n_{x_0})$ possesses the canonical semimartingale decomposition*

$$\Lambda_t = \delta_{x_0} + \int_0^t \mathscr{A}_s(\Lambda_s)\,\mathrm{d}s + M_t, \qquad t \in [0, T], \tag{2.29}$$

*with $M$ in $\mathscr{M}^c_2(\mathbb{F}, \mathbf{P}^n_{x_0}; \mathcal{S}')$ and with tensor quadratic variation*

$$\langle\!\langle M \rangle\!\rangle_t = \int_0^t \mathscr{Q}_s(\Lambda_s)\,\mathrm{d}s + \frac{1}{n}\int_0^t \mathscr{C}_s(\Lambda_s)\,\mathrm{d}s, \qquad t \in [0, T]. \tag{2.30}$$

*Moreover, there exists $p_0 > 0$ independent of $n$ such that*

$$\mathbf{P}^n_{x_0}\Big[\Lambda \in \mathbf{C}([0, T];\, \mathscr{H}_{-p})\Big] = 1 \quad \textit{and} \quad M \in \mathscr{M}^c_2(\mathbf{P}^n_{x_0}; \mathscr{H}_{-p}) \tag{2.31}$$

*for any $p \geq p_0$.*

***Proof*** See Appendix B. □

The integrals in (2.29) and (2.30) are understood in the Bochner sense. Proposition 2.6 establishes Step 2 discussed in Section 1.4, where we outlined the strategy of proof for our key results.

**Remark 2.7** | Let $\phi \in \mathcal{S}$ and set $\|\phi\|_m^* = \max_{|\alpha|\le m} \sup_{x\in\mathbb{R}^d} |(1+|x|^2)^{m/2}\partial^\alpha\phi(x)|$ for $m \in \mathbb{N}_0$; see (A.2) for more on the definition of the family of seminorms $\|\cdot\|_m^*$ with $m \in \mathbb{N}$. Under Assumptions 1.1 and 1.2, we obtain for all $\lambda, \lambda' \in \mathbf{M}_1^+$ and $t \in [0,T]$ the bound

$$|\lambda[\mathscr{L}_t(\lambda')\phi]| \le \|b\|_\infty\|\phi\|_1^* + \frac{1}{2}(\|\sigma\|_\infty^2 + \|\bar\sigma\|_\infty^2)\|\phi\|_2^* \le c_{b,\sigma,\bar\sigma}\|\phi\|_2^* \tag{2.32}$$

with $c_{b,\sigma,\bar\sigma} \coloneqq \|b\|_\infty + \frac{1}{2}(\|\sigma\|_\infty^2 + \|\bar\sigma\|_\infty^2) < \infty$. This is an estimate that we rely on repeatedly, for instance in the following form. By construction, $\mathbf{P}_{x_0}^n$ puts all mass on the set $\mathbf{C}([0,T];\mathcal{S}'\cap\mathbf{M}_1^+)$. Using Proposition 2.6 and (2.26), we can evaluate (2.29) under $\mathbf{P}_{x_0}^n$ to get for each $t \in [0,T]$ that

$$M_t[\phi] = \Lambda_t[\phi] - \delta_{x_0}[\phi] - \Big(\int_s^t \mathscr{A}_u(\Lambda_u)\Big)[\phi]\,\mathrm{d}u = \Lambda_t[\phi] - \delta_{x_0}[\phi] - \int_s^t \Lambda_u[\mathscr{L}_u(\Lambda_u)\phi]\,\mathrm{d}u\,.$$

Under Assumptions 1.1 and 1.2, the estimate (2.32) now shows that we can bound the process $M[\phi] \coloneqq (M_t[\phi])_{t\in[0,T]}$ uniformly on the set $[0,T]\times\mathbf{C}([0,T];\mathcal{S}'\cap\mathbf{M}_1^+)$ by

$$\sup_{t\in[0,T]} |M_t[\phi]| \le 2\|\phi\|_0^* + Tc_{b,\sigma,\bar\sigma}\|\phi\|_2^* \le (2+Tc_{b,\sigma,\bar\sigma})\|\phi\|_2^*. \tag{2.33}$$

As a consequence, $M[\phi]$ is a true $(\mathbb{F},\mathbf{P}_{x_0}^n)$-martingale for each $n \in \mathbb{N}$. In fact, using (2.1) and Lemma A.2 we show in in Step 6 of the proof of Proposition 2.6 that for an appropriate $p > 0$ we get that

$$\|M_T\|_{\mathscr{H}_{-p}}^2 \le c(2+Tc_{b,\sigma,\bar\sigma})\,. \tag{2.34}$$

Together with the regularization result in Lemma 2.2 this implies that $M \in \mathscr{M}_2^{\mathrm{c}}(\mathbf{P}_{x_0}^n;\mathscr{H}_{-p})$.

Since by (1.6) in Lemma 1.4, any narrow cluster point $\mathbf{P}_{x_0}^\infty$ of $(\mathbf{P}_{x_0}^n)_{n\in\mathbb{N}}$ puts all mass on the set $\mathbf{C}([0,T];\mathcal{S}'\cap\mathbf{M}_1^+)$, the estimates in (2.32), (2.33) remain valid for $\mathbf{P}_{x_0}^\infty$. Thus if $M[\phi]$ happens to have the $(\mathbb{F},\mathbf{P}_{x_0}^\infty)$-local martingale property, then $M[\phi]$ is also a bounded and thus a true martingale.

# 3 REGULARITY ESTIMATES

In this section, we revisit the *emergence of regularity*, which is the focus of the first paper in this series; see Crowell [12]. There, we showed that for every $t \in (0,T]$, the random probability measure $\Lambda_t$ under $\mathbf{P}_{x_0}^\infty$ is almost surely absolutely continuous with respect to Lebesgue measure, with a density that possesses certain analytic regularity properties. This insight is an essential ingredient for our study here.

Next to the emergence of regularity result, we also require certain uniform regularity estimates for the time-marginal density of a single particle in the $n$-particle system (1.1), (1.2); see Lemma 3.3 below. This result is technically simpler because the time-marginals are nonrandom. In principle, such estimates can be obtained in different ways; see the comments after the proof of Lemma 3.3. For reasons of expositional continuity, we derive it using the techniques familiar from [12]. For this we recall the abstract interpolation result in Proposition 3.2 below.

### 3.1 Revisiting the emergence of regularity

Let $\mathsf{H}^s_r(\mathbb{R}^d)$ for $r \in (1,\infty)$ and $s \in \mathbb{R}$ be the Bessel potential spaces; see e.g. Bergh and Löfström [6, Ch. 6] or [12, Appendix A]. For $u, s > 0$ and $w \in (-u, s)$, we have the chain of continuous inclusions

$$\mathcal{S}(\mathbb{R}^d) \hookrightarrow \mathsf{H}^s_r(\mathbb{R}^d) \hookrightarrow \mathsf{H}^w_r(\mathbb{R}^d) \hookrightarrow \mathsf{H}^{-u}_r(\mathbb{R}^d) \hookrightarrow \mathcal{S}'(\mathbb{R}^d). \tag{3.1}$$

In this setting, $\mathsf{H}^w_r(\mathbb{R}^d)$ is called an *intermediate space* between $\mathsf{H}^s_r(\mathbb{R}^d)$ and $\mathsf{H}^{-u}_r(\mathbb{R}^d)$. Recall also that $\mathsf{H}^0_r(\mathbb{R}^d)$ equals $\mathbb{L}^r(\mathbb{R}^d, \mathrm{d}x)$ with equivalent norms; so for $w \geq 0$, elements of $\mathsf{H}^w_r(\mathbb{R}^d)$ are functions, not merely distributions.

In [12, Thm. 1], we obtained the following result.

**Theorem 3.1** | *Impose Assumptions 1.1 and 1.2. Then there exist real numbers $w > 0$ and $r > 1$ such that for any cluster point $\mathbf{P}^\infty_{x_0}$ of the sequence $(\mathbf{P}^n_{x_0})_{n\in\mathbb{N}}$, we have*

$$\mathbf{P}^\infty_{x_0}\left[\left\{\omega \in \Omega \,:\, \frac{\mathrm{d}\Lambda_t(\omega)}{\mathrm{d}x} = \bar{p}(t,\omega) \text{ is in } \mathsf{H}^w_r(\mathbb{R}^d) \text{ for almost all } t \in (0,T]\right\}\right] = 1, \tag{3.2}$$

*where $\bar{p} : (0,T] \times \Omega \to \mathsf{H}^w_r(\mathbb{R}^d)$ is a strongly measurable function.*

*More precisely, for each $t \in (0,T]$, we have $\mathbf{P}^\infty_{x_0}$-a.s. that $\Lambda_t \ll \mathrm{d}x$, and for any $q \geq 1$, $w$ and $r$ can be chosen in such a way that there exists $1 > \gamma > 0$ such that the function $t \mapsto (\omega \mapsto \mathrm{d}\Lambda_t(\omega)/\mathrm{d}x)$ defines a strongly measurable map*

$$p : (0,T] \to \mathbb{L}^q\Big((\Omega, \mathcal{F}, \mathbf{P}^\infty_{x_0});\, \mathsf{H}^w_r(\mathbb{R}^d)\Big) \tag{3.3}$$

*which satisfies for some constant $c_{\text{Thm. 3.1}} < \infty$ the bound*

$$\Big\| \|p(t)\|_{\mathsf{H}^w_r(\mathbb{R}^d)} \Big\|_{\mathbb{L}^q(\mathbf{P}^\infty_{x_0})} \leq c_{\text{Thm. 3.1}} (1 \wedge t)^{-\gamma} \tag{3.4}$$

*for all $t \in (0,T]$. In (3.4), $c_{\text{Thm. 3.1}}$ depends on $r$, $w$, $q$, $\gamma$, but not on $t$ and neither on $\mathbf{P}^\infty_{x_0}$.*

*Finally, $\bar{p}$ is unique up to $(\mathrm{d}t \otimes \mathrm{d}\mathbf{P}^\infty_{x_0})$-a.e. equality and we have that $\bar{p}(t,\,\cdot\,) = p(t)$ in $\mathbb{L}^q((\Omega, \mathcal{F}, \mathbf{P}^\infty_{x_0}); \mathsf{H}^w_r(\mathbb{R}^d, \mathrm{d}x))$, for almost every $t \in (0,T]$. In particular, $\bar{p}(t,\omega)$ is the density of $\Lambda_t(\omega)$ $(\mathrm{d}t \otimes \mathrm{d}\mathbf{P}^\infty_{x_0})$-a.e., and $\bar{p}(t,\,\cdot\,)$ satisfies the bound (3.4) for almost every $t \in (0,T]$.*

One of the difficulties encountered in the proof of Theorem 3.1 is the randomness of $\Lambda$. In [12], we devise a *direct strategy* to work with the sequence $(\mathbf{P}^n_{x_0})_{n\in\mathbb{N}}$ via an *approximation and interpolation* scheme, which has its foundations in the seminal contribution of Fournier and Printems [23], with subsequent works such as Debussche and Romito [18] or Bally and Caramellino [4]; see [12] for a fuller discussion and additional references. In broad terms, we start with an Euler-type approximation of the particle systems (1.1), (1.2), and use this to define an approximated empirical measure. By jointly passing the empirical measure and its approximation to narrow limits, we can decompose $(\Lambda_t)_{t\in[0,T]}$ into a sum of two-parameter processes $(A_{\varepsilon,t}+E_{\varepsilon,t})_{0\le\varepsilon\le t\le T}$, consisting of an approximation and an error term, respectively. To derive regularity estimates for the terms $A$ and $E$, we need probabilistic arguments building on the particle systems and a limiting procedure. The approximation $A_{\varepsilon,t}$ is a probability measure that, for elliptic $\sigma$ and $\varepsilon > 0$, has a smooth density. However, as $\varepsilon \to 0$, it becomes singular, but the rate at which this occurs can be bounded. The error $E_{\varepsilon,t} = \Lambda_t - A_{\varepsilon,t}$, on the other hand, is a signed measure for which we can bound the rate at which $E_{\varepsilon,t}$ vanishes as $\varepsilon \to 0$. To obtain an estimate for $\Lambda_t$, we then apply an interpolation argument from harmonic analysis to trade off the rate at which the smooth part $A_{\varepsilon,t}$ becomes singular against the rate at which $E_{\varepsilon,t}$ vanishes as $\varepsilon \to 0$. The randomness introduced by the common noise makes the estimates for $E_{\varepsilon,t}$ random and the estimates for $A_{\varepsilon,t}$ and $E_{\varepsilon,t}$ available only for certain pairs $(\varepsilon, t)$. This restricts how we can combine the estimates. To overcome this difficulty, we built in [12, Prop. 5.1] the following workhorse interpolation result.

**Proposition 3.2** | *Consider reals $s > 0$, $1 < r < 2$ and $\xi > 0$, take $r'$ conjugate to $r$ and set $u := 1 + d/r'$. Let $\lambda \in \mathcal{S}'(\mathbb{R}^d)$ be a distribution and assume that there exists $\varepsilon_0 > 0$ such that for each $\varepsilon \in (0, \varepsilon_0]$, we have a decomposition*

$$\lambda = a_\varepsilon + e_\varepsilon \tag{3.5}$$

*for some $a_\varepsilon$ and $e_\varepsilon$ in $\mathcal{S}'(\mathbb{R}^d)$ satisfying*

$$\|a_\varepsilon\|_{\mathsf{H}^s_r(\mathbb{R}^d)} \le c_a \varepsilon^{-(d/r'+s)/2} \quad \text{and} \quad \|e_\varepsilon\|_{\mathsf{H}^{-u}_r(\mathbb{R}^d)} \le c_e \varepsilon^{(1+\xi)/2} \tag{3.6}$$

*with constants $c_a$, $c_e < \infty$ that may depend on $\varepsilon_0$, but are independent of $\varepsilon \neq \varepsilon_0$. Then the following hold:*

1) *There exists $w_0 = w_0(d, \xi, s, r) \in (-u, s)$ such that $\lambda \in \mathsf{H}^w_r(\mathbb{R}^d)$ for all $w < w_0$. More specifically, we have*

$$\|\lambda\|_{\mathsf{H}^w_r(\mathbb{R}^d)} \le c\big(c_a \varepsilon_0^{-(d/r'+s)/2} + c_e \varepsilon_0^{(1+\xi)/2}\big) \tag{3.7}$$

*for a constant $c = c(w, d, \xi, s, r) < \infty$ independent of $\varepsilon_0$.*

2) *For all $\xi > 0$, there exist combinations of parameters $s > 0$ and $1 < r < 2$ satisfying*

$$\gamma := (d/r' + s)/2 < 1$$

*and such that if the decomposition (3.5) and the bounds in (3.6) hold, then we have $w_0(d, \xi, s, r) \in (0, s)$ in part 1). In particular, $\lambda$ is then a function in some $\mathsf{H}^w_r(\mathbb{R}^d)$ with $w > 0$.*

3) *For all $\xi > 0$, $s > 0$ and $1 < r < 2$, there exists a number $\alpha_0 = \alpha_0(\xi, r, s) > 0$ such that if the decomposition in (3.5) and the norm bounds in (3.6) are valid only for*

$$\varepsilon \in \{2^{-\alpha_0 n}\varepsilon_0 \,:\, n \in \mathbb{N}_0\} =: (0, \varepsilon_0]_{\alpha_0} \subseteq (0, \varepsilon_0]\,,$$

*then the conclusions from parts 1) and 2) already hold.*

Part 3) of the previous result is used in [12] to prove Theorem 3.1.

### 3.2 An additional uniform regularity estimate for the proof of Theorem 1.6

To prove Theorem 1.6 we require an additional egularity estimate for the time-marginale density of the a single particle in the $n$-particle systems, which is uniform in $n$. To obtain this estimate we follow the same broad strategy that we just outlined for the proof of Theorem 3.1. The difference now is that the time-marginal density is not random, so the difficulty we mentioned above can be avoided. This lets us use directly part 2) of Proposition 3.2.

**Lemma 3.3** | *Impose Assumptions 1.1 and 1.2. Then for the weak solution $\boldsymbol{X}^n$ of* (1.1), (1.2)*, we have for any $t \in (0,T]$ and $n \in \mathbb{N}$ that $\mathrm{Law}_{\mathbb{P}^n_{x_0}}(X^{1,n}_t) \ll \mathrm{d}x$. In fact, there exist exponents*

$$1 < r < \infty \text{ and } 0 < \gamma < 1\,,$$

*a constant $c = c_{\text{Lem. 3.3}} < \infty$ and for each $n \in \mathbb{N}$ a measurable*

$$p^{1,n} : (0,T] \to \mathbb{L}^r(\mathrm{d}x) \tag{3.8}$$

*which jointly satisfy the bound*

$$\sup_{n\in\mathbb{N}} \|p^{1,n}(t)\|_{\mathbb{L}^r(\mathbb{R}^d)} \le c(1 \wedge t)^{-\gamma} \quad \text{for all } t \in (0,T] \tag{3.9}$$

*and for which $\mathrm{dLaw}_{\mathbb{P}^n_{x_0}}(X^{1,n}_t)/\mathrm{d}x = p^{1,n}(t)$ in $\mathbb{L}^r(\mathrm{d}x)$ for $t \in (0,T]$.*

***Proof*** Starting from the weak solution $\boldsymbol{X}^n$ of (1.1), (1.2) defined on $(\Omega^n_{\text{par}}, \mathcal{G}^n, \mathbb{P}^n_{x_0})$, consider for $0 < \varepsilon \le t$ the auxiliary random variables

$$Y^{1,n;\varepsilon,t}_t = X^{1,n}_{t-\varepsilon} + \int_{t-\varepsilon}^t \sigma_r(X^{1,n}_{t-\varepsilon}, \mu^n_{t-\varepsilon})\,\mathrm{d}B^{1,n}_r + \bar{\sigma}_{t-\varepsilon}(X^{1,n}_{t-\varepsilon}, \mu^n_{t-\varepsilon})(Z^n_t - Z^n_{t-\varepsilon})\,. \tag{3.10}$$

This gives a special case of the auxiliary process used in [12, Eqn. (1.17)]. We want to apply Proposition 3.2 with $\varepsilon_0 \coloneqq 1 \wedge t$ to the decomposition $\lambda = a_\varepsilon + e_\varepsilon$ from (3.5) using

$$\begin{aligned}
\lambda &\coloneqq \mathrm{Law}_{\mathbb{P}^n_{x_0}}(X^{1,n}_t)\,,\\
a_\varepsilon &\coloneqq \mathrm{Law}_{\mathbb{P}^n_{x_0}}(Y^{1,n;\varepsilon,t}_t)\,,\\
e_\varepsilon &\coloneqq \mathrm{Law}_{\mathbb{P}^n_{x_0}}(X^{1,n}_t) - \mathrm{Law}_{\mathbb{P}^n_{x_0}}(Y^{1,n;\varepsilon,t}_t)\,.
\end{aligned}$$

Since each of the above terms is a finite signed measure on $\mathbb{R}^d$, $\lambda = a_\varepsilon + e_\varepsilon$ is indeed a decomposition consisting of elements in $\mathcal{S}'$; they act on functions $\phi \in \mathcal{S}$ by

$$\begin{aligned}
a_\varepsilon[\phi] &= \mathbf{E}^{\mathbb{P}^n_{x_0}}[\phi(Y^{1,n;\varepsilon,t}_t)]\,,\\
e_\varepsilon[\phi] &= \mathbf{E}^{\mathbb{P}^n_{x_0}}[\phi(X^{1,n}_t) - \phi(Y^{1,n;\varepsilon,t}_t)]\,.
\end{aligned}$$

To verify (3.6), we next derive quantitative estimates for the relevant norms. First, by the definition of $a_\varepsilon$ and the tower property,

$$a_\varepsilon[\phi] = \mathbf{E}^{\mathbb{P}^n_{x_0}}[\phi(Y^{1,n;\varepsilon,t}_t)] = \mathbf{E}^{\mathbb{P}^n_{x_0}}\left[\mathbf{E}^{\mathbb{P}^n_{x_0}}[\phi(Y^{1,n;\varepsilon,t}_t) \,|\, \mathcal{F}^{X^{1,n},\mu^n}_{t-\varepsilon}]\right], \tag{3.11}$$

where

$$\mathcal{F}_{t-\varepsilon}^{X^{1,n},\mu^n} = \sigma\Big((X_s^{1,n},\mu_s^n) : 0 \le s \le t-\varepsilon\Big).$$

Recall that $\mathcal{N}(m,v)$ denotes the $d$-dimensional normal distribution with mean vector $m$ and covariance matrix $v$, and $g(y;m,v)$ is its density on $\mathbb{R}^d$. We note from (3.10) that

$$Y_t^{1,n;\varepsilon,t} \,|\, \mathcal{F}_{t-\varepsilon}^{X^{1,n},\mu^n} \overset{\mathrm{d}}{=} \mathcal{N}\Big(X_{t-\varepsilon}^{1,n}, \Sigma_{\varepsilon,t}^{1,n}(X_{t-\varepsilon}^{1,n},\mu_{t-\varepsilon}^n)\Big)$$

with covariance matrix given by

$$\Sigma_{\varepsilon,t}^{1,n}(X_{t-\varepsilon}^{1,n},\mu_{t-\varepsilon}^n) := \int_{t-\varepsilon}^t (\sigma_r\sigma_r^\mathsf{T})(X_{t-\varepsilon}^{1,n},\mu_{t-\varepsilon}^n)\,\mathrm{d}r + \varepsilon(\bar\sigma_{t-\varepsilon}\bar\sigma_{t-\varepsilon}^\mathsf{T})(X_{t-\varepsilon}^{1,n},\mu_{t-\varepsilon}^n)\,.$$

So we can write the conditional expectation in (3.11) $\mathbb{P}_{x_0}^n$-a.s. as

$$\mathbf{E}^{\mathbb{P}_{x_0}^n}[\phi(Y_t^{1,n;\varepsilon,t}) \,|\, \mathcal{F}_{t-\varepsilon}^{X^{1,n},\mu^n}] = \int \phi(y) g\Big(y; X_{t-\varepsilon}^{1,n}, \Sigma_{\varepsilon,t}^{1,n}(X_{t-\varepsilon}^{1,n},\mu_{t-\varepsilon}^n)\Big)\,\mathrm{d}y.$$

Taking expectations and using (3.11) and the duality between $\mathsf{H}_r^s$ and $\mathsf{H}_{r'}^{-s}$ from [12, Lem. 2.5] for $1/r + 1/r' = 1$ gives us the bound

$$\begin{aligned} |a_\varepsilon[\phi]| &= \left| \mathbf{E}^{\mathbb{P}_{x_0}^n}\left[ \int_{\mathbb{R}^d} \phi(y) g\Big(y; X_{t-\varepsilon}^{1,n}, \Sigma_{\varepsilon,t}^{1,n}(X_{t-\varepsilon}^{1,n},\mu_{t-\varepsilon}^n)\Big)\,\mathrm{d}y\right]\right| \\ &\le c(r,s)\|\phi\|_{\mathsf{H}_{r'}^{-s}}\, \mathbf{E}^{\mathbb{P}_{x_0}^n}\left[\Big\| g\Big(\,\cdot\,; X_{t-\varepsilon}^{1,n}, \Sigma_{\varepsilon,t}^{1,n}(X_{t-\varepsilon}^{1,n},\mu_{t-\varepsilon}^n)\Big)\Big\|_{\mathsf{H}_r^s}\right]. \end{aligned} \tag{3.12}$$

Our next goal is to obtain a bound for the $\mathsf{H}_r^s$-norm inside the expectation. Since the $\mathsf{H}_r^s$-norm is translation invariant, see [12, Lem. 3.3], we can reduce to the case of centered Gaussians, i.e.

$$\Big\| g\Big(\,\cdot\,; X_t^{i,n}, \Sigma_{\varepsilon,t}^{1,n}(X_{t-\varepsilon}^{1,n},\mu_{t-\varepsilon}^n)\Big)\Big\|_{\mathsf{H}_r^s} = \Big\| g\Big(\,\cdot\,; 0, \Sigma_{\varepsilon,t}^{1,n}(X_{t-\varepsilon}^{1,n},\mu_{t-\varepsilon}^n)\Big)\Big\|_{\mathsf{H}_r^s}. \tag{3.13}$$

To bound the right-hand side of (3.13), we next use the scaling properties of the Gaussian density. For this, we write

$$\begin{aligned} \tilde\Sigma_{\varepsilon,t}^{1,n}(X_{t-\varepsilon}^{1,n},\mu_{t-\varepsilon}^n) &:= \frac{1}{\varepsilon}\Sigma_{\varepsilon,t}^{1,n}(X_{t-\varepsilon}^{1,n},\mu_{t-\varepsilon}^n) \\ &= \frac{1}{\varepsilon}\int_{t-\varepsilon}^t (\sigma_r\sigma_r^\mathsf{T})(X_{t-\varepsilon}^{1,n},\mu_{t-\varepsilon}^n)\,\mathrm{d}r + (\bar\sigma_{t-\varepsilon}\bar\sigma_{t-\varepsilon}^\mathsf{T})(X_{t-\varepsilon}^{1,n},\mu_{t-\varepsilon}^n) \end{aligned}$$

and note that since a sum of symmetric and uniformly elliptic matrices is again symmetric and uniformly elliptic, the same computation as in [12, Proof of Prop. 3.2] shows that the family

$$\mathcal{V} = \{\tilde\Sigma_{\varepsilon,t}(x,\lambda) \,:\, 0 < \varepsilon \le t, x \in \mathbb{R}^d, \lambda \in \mathbf{M}_1^+\}$$

is uniformly elliptic and bounded. This is the ingredient needed to apply [12, Lem. 3.4], which gives us a constant $c(r,s) < \infty$ such that

$$\Big\| g\Big(\,\cdot\,; 0, \tilde\Sigma_{\varepsilon,t}^{1,n}(X_{t-\varepsilon}^{1,n},\mu_{t-\varepsilon}^n)\Big)\Big\|_{\mathsf{H}_r^s} \le c(r,s)\,.$$

Returning to (3.13), we can follow the same steps leading to [12, En. (3.31)]. Indeed, using first the definition of the normal density $g$, then the scaling property of the $\mathsf{H}^s_r$-norm from [12, Lem. 3.3] and finally the bound in the display just above leads to

$$\begin{aligned}\Big\|g\Big(y;0,\Sigma^{1,n}_{\varepsilon,t}(X^{1,n}_{t-\varepsilon},\mu^n_{t-\varepsilon})\Big)\Big\|_{\mathsf{H}^s_r(\mathrm{d}y)} &= \varepsilon^{-d/2}\Big\|g\Big(\varepsilon^{-1/2}y;0,\tilde\Sigma^{1,n}_{\varepsilon,t}(X^{1,n}_{t-\varepsilon},\mu^n_{t-\varepsilon})\Big)\Big\|_{\mathsf{H}^s_r(\mathrm{d}y)}\\ &= \varepsilon^{-(s-d/r)/2}\varepsilon^{-d/2}\Big\|g\Big(y;0,\tilde\Sigma^{1,n}_{\varepsilon,t}(X^{1,n}_{t-\varepsilon},\mu^n_{t-\varepsilon})\Big)\Big\|_{\mathsf{H}^s_r(\mathrm{d}y)}\\ &\le c(r,s)\varepsilon^{-(s+d/r')/2}\,, \end{aligned}\tag{3.14}$$

provided $0<\varepsilon\le 1\wedge t$, which is required to use the scaling for the $\mathsf{H}^s_r$-norm. The bound in (3.14) is deterministic. Inserting it into (3.12) produces the deterministic estimate

$$|a_\varepsilon[\phi]|\le \|\phi\|_{\mathsf{H}^{-s}_{r'}}c(r,s)\varepsilon^{-(s+d/r')/2}\,,$$

and the dual characterization of the $\mathsf{H}^s_r$-norm from [12, Lem. 2.5] gives

$$\|a_\varepsilon\|_{\mathsf{H}^s_r}=\sup\{|a_\varepsilon[\phi]|\,:\,\|\phi\|_{\mathsf{H}^{-s}_{r'}}\le 1\}\le c(r,s)\varepsilon^{-(d/r'+s)/2}\,.\tag{3.15}$$

We have thus established the first inequality in (3.6).

Let us now recall that $[\phi]_{\mathsf{H\ddot{o}l}_\rho}=\sup_{x\neq y}(|\phi(x)-\phi(y)|)/(|x-y|^\rho)$ with $0<\rho<1$ denotes the Hölder-seminorm, and from the Sobolev embedding in [12, Lem. 2.4] that $\|\cdot\|_\infty+[\,\cdot\,]_{\mathsf{H\ddot{o}l}_\rho}\le c_\rho\|\cdot\|_{\mathsf{H}^u_{r'}}$, where we set $u:=1+d/r'$. For each $\phi\in\mathcal{S}$, we then see that

$$\begin{aligned}|e_\varepsilon[\phi]| &= |\,\mathbf{E}^{\mathbb{P}^n_{x_0}}[\phi(X^{1,n}_t)-\phi(Y^{1,n;\varepsilon,t}_t)]|\\ &\le [\phi]_{\mathsf{H\ddot{o}l}_\rho}\mathbf{E}^{\mathbb{P}^n_{x_0}}[|X^{1,n}_t-Y^{1,n;\varepsilon,t}_t|^\rho]\\ &\le c(\rho,r,u)\|\phi\|_{\mathsf{H}^u_{r'}}\mathbf{E}^{\mathbb{P}^n_{x_0}}[|X^{1,n}_t-Y^{1,n;\varepsilon,t}_t|^\rho]\,.\end{aligned}$$

To bound the expectation in the last line, we apply Jensen's inequality to find for any $q\ge 1$ that

$$\mathbf{E}^{\mathbb{P}^n_{x_0}}[|X^{1,n}_t-Y^{1,n;\varepsilon,t}_t|^\rho]\le \|X^{1,n}_t-Y^{1,n;\varepsilon,t}_t\|^\rho_{\mathbb{L}^q(\mathbb{P}^n_{x_0})}\,.$$

Next, recall that $\beta$ denotes the Hölder-regularity of $\sigma$, $\bar\sigma$ from Assumption 1.2. We now apply [12, Prop. 2.2], which established that for any $q\ge 1$, there exists $c(q,b,\sigma,\bar\sigma,T,\beta)<\infty$ with

$$\|X^{1,n}_t-Y^{1,n;\varepsilon,t}_t\|_{\mathbb{L}^q(\mathbb{P}^n_{x_0})}\le c(q,b,\sigma,\bar\sigma,T,\beta)\varepsilon^{\frac12+\frac\beta2}$$

for all $n\in\mathbb{N}$ and $0\le\varepsilon\le t$. All this gives us for $q\ge 1$ that

$$|e_\varepsilon[\phi]|\le c(q,b,\sigma,\bar\sigma,T,\beta,\rho,r,u)\|\phi\|_{\mathsf{H}^u_{r'}}\varepsilon^{\rho(1+\beta)/2}\,.$$

Since $\beta>0$ from Assumption 1.2, we can choose $\rho<1$ sufficiently large guarantee that we have $\xi:=\rho(1+\beta)-1>0$. With this choice, the dual characterization of the $\mathsf{H}^{-u}_r$-norm from [12, Lem. 2.5] finally gives

$$\|e_\varepsilon\|_{\mathsf{H}^{-u}_r}\le c(q,b,\sigma,\bar\sigma,T,\beta,\rho,r,u)\varepsilon^{(1+\xi)/2}\tag{3.16}$$

for all $0<\varepsilon\le t$, and so the second inequality in (3.6) holds.

The quantitative estimates in (3.15) and (3.16) match precisely those appearing in Proposition 3.2, by part 2) of which we can find real numbers $r > 1$ and $s > 0$ to have $\gamma := (d/r' + s)/2 < 1$ and

$$\|\lambda\|_{\mathsf{H}^w_r(\mathbb{R}^d)} \le c\Big(c(q, b, \sigma, \bar{\sigma}, T, \beta, \rho, r, u)(1 \wedge t)^{(1+\xi)/2} + c(r, s)(1 \wedge t)^{-\gamma}\Big)$$

for $w \in (0, w_0)$ with $w_0 = w_0(d, \xi, s, r)$ and $c = c(w, d, \xi, s, r) < \infty$ as in Proposition 3.2. Note that the choice of $w$ above can be made independent of $t \in (0, T]$ and $n \in \mathbb{N}$ since the bounds in (3.15) and (3.16) are uniform in those quantities. Hence also $c$ does not depend on $t$ nor $n$. Since (3.1) gives $\|\lambda\|_{\mathbb{L}^r(\mathbb{R}^d, \mathrm{d}x)} \le c(r, w, d)\|\lambda\|_{\mathsf{H}^w_r(\mathbb{R}^d)}$ and since $(1 \wedge t)^{(1+\xi)/2} \le (1 \wedge t)^{-\gamma}$, we get

$$\|\lambda\|_{\mathbb{L}^r(\mathbb{R}^d, \mathrm{d}x)} \le c_{\text{Lem. 3.3}}(1 \wedge t)^{-\gamma} \tag{3.17}$$

for all $t \in (0, T]$ and $n \in \mathbb{N}$, where $c_{\text{Lem. 3.3}} < \infty$ is an appropriate global constant that is independent of $t$ and $n$. The estimate (3.17) is similar to [12, Eqn. (6.14)] in the proof of Theorem 3.1. In fact, if $\Lambda$ in [12, Eqn. (6.14)] does not depend on $\omega$, then that bound reduces exactly to (3.17). With this observation, we can closely follow Steps 3–5 in [12, Proof of Thm. 1.5] with straightforward adaptations to conclude the proof. □

Since the proof of Lemma 3.3 uses only part 2) of Proposition 3.2, but not part 3), Lemma 3.3 can be proved by appealing to results that predate Proposition 3.2, e.g. those developed in Bally and Caramellino [4] or Romito [55]. In fact, there is no need to use the approximation and interpolation scheme. Other common techniques to obtain estimates for the time marginals of SDEs, such as those in Stroock and Varadhan [63] or Aronson [3], see also Porper and Èidel'man [51], can be used. We choose to use Proposition 3.2 for reasons of expositional continuity.

# 4 FORMALLY LIMITING MARTINGALE PROBLEM

In this section we prove our main result which we repeat here for sake of reference.

**Theorem 1.6** | *Under Assumptions 1.2 and 1.5, any narrow cluster point $\mathbf{P}^\infty_{x_0}$ of the sequence $(\mathbf{P}^n_{x_0})_{n\in\mathbb{N}}$ solves the martingale problem $\mathrm{MP}(\delta_{x_0}, \mathscr{A}, \mathscr{Q})$. In particular, the coordinate process $\Lambda$ on $(\Omega, \mathcal{F}, \mathbf{P}^\infty_{x_0})$ possesses the $\mathcal{S}'$-valued canonical semimartingale decomposition*

$$\Lambda_t = \delta_{x_0} + \int_0^t \mathscr{A}_s(\Lambda_s)\,\mathrm{d}s + M_t\,, \qquad t \in [0, T]\,, \tag{1.13}$$

*where $M = (M_t)_{t\in[0,T]}$ is in $\mathscr{M}^c_2(\mathbb{F}, \mathbf{P}^\infty_{x_0}; \mathcal{S}')$ with quadratic variation*

$$\langle\!\langle M \rangle\!\rangle_t = \int_0^t \mathscr{Q}_s(\Lambda_s)\,\mathrm{d}s\,, \qquad t \in [0, T]\,. \tag{1.14}$$

*Finally, the* empirical propagation of chaos property *holds, that is,*

$$\mathrm{Law}_{\mathbb{P}^{n_k}_{x_0}}(\mu^{n_k}) \longrightarrow \mathrm{Law}_{\mathbf{P}^\infty_{x_0}}(\Lambda) \textit{ narrowly in } \mathbf{M}^+_1\Big(\mathbf{C}([0, T]; \mathcal{S}')\Big) \tag{1.15}$$

*as $k \to \infty$ for any subsequence $(n_k)_{k\in\mathbb{N}}$ such that $\mathbf{P}^{n_k}_{x_0} \to \mathbf{P}^\infty_{x_0}$.*

The strategy of proof was already touched upon in Section 1.3. There we noted that in order to prove the $(\mathbb{F}, \mathbf{P}^{\infty}_{x_0})$-local martingale property of $M$ in (1.13) we should aim to establish (1.22). The following technical proposition gives a more general result from which (1.22) is obtained; see (4.4) below.

**Proposition 4.2** | *Under Assumptions 1.5 and 1.2, let $(n_k)_{k\in\mathbb{N}}$ be a sequence such that $(\mathbf{P}^{n_k}_{x_0})_{k\in\mathbb{N}}$ converges to $\mathbf{P}^{\infty}_{x_0}$ in the narrow topology. Then we have*

$$\lim_{k\to\infty} \mathbf{E}^{\mathbf{P}^{n_k}_{x_0}}\Big[\int_s^t g_u(\Lambda)\mathscr{A}_u(\Lambda_u)[\phi]\,\mathrm{d}u\Big] = \mathbf{E}^{\mathbf{P}^{\infty}_{x_0}}\Big[\int_s^t g_u(\Lambda)\mathscr{A}_u(\Lambda_u)[\phi]\,\mathrm{d}u\Big] \tag{4.1}$$

*for any $\phi \in \mathcal{S}$, all $s,t \in [0,T]$ with $s \le t$, and all bounded measurable functions $g : [0,T]\times\Omega \to \mathbb{R}$ such that for each $u \in [0,T]$, the real-valued function $g_u(\,\cdot\,) := g(u,\,\cdot\,)$ on $\Omega$ is continuous.*

The proof of this result will occupy us for the better half of this section. For it, we present two technical ingredients.

**Almost sure representation** Consider a narrow cluster point $\mathbf{P}^{\infty}_{x_0}$ of $(\mathbf{P}^{n}_{x_0})_{n\in\mathbb{N}}$ and a subsequence $(\mathbf{P}^{n_k}_{x_0})_{k\in\mathbb{N}}$ with $\mathbf{P}^{n_k}_{x_0} \to \mathbf{P}^{\infty}_{x_0}$ narrowly in $\mathbf{M}_1^+(\mathbf{C}([0,T];\mathcal{S}'))$. To simplify notation, we slightly abuse it and relabel the sequence by setting $\mathbf{P}^{k}_{x_0} := \mathbf{P}^{n_k}_{x_0}$ for all $k \in \mathbb{N}$. Next, take $p > p_0$ with $p_0$ as in the proof of Proposition 2.6. This ensures that all elements of the sequence $(\mathbf{P}^{k}_{x_0})_{k\in\mathbb{N}}$ and $\mathbf{P}^{\infty}_{x_0}$ concentrate all their mass on $\mathbf{C}([0,T];\,\mathscr{H}_{-p})$, which unlike $\mathbf{C}([0,T];\mathcal{S}')$ is separable. By the Skorohod representation theorem, see e.g. Pollard [50, Thm. IV.13], we can thus produce an appropriate probability space $(\tilde{\Omega}, \tilde{\mathcal{F}}, \tilde{\mathbf{P}})$ and a sequence of random variables $\tilde{\Lambda}^k$ for $k \in \mathbb{N}\cup\{\infty\}$, all valued in $\mathbf{C}([0,T];\mathscr{H}_{-p})$ and defined on $(\tilde{\Omega}, \tilde{\mathcal{F}}, \tilde{\mathbf{P}})$, satisfying

$$\lim_{k\to\infty} \tilde{\Lambda}^k = \tilde{\Lambda}^\infty \qquad \text{in } \mathbf{C}([0,T];\,\mathscr{H}_{-p}), \tilde{\mathbf{P}}\text{-a.s.} \tag{4.2}$$

as well as

$$\mathrm{Law}_{\tilde{\mathbf{P}}}(\tilde{\Lambda}^k) = \mathrm{Law}_{\mathbf{P}^k_{x_0}}(\Lambda) \qquad \text{for all } k \in \mathbb{N}\cup\{\infty\}\,. \tag{4.3}$$

From (1.5) and (1.6), we see that $\tilde{\mathbf{P}}$-a.s. for all $k \in \mathbb{N}\cup\{\infty\}$, the process $\tilde{\Lambda}^k$ is $\mathbf{M}_1^+$-valued, and a particular consequence of (4.3) is that for each $\phi \in \mathcal{S}$ and all $k \in \mathbb{N}\cup\{\infty\}$, we have the identity

$$\mathbf{E}^{\mathbf{P}^{k}_{x_0}}\Big[\int_s^t g_u(\Lambda)\mathscr{A}_u(\Lambda_u)[\phi]\,\mathrm{d}u\Big] = \mathbf{E}^{\tilde{\mathbf{P}}}\Big[\int_s^t g_u(\tilde{\Lambda}^k)\mathscr{A}_u(\tilde{\Lambda}^k_u)[\phi]\,\mathrm{d}u\Big]\,.$$

So instead of (4.1) in Proposition 4.2, we can aim to prove

$$\lim_{k\to\infty} \mathbf{E}^{\tilde{\mathbf{P}}}\Big[\int_s^t g_u(\tilde{\Lambda}^k)\mathscr{A}_u(\tilde{\Lambda}^k_u)[\phi]\,\mathrm{d}u\Big] = \mathbf{E}^{\tilde{\mathbf{P}}}\Big[\int_s^t g_u(\tilde{\Lambda}^\infty)\mathscr{A}_u(\tilde{\Lambda}^\infty_u)[\phi]\,\mathrm{d}u\Big]\,. \tag{4.4}$$

This involves a single probability measure and an almost surely convergent sequence $(\tilde{\Lambda}^k)_{k\in\mathbb{N}\cup\{\infty\}}$ of random variables. Note that to establish (4.4), we no longer need global continuity of $\lambda \mapsto \int_s^t \mathscr{A}_u(\lambda_u)[\phi]\,\mathrm{d}u$ as a map $\mathbf{C}([0,T];\mathcal{S}') \to \mathbb{R}$, or rather continuity on a measurable set of full $\mathbf{P}^{\infty}_{x_0}$-measure with $\mathbf{P}^{\infty}_{x_0}$-negligible boundary, which would be sufficient by the Portemanteau theorem; we only need continuity for the family of convergent sequences $(\tilde{\Lambda}^k(\tilde{\omega}))_{k\in\mathbb{N}\cup\{\infty\}}$ for each $\tilde{\omega} \in \tilde{\Omega}_0$, where $\tilde{\Omega}_0 \in \tilde{\mathcal{F}}$ is some set of full $\tilde{\mathbf{P}}$-measure.

**Smoothing** To establish the above continuity property, we use convolutional smoothing. This allows us to approximate the non-smooth drift function $b$ with a regularized function; see Lemma 4.3 below. In more detail, this goes as follows.

Let $h \in \mathbf{C}_c^\infty(\mathbb{R}^d)$ be such that $h \geq 0$ and $\int_{\mathbb{R}^d} h(x)\,\mathrm{d}x = 1$, and let $h_\delta := \delta^d h(\delta x)$ for $\delta > 0$. Then $(h_\delta)_{\delta>0}$ is a compactly supported approximate identity; see Rudin [56, Def. 6.31]. We then define

$$b^\delta : [0,T] \times \mathbb{R}^d \times \mathbf{M}_1^+ \to \mathbb{R}^d$$

for each $\delta > 0$ by setting

$$b_t^\delta(x,\lambda) := \Big(b_t(\lambda) * h_\delta\Big)(x) := \int_{\mathbb{R}^d} b_t(y,\lambda) h_\delta(x-y)\,\mathrm{d}y\,, \tag{4.5}$$

where we denote by $b_t(\lambda)$ for $(t,\lambda) \in [0,T] \times \mathbf{M}_1^+$ the function $b_t(\,\cdot\,,\lambda) : \mathbb{R}^d \to \mathbb{R}^d$. By standard results, the integral in (4.5) is well defined because $b$ is bounded by Assumption 1.5; see e.g. Makarov and Podkorytov [45, Cor. 9.3.2].

For the next lemma, we recall from before Assumption 1.2 that $\mathcal{E}$ denotes the product-$\sigma$-algebra on $[0,T] \times \mathbb{R}^d \times \mathbf{M}_1^+$ generated by the Lebesgue-measurable sets of $[0,T] \times \mathbb{R}^d$ and the Borel-measurable sets of $\mathcal{P}_{\mathrm{wk}^*}$, the set of probability measures on $\mathcal{B}(\mathbb{R}^d)$ endowed with the narrow topology.

The function $b^\delta$ enjoys the following properties.

**Lemma 4.3** | *Let $b : [0,T] \times \mathbb{R}^d \times \mathbf{M}_1^+ \to \mathbb{R}^d$ satisfy Assumption 1.5 and fix $\delta > 0$. Then $b^\delta$ is $\mathcal{E}/\mathcal{B}(\mathbb{R}^d)$-measurable. Moreover:*

1) *For each $\lambda \in \mathbf{M}_1^+$ and $t \in [0,T]$, $x \mapsto b_t^\delta(x,\lambda)$ is continuous.*

2) *If $K \subseteq \mathbb{R}^d$ is compact and $1 \leq r' < \infty$, then for each $\lambda \in \mathbf{M}_1^+$ and $t \in [0,T]$, we have $\lim_{\delta\to 0} \|b_t(\lambda) - b_t^\delta(\lambda)\|_{\mathbb{L}^{r'}(K,\,\mathrm{d}x)} = 0$.*

3) *For each $\lambda \in \mathbf{M}_1^+$ and $t \in [0,T]$, $\|b_t^\delta(\lambda)\|_\infty \leq \|b_t(\lambda)\|_{\mathbb{L}^\infty(\mathbb{R}^d,\mathrm{d}x)}$.*

4) *Let $K_h := \mathrm{cl}_{(\mathbb{R}^d,|\,\cdot\,|)}\{x \in \mathbb{R}^d : |h(x)| > 0\}$, and $K \subseteq \mathbb{R}^d$ be compact. Then*

$$K' := K - K_h = \{x - x' : x \in K, x' \in K_h\}$$

*is compact, and for all $\lambda \in \mathbf{M}_1^+$, $t \in [0,T]$ and $r' \in [1,\infty]$, we have*

$$\|b_t^\delta(\lambda)\|_{\mathbb{L}^{r'}(K,\,\mathrm{d}x)} \leq \|b_t(\lambda)\|_{\mathbb{L}^{r'}(K',\mathrm{d}x)}\,.$$

***Proof*** See Appendix B. □

**Comparing terms** We now consider the *smoothed drift coefficient* $b_t^\delta(\,\cdot\,,\lambda)$ and define for each $\delta > 0$, $(t,\lambda) \in [0,T] \times \mathbf{M}_1^+$ and $\phi \in \mathcal{S}$ the function

$$\Big(\mathscr{L}_t^\delta(\lambda)\phi\Big)(x) := b_t^\delta(x,\lambda) \cdot \nabla\phi(x) + \frac{1}{2} a_t(x,\lambda) : \nabla^2\phi(x) \tag{4.6}$$

with $x \in \mathbb{R}^d$. This is entirely analogous to the definition of $\mathscr{L}_t(\lambda)\phi$ in (2.8), just replacing $b_t(x,\lambda)$ by $b_t^\delta(x,\lambda)$. Note in particular for later use that (2.8) and (4.6) give

$$\mathscr{L}_t(\lambda)\phi - \mathscr{L}_t^\delta(\lambda)\phi = \Big(b_t(\lambda) - b_t^\delta(\lambda)\Big) \cdot \nabla\phi. \tag{4.7}$$

We can now use the definition (2.26) of $\mathscr{A}_u(\lambda)[\phi] = \lambda[\mathscr{L}_u(\lambda)\phi]$ for $\lambda \in \mathbf{M}_1^+$ to write

$$\begin{aligned}
&\int_s^t g_u(\tilde\Lambda^\infty)\mathscr{A}_u(\tilde\Lambda_u^\infty)[\phi]\,\mathrm{d}u - \int_s^t g_u(\tilde\Lambda^k)\mathscr{A}_u(\tilde\Lambda_u^k)[\phi]\,\mathrm{d}u \\
&= \int_s^t \Big(g_u(\tilde\Lambda^\infty)\tilde\Lambda_u^\infty[\mathscr{L}_u(\tilde\Lambda_u^\infty)\phi] - g_u(\tilde\Lambda^k)\tilde\Lambda_u^k[\mathscr{L}_u(\tilde\Lambda_u^k)\phi]\Big)\,\mathrm{d}u \\
&=: \tilde A_\delta + \tilde B_{\delta,k} + \tilde C_{\delta,k} + \tilde D_k\,,
\end{aligned} \tag{4.8}$$

where

$$\tilde A_\delta := \int_s^t g_u(\tilde\Lambda^\infty)\Big(\tilde\Lambda_u^\infty[\mathscr{L}_u(\tilde\Lambda_u^\infty)\phi] - \tilde\Lambda_u^\infty[\mathscr{L}_u^\delta(\tilde\Lambda_u^\infty)\phi]\Big)\,\mathrm{d}u\,, \tag{4.9}$$

$$\tilde B_{\delta,k} := \int_s^t g_u(\tilde\Lambda^\infty)\Big(\tilde\Lambda_u^\infty[\mathscr{L}_u^\delta(\tilde\Lambda_u^\infty)\phi] - \tilde\Lambda_u^k[\mathscr{L}_u^\delta(\tilde\Lambda_u^\infty)\phi]\Big)\,\mathrm{d}u\,, \tag{4.10}$$

$$\tilde C_{\delta,k} := \int_s^t \Big(g_u(\tilde\Lambda^\infty)\tilde\Lambda_u^k[\mathscr{L}_u^\delta(\tilde\Lambda_u^\infty)\phi] - g_u(\tilde\Lambda^k)\tilde\Lambda_u^k[\mathscr{L}_u(\tilde\Lambda_u^\infty)\phi]\Big)\,\mathrm{d}u, \tag{4.11}$$

$$\tilde D_k := \int_s^t g_u(\tilde\Lambda^k)\Big(\tilde\Lambda_u^k[\mathscr{L}_u(\tilde\Lambda_u^\infty)\phi] - \tilde\Lambda_u^k[\mathscr{L}_u(\tilde\Lambda_u^k)\phi]\Big)\,\mathrm{d}u\,. \tag{4.12}$$

Let us also define

$$A_\delta := \mathbf{E}^{\tilde{\mathbf{P}}}[\tilde A_\delta]\,,\quad B_{\delta,k} := \mathbf{E}^{\tilde{\mathbf{P}}}[\tilde B_{\delta,k}]\,,\quad C_{\delta,k} := \mathbf{E}^{\tilde{\mathbf{P}}}[\tilde C_{\delta,k}]\,,\quad D_k := \mathbf{E}^{\tilde{\mathbf{P}}}[\tilde D_k]\,. \tag{4.13}$$

Via (4.8) and (4.13), we see that in order to establish (4.4), it is sufficient to show that

$$\lim_{\delta\to 0}|A_\delta| + \lim_{\delta\to 0}\lim_{k\to\infty}|B_{\delta,k}| + \lim_{\delta\to 0}\lim_{k\to\infty}|C_{\delta,k}| + \lim_{k\to\infty}|D_k| = 0\,.$$

## 4.1 Bounding the individual contributions

In this section, we establish four auxiliary results to handle the terms $A_\delta$, $B_{\delta,k}$, $C_{\delta,k}$, $D_k$ in (4.13). We start by controlling the term $A_\delta$.

**Lemma 4.4** | *Under Assumptions 1.2 and 1.5, fix $\phi \in \mathbf{C}_c^\infty(\mathbb{R}^d)$ and a bounded measurable $g : [0,T]\times\Omega \to \mathbb{R}$. Then for any $\varepsilon > 0$, there exists $\delta_A(\varepsilon) > 0$ such that $|A_\delta| < \varepsilon$ for all $\delta < \delta_A(\varepsilon)$. In other words, we have $\lim_{\delta\to 0}|A_\delta| = 0$.*

***Proof*** As $A_\delta = \mathbf{E}^{\tilde{\mathbf{P}}}[\tilde A_\delta]$, recall from (4.9) and (4.7) that

$$\begin{aligned}
\tilde A_\delta &= \int_s^t g_u(\tilde\Lambda^\infty)\tilde\Lambda_u^\infty\Big[\Big(\mathscr{L}_u(\tilde\Lambda_u^\infty) - \mathscr{L}_u^\delta(\tilde\Lambda_u^\infty)\Big)\phi\Big]\,\mathrm{d}u \\
&= \int_s^t g_u(\tilde\Lambda^\infty)\tilde\Lambda_u^\infty\Big[\Big(b_u(\tilde\Lambda_u^\infty) - b_u^\delta(\tilde\Lambda_u^\infty)\Big)\cdot\nabla\phi\Big]\,\mathrm{d}u\,.
\end{aligned} \tag{4.14}$$

Consider first the case $s > 0$. Since $\mathrm{Law}_{\tilde{\mathbf{P}}}(\tilde\Lambda^\infty) = \mathbf{P}_{x_0}^\infty = \mathrm{Law}_{\mathbf{P}_{x_0}^\infty}(\Lambda)$ by the Skorohod construction (4.3), the process $\tilde\Lambda^\infty$ under $\tilde{\mathbf{P}}$ has the same probabilistic properties as $\Lambda$ under $\mathbf{P}_{x_0}^\infty$. Therefore, by Theorem 3.1, we have $\tilde\Lambda_t \ll \mathrm{d}x$ ($\mathrm{d}t\otimes\mathrm{d}\tilde{\mathbf{P}}$)-a.e. In fact, we can find

$$1 < q < \infty \text{ and } 1 < r < \infty,\ 0 < \gamma < 1\,, \tag{4.15}$$

a constant $c_{\mathrm{Thm.\,3.1}} < \infty$ along with a measurable map $\tilde{p} : (0,T] \to \mathbb{L}^q((\tilde{\Omega}, \tilde{\mathfrak{F}}, \tilde{\mathbf{P}}); \mathbb{L}^r(\mathrm{d}x))$ which with $\tilde{p}_t(\,\cdot\,) := \tilde{p}(t)(\,\cdot\,)$ satisfies the bound

$$\Big\| \|\tilde{p}_t\|_{\mathbb{L}^r(\mathrm{d}x)} \Big\|_{\mathbb{L}^q(\tilde{\mathbf{P}})} \le c_{\mathrm{Thm.\,3.1}}(1 \wedge t)^{-\gamma} \quad \text{for all } t \in (0,T]\,, \tag{4.16}$$

and for which we have

$$\tilde{p}_t(\tilde{\omega}) = \mathrm{d}\tilde{\Lambda}^\infty_t(\tilde{\omega})/\mathrm{d}x \quad \text{in } \mathbb{L}^r(\mathrm{d}x)\ (\mathrm{d}t \otimes \mathrm{d}\tilde{\mathbf{P}})\text{-a.e.} \tag{4.17}$$

So from (4.17), we have $(\mathrm{d}t \otimes \mathrm{d}\tilde{\mathbf{P}})$-a.e. that

$$\tilde{\Lambda}^\infty_u\Big[\Big(b_u(\tilde{\Lambda}^\infty_u) - b^\delta_u(\tilde{\Lambda}^\infty_u)\Big) \cdot \nabla\phi\Big] = \int_{\mathbb{R}^d} \Big(b_u(x, \tilde{\Lambda}^\infty_u) - b^\delta_u(x, \tilde{\Lambda}^\infty_u)\Big) \cdot \nabla\phi(x)\, \tilde{p}_u(x)\, \mathrm{d}x\,. \tag{4.18}$$

To estimate $\tilde{A}_\delta$, we denote the support of the gradient of $\phi$ by

$$K_\phi := \mathrm{cl}_{(\mathbb{R}^d, |\,\cdot\,|)}\{x \in \mathbb{R}^d : |\nabla\phi(x)| > 0\}\,.$$

Then (4.14), (4.18) and the fact that $\nabla\phi$ is supported in $K_\phi$ yield with $c_{g,\phi} := \|g\|_\infty \|\phi\|^*_1$ that

$$\begin{aligned} |\tilde{A}_\delta| &= \Big| \int_s^t g_u(\tilde{\Lambda}^\infty)\tilde{\Lambda}^\infty_u\Big[\Big(b_u(\tilde{\Lambda}^\infty_u) - b^\delta_u(\tilde{\Lambda}^\infty_u)\Big) \cdot \nabla\phi\Big]\,\mathrm{d}u\Big| \\ &\le c_{g,\phi} \int_s^t \int_{\mathbb{R}^d} |b_u(x, \tilde{\Lambda}^\infty_u) - b^\delta_u(x, \tilde{\Lambda}^\infty_u)| \mathbb{1}_{K_\phi}(x) \tilde{p}_u(x)\,\mathrm{d}x\,\mathrm{d}u\,. \end{aligned} \tag{4.19}$$

Now take $r$ and $q$ as in (4.15) with conjugate exponents $r'$ and $q'$, i.e. satisfying the identity $1/r + 1/r' = 1 = 1/q + 1/q'$. Recalling that $A_\delta = \mathbf{E}^{\tilde{\mathbf{P}}}[\tilde{A}_\delta]$, taking expectations in (4.19) and using Hölder's inequality in the $\mathrm{d}x$-integral, then Tonelli's theorem to interchange the $\mathrm{d}u$-integral and $\mathbf{E}^{\tilde{\mathbf{P}}}$ and again Hölder's inequality for the expectation gives

$$\begin{aligned} |A_\delta| &\le c_{g,\phi} \int_s^t \mathbf{E}^{\tilde{\mathbf{P}}}[\|b^\delta_u(\tilde{\Lambda}^\infty_u) - b_u(\tilde{\Lambda}^\infty_u)\|_{\mathbb{L}^{r'}(K_\phi, \mathrm{d}x)} \|\tilde{p}_u\|_{\mathbb{L}^r(\mathbb{R}^d, \mathrm{d}x)}]\,\mathrm{d}u \\ &\le c_{g,\phi} \int_s^t \Big\| \|b^\delta_u(\tilde{\Lambda}^\infty_u) - b_u(\tilde{\Lambda}^\infty_u)\|_{\mathbb{L}^{r'}(K_\phi, \mathrm{d}x)} \Big\|_{\mathbb{L}^{q'}(\tilde{\mathbf{P}})} \Big\| \|\tilde{p}_u\|_{\mathbb{L}^r(\mathrm{d}x)} \Big\|_{\mathbb{L}^q(\tilde{\mathbf{P}})}\,\mathrm{d}u\,. \end{aligned}$$

But now we can use (4.16), so that with $c_{g,\phi,r,q,\gamma} := c_{g,\phi} c_{\mathrm{Thm.\,3.1}}$, we get the bound

$$|A_\delta| \le c_{g,\phi,r,q,\gamma} \int_s^t \Big\| \|b^\delta_u(\tilde{\Lambda}^\infty_u) - b_u(\tilde{\Lambda}^\infty_u)\|_{\mathbb{L}^{r'}(K_\phi, \mathrm{d}x)} \Big\|_{\mathbb{L}^{q'}(\tilde{\mathbf{P}})} (1 \wedge u)^{-\gamma}\,\mathrm{d}u.$$

Thus a final application of Hölder's inequality with conjugate exponents $1 < w, w' < \infty$, which we chose momentarily, yields

$$\begin{aligned} |A_\delta| &\le c_{g,\phi,r,q,\gamma} \Big( \int_s^t (1 \wedge u)^{-\gamma w}\,\mathrm{d}u \Big)^{\frac{1}{w}} \\ &\quad \times \Big( \int_s^t \Big\| \|b^\delta_u(\tilde{\Lambda}^\infty_u) - b_u(\tilde{\Lambda}^\infty_u)\|_{\mathbb{L}^{r'}(K_\phi, \mathrm{d}x)} \Big\|^{w'}_{\mathbb{L}^{q'}(\tilde{\mathbf{P}})}\,\mathrm{d}u \Big)^{\frac{1}{w'}}. \end{aligned} \tag{4.20}$$

Since $\gamma < 1$ by (4.15), we can find $w > 1$ small enough to have

$$\int_s^t (1 \wedge u)^{-\gamma w}\,\mathrm{d}u \le \int_0^T (1 \wedge u)^{-\gamma w}\,\mathrm{d}u =: c_{w,\gamma} < \infty.$$

With this, (4.20) then implies

$$|A_\delta| \le c_A \Big( \int_s^t \Big\| \|b_u^\delta(\tilde\Lambda_u^\infty) - b_u(\tilde\Lambda_u^\infty)\|_{\mathbb{L}^{r'}(K_\phi, \mathrm{d}x)} \Big\|_{\mathbb{L}^{q'}(\tilde{\mathbf{P}})}^{w'} \mathrm{d}u \Big)^{\frac{1}{w'}} \tag{4.21}$$

with a constant $c_A := c_{g,\phi,r,q,\gamma} c_{w,\gamma}$ which is independent of $s$ and $t$. From part 2) of Lemma 4.3, we have $\tilde{\mathbf{P}}$-a.s. that

$$\lim_{\delta\to 0} \|b_u^\delta(\tilde\Lambda_u^\infty) - b_u(\tilde\Lambda_u^\infty)\|_{\mathbb{L}^{r'}(K_\phi, \mathrm{d}x)} = 0 \quad \text{for all } u \in [0,T]. \tag{4.22}$$

As both $b$ and $b^\delta$ are bounded by Assumption 1.5 and part 3) of Lemma 4.3, dominated convergence, (4.21) and (4.22) allow us to conclude that for any $\varepsilon > 0$, there exists $\delta_A(\varepsilon)$ such that for all $\delta < \delta_A(\varepsilon)$, we have $|A_\delta| < \varepsilon$.

To deal with the case $s = 0$, consider $s' \in (0,t)$. Splitting the $\mathrm{d}u$-integral in $\tilde A_\delta$ into $\int_0^{s'} + \int_{s'}^t$ yields via (4.14) the estimate

$$\begin{aligned} |\tilde A_\delta| \le \|g\|_\infty \int_0^{s'} \Big| \tilde\Lambda_u^\infty \Big[ \Big( b_u(\tilde\Lambda_u^\infty) - b_u^\delta(\tilde\Lambda_u^\infty)\Big) \cdot \nabla\phi \Big] \Big| \, \mathrm{d}u \\ + \|g\|_\infty \int_{s'}^t \Big| \tilde\Lambda_u^\infty \Big[ \Big( b_u(\tilde\Lambda_u^\infty) - b_u^\delta(\tilde\Lambda_u^\infty)\Big) \cdot \nabla\phi \Big] \Big| \, \mathrm{d}u \, . \end{aligned}$$

Because $\|b_u^\delta(\tilde\Lambda_u^\infty)\|_\infty \le \|b_u(\tilde\Lambda_u^\infty)\|_{\mathbb{L}^\infty}$ by part 3) of Lemma 4.3 and because $\tilde\Lambda^\infty$ is $\tilde{\mathbf{P}}$-a.s. $\mathbf{M}_1^+$-valued, we get

$$|A_\delta| \le 2s'\|g\|_\infty \|b\|_\infty \|\phi\|_1^* + \|g\|_\infty \, \mathbf{E}^{\tilde{\mathbf{P}}} \Big[ \int_{s'}^t \Big| \tilde\Lambda_u^\infty \Big[ \Big( b_u^\delta(\tilde\Lambda_u^\infty) - b_u(\tilde\Lambda_u^\infty)\Big) \cdot \nabla\phi \Big] \Big| \, \mathrm{d}u \Big].$$

Both terms on the right-hand side can be made smaller than $\varepsilon/2$—the first by choosing $s'$ small enough and the second by our estimate for the case $s' > 0$. So again we find for any $\varepsilon > 0$ some $\delta_A(\varepsilon)$ such that for all $\delta < \delta_A(\varepsilon)$, we have $|A_\delta| < \varepsilon$. This completes the proof. ☐

We now turn our attention to the term $B_{\delta,k}$.

**Lemma 4.5** | *Under Assumptions 1.2 and 1.5, fix $\phi \in \mathbf{C}_c^\infty(\mathbb{R}^d)$ and a bounded measurable $g : [0,T] \times \Omega \to \mathbb{R}$. Then for any $\varepsilon > 0$ and $\delta > 0$, there exists $k_B(\delta,\varepsilon) > 0$ such that $|B_{\delta,k}| < \varepsilon$ for all $k > k_B(\varepsilon,\delta)$. In other words, $\lim_{k\to\infty} |B_{\delta,k}| = 0$ for all $\delta > 0$.*

***Proof*** The definition of $\tilde B_{\delta,k}$ in (4.10) gives

$$\begin{aligned} |B_{\delta,k}| = |\,\mathbf{E}^{\tilde{\mathbf{P}}}[\tilde B_{\delta,k}]| &= \Big| \mathbf{E}^{\tilde{\mathbf{P}}} \Big[ \int_s^t g_u(\tilde\Lambda^\infty)(\tilde\Lambda_u^\infty - \tilde\Lambda_u^k)[\mathscr{L}_u^\delta(\tilde\Lambda_u^\infty)\phi] \, \mathrm{d}u \Big] \Big| \\ &\le \|g\|_\infty \mathbf{E}^{\tilde{\mathbf{P}}} \Big[ \int_s^t \Big| (\tilde\Lambda_u^k - \tilde\Lambda_u^\infty)[\mathscr{L}_u^\delta(\tilde\Lambda_u^\infty)\phi] \Big| \, \mathrm{d}u \Big]. \end{aligned} \tag{4.23}$$

Fix a measurable subset $\tilde\Omega_0 \subseteq \tilde\Omega$ with $\tilde{\mathbf{P}}[\tilde\Omega_0] = 1$ such that $\tilde\Lambda^k(\tilde\omega)$ is $\mathbf{M}_1^+$-valued for all $k \in \mathbb{N} \cup \{\infty\}$ and $\tilde\omega \in \tilde\Omega_0$. Using parts 1) and 3) of Lemma 4.3 together with the fact that $\phi \in \mathbf{C}_c^\infty$ and $b$ is bounded by Assumption 1.5, we see that $x \mapsto b_u^\delta(x, \tilde\Lambda_u^\infty(\tilde\omega)) \cdot \nabla\phi(x)$ is continuous for each $(u,\tilde\omega) \in [s,t] \times \tilde\Omega_0$ and bounded uniformly in $(u,x,\tilde\omega) \in [s,t] \times \mathbb{R}^d \times \tilde\Omega_0$.

By Assumption 1.2, the function $x \mapsto a_u(x, \tilde{\Lambda}^\infty_u(\tilde{\omega})) : \nabla^2\phi(x)$ satisfies the same properties. From its definition in (4.6), we see that the same is true for the map $x \mapsto \mathscr{L}^\delta_u(\tilde{\Lambda}^\infty_u(\tilde{\omega}))\phi(x)$. Since $\lim_{k\to\infty} \tilde{\Lambda}^k(\tilde{\omega}) = \tilde{\Lambda}^\infty(\tilde{\omega})$ in $\mathbf{C}([0,T]; \mathscr{H}_{-p})$ and $\tilde{\Lambda}^\infty(\tilde{\omega})$ is $\mathbf{M}_1^+$-valued for each $\tilde{\omega} \in \tilde{\Omega}_0$, Lemma A.3 applied for each $(u, \tilde{\omega}) \in [0,T] \times \tilde{\Omega}_0$ shows that $\lim_{k\to\infty} \tilde{\Lambda}^k_u(\tilde{\omega}) = \tilde{\Lambda}^\infty_u(\tilde{\omega})$ in $\mathcal{P}_{\mathrm{wk}^*}$, i.e. we have convergence in the narrow topology. So

$$\lim_{k\to\infty} \Big|(\tilde{\Lambda}^k_u - \tilde{\Lambda}^\infty_u)[\mathscr{L}^\delta_u(\tilde{\Lambda}^\infty_u)\phi]\Big| = 0 \quad \text{for all } u \in [s,t],\ \tilde{\mathbf{P}}\text{-a.s.},$$

and using dominated convergence, (4.23) shows that for all $\delta > 0$, we can find $k_B(\delta, \varepsilon)$ such that $k > k_B(\delta, \varepsilon)$ implies $|B_{\delta,k}| < \varepsilon$. □

Obtaining control over $C_{\delta,k} = \mathbf{E}^{\tilde{\mathbf{P}}}[\tilde{C}_{\delta,k}]$ is more intricate than for $A_\delta$, $B_{\delta,k}$ or $D_{\delta,k}$ below.

**Lemma 4.6** | *Under Assumptions 1.2 and 1.5, fix $\phi \in \mathbf{C}^\infty_c(\mathbb{R}^d)$ and a bounded measurable $g : [0,T] \times \Omega \to \mathbb{R}$ such that for each $u \in [0,T]$, the real-valued function $g_u(\,\cdot\,) = g(u, \,\cdot\,)$ is continuous. Then for any $\varepsilon > 0$, there exist $\delta_C(\varepsilon) > 0$ and $k_C(\varepsilon) > 0$ such that $|C_{k,\delta}| < \varepsilon$ for all $0 < \delta < \delta_C(\varepsilon)$ and $k > k_C(\varepsilon)$. In other words, we have* $\lim_{\delta\to 0} \lim_{k\to\infty} |C_{\delta,k}| = 0$.

To control the term $C_{\delta,k}$, we use the regularity properties of the finite particle system $\boldsymbol{X}^n$ in (1.1), (1.2), which we must established in (3.3). The estimate for $\tilde{C}_{\delta,k}$ is complicated by the fact that it involves both $\tilde{\Lambda}^k$ and $\tilde{\Lambda}^\infty$ and both the smoothed and original generators $\mathscr{L}^\delta$ and $\mathscr{L}$. More specifically, in (4.24) and (4.27) below, we encounter a term of the form

$$\int_s^t g_u(\tilde{\Lambda}^k)\tilde{\Lambda}^k_u\Big[\Big(b^\delta_u(\tilde{\Lambda}^\infty_u) - b_u(\tilde{\Lambda}^\infty_u)\Big) \cdot \nabla\phi\Big]\,\mathrm{d}u \tag{4.27}$$

in which the integral of the random measurable function $x \mapsto b^\delta_u(x, \tilde{\Lambda}^\infty_u) - b_u(x, \tilde{\Lambda}^\infty_u)$ against the random and purely atomic measure $\tilde{\Lambda}^k_u$ appears. The term $\tilde{C}_{\delta,k}$ thus involves the action of an *irregular measure on an irregular function*, and controlling this pairing is not straightforward. Moreover, as a consequence of the Skorohod construction, $\tilde{\Lambda}^k_u$ and $\tilde{\Lambda}^\infty_u$ are coupled and necessarily dependent, which severely complicates matters since both appear in the pairing just mentioned.

Our strategy in the proof of Lemma 4.6 is to first use a compactness argument to get rid of the randomness stemming from $\tilde{\Lambda}^\infty_u$, and this exploits Assumption 1.5. The foundation for this strategy is laid in Step 2.2 of the proof and culminates in the bound (4.34) in Step 2.4. It lets us control $\tilde{C}_{\delta,k}$, or more precisely (4.27), by a finite number of terms of the form $\tilde{\Lambda}^k_u[|(b^\delta_u(\lambda) - b_u(\lambda)) \cdot \nabla\phi|]$, where $\lambda \in \mathbf{M}_1^+$ is now a *fixed* probability measure. Since this is independent of $\tilde{\Lambda}^\infty_u$, we can then in Step 2.5 use the fact that $\tilde{\Lambda}^k_u$ mimics the empirical measure of the *exchangeable* $k$-particle McKean–Vlasov system. Indeed, upon taking expectations, we can reduce our problem to bounding terms involving only the law of a single particle, namely $\mathbf{E}^{\mathbb{P}^k_{x_0}}[|(b^\delta_u(X^{1,k}_u, \lambda) - b_u(X^{1,k}_u, \lambda)) \cdot \nabla\phi(X^{1,k}_u)|]$; see (4.36) below. Thanks to the work we put into Lemma 3.3, we already know that the law of $X^{1,k}_u$ possesses a density with respect to Lebesgue measure, whose $\mathbb{L}^r(\mathrm{d}x)$-norm we can control uniformly over $k \in \mathbb{N}$. In other words, with this reduction, we find ourselves in (4.36) in a situation where we have the action of a *regular measure on an irregular function*, and controlling this pairing is straightforward. Indeed, using Hölder's inequality, we can finally use the fact that $\|(b^\delta_u(\lambda) - b_u(\lambda)) \cdot \nabla\phi\|_{\mathbb{L}^{r'}(\mathrm{d}x)}$ vanishes as $\delta \to 0$; see (4.38) below.

***Proof of Lemma 4.6*** Take $\phi \in \mathbf{C}_c^\infty(\mathbb{R}^d)$ and a measurable bounded $g : [0,T] \times \Omega \to \mathbb{R}$ with the property that for each $u \in [0,T]$, the real-valued function $g_u(\,\cdot\,) \coloneqq g(u,\,\cdot\,)$ on $\Omega$ is continuous. We begin by using the definition of $\tilde{C}_{\delta,k}$ from (4.11) to obtain

$$\begin{aligned}\tilde{C}_{\delta,k} &= \int_s^t \Big(g_u(\tilde{\Lambda}^\infty)\tilde{\Lambda}_u^k[\mathscr{L}_u^\delta(\tilde{\Lambda}_u^\infty)\phi] - g_u(\tilde{\Lambda}^k)\tilde{\Lambda}_u^k[\mathscr{L}_u(\tilde{\Lambda}_u^\infty)\phi]\Big)\,\mathrm{d}u \\ &= \int_s^t \Big(g_u(\tilde{\Lambda}^\infty) - g_u(\tilde{\Lambda}^k)\Big)\tilde{\Lambda}_u^k[\mathscr{L}_u^\delta(\tilde{\Lambda}_u^\infty)\phi]\,\mathrm{d}u \\ &\quad + \int_s^t g_u(\tilde{\Lambda}^k)\Big(\tilde{\Lambda}_u^k[\mathscr{L}_u^\delta(\tilde{\Lambda}_u^\infty)\phi] - \tilde{\Lambda}_u^k[\mathscr{L}_u(\tilde{\Lambda}_u^\infty)\phi]\Big)\,\mathrm{d}u \\ &=: \tilde{C}_{\delta,k,1} + \tilde{C}_{\delta,k,2}\,. \end{aligned} \tag{4.24}$$

Steps 1 and 2 below consider the last two terms in turn. In Step 3, we then combine our findings to complete the proof.

***Step 1*** For the term $\tilde{C}_{\delta,k,1}$, we copy its definition from (4.24), namely

$$\tilde{C}_{\delta,k,1} = \int_s^t \Big(g_u(\tilde{\Lambda}^\infty) - g_u(\tilde{\Lambda}^k)\Big)\tilde{\Lambda}_u^k[\mathscr{L}_u^\delta(\tilde{\Lambda}_u^\infty)\phi]\,\mathrm{d}u\,.$$

Assumptions 1.2 and 1.5 allow us to use the estimate for $|\tilde{\Lambda}_t^k[\mathscr{L}_u^\delta(\tilde{\Lambda}_u^\infty)\phi]|$ from (2.32) in Remark 2.7 to get

$$|\mathbf{E}^{\tilde{\mathbf{P}}}[\tilde{C}_{\delta,k,1}]| \le c_{b,\sigma,\bar{\sigma}}\|\phi\|_2^*\,\mathbf{E}^{\tilde{\mathbf{P}}}\Big[\int_s^t |g_u(\tilde{\Lambda}^\infty) - g_u(\tilde{\Lambda}^k)|\,\mathrm{d}u\Big]\,, \tag{4.25}$$

where we also used that by part 3) of Lemma 4.3, we have

$$c_{b^\delta,\sigma,\bar{\sigma}} = \|b^\delta\|_\infty + \frac{1}{2}(\|\sigma\|_\infty^2 + \|\bar{\sigma}\|_\infty^2) \le c_{b,\sigma,\bar{\sigma}}\,.$$

Because $g_u \in \mathbf{C}_\mathrm{b}(\Omega)$ and $\Omega = \mathbf{C}([0,T];\mathcal{S}')$, the continuity of the embedding $\mathscr{H}_{-p} \hookrightarrow \mathcal{S}'$ in (A.9) implies that for each $u \in [s,t]$, the function $g_u$ on $\mathbf{C}([0,T];\mathscr{H}_{-p})$ is continuous and bounded. By the Skorohod construction (4.2), we have $\tilde{\mathbf{P}}$-a.s. that $\tilde{\Lambda}^k \to \tilde{\Lambda}^\infty$ in $\mathbf{C}([0,T];\mathscr{H}_{-p})$ as $k \to \infty$, and so we get that $\tilde{\mathbf{P}}$-a.s.,

$$\lim_{k\to\infty} g_u(\tilde{\Lambda}^k) = g_u(\tilde{\Lambda}^\infty) \quad \text{for all } u \in [0,T].$$

Now fix $\varepsilon > 0$. Since $g$ is bounded, we can apply dominated convergence in (4.25) to get

$$|\mathbf{E}^{\tilde{\mathbf{P}}}[\tilde{C}_{k,\delta,1}]| < \varepsilon/2 \tag{4.26}$$

for $k > k_C(\varepsilon)$ large enough, independently of $\delta > 0$.

***Step 2*** For the term $\tilde{C}_{\delta,k,2}$ and its expectation, more work is required. We first copy its definition from (4.24) and then use (4.7) to obtain

$$\begin{aligned}\tilde{C}_{\delta,k,2} &= \int_s^t g_u(\tilde{\Lambda}^k)\Big(\tilde{\Lambda}_u^k[\mathscr{L}_u^\delta(\tilde{\Lambda}_u^\infty)\phi] - \tilde{\Lambda}_u^k[\mathscr{L}_u(\tilde{\Lambda}_u^\infty)\phi]\Big)\,\mathrm{d}u \\ &= \int_s^t g_u(\tilde{\Lambda}^k)\tilde{\Lambda}_u^k\Big[\Big(b_u^\delta(\tilde{\Lambda}_u^\infty) - b_u(\tilde{\Lambda}_u^\infty)\Big)\cdot\nabla\phi\Big]\,\mathrm{d}u\,. \end{aligned} \tag{4.27}$$

Our aim is to establish for $\mathbf{E}^{\tilde{\mathbf{P}}}[\tilde{C}_{\delta,k,2}]$ a bound of $\varepsilon/2$. This requires some intermediate steps, with the key steps being 2.2 and 2.4 below. Step 2.2 contains a basic reduction using a compactness argument. It forms the basis for Step 2.4 in which we obtain in (4.34) a bound on $\mathbf{E}^{\tilde{\mathbf{P}}}[\tilde{C}_{\delta,k,2}]$ that is independent of $\tilde{\Lambda}^\infty$.

***Step 2.1*** As in the proof of Lemma 4.4 and part 4) of Lemma 4.3, we recall $h$ from before (4.5) and define

$$K_\phi := \mathrm{cl}_{(\mathbb{R}^d, |\cdot|)}\{x \in \mathbb{R}^d : |\nabla\phi(x)| > 0\}\,,$$
$$K_h := \mathrm{cl}_{(\mathbb{R}^d, |\cdot|)}\{x \in \mathbb{R}^d : |h(x)| > 0\}$$

and set

$$K'_\phi := K_\phi - K_h\,.$$

Since $0 \in K_h$, we have $K_\phi \subseteq K'_\phi$. With this notation, we note that for each pair $\lambda, \lambda' \in \mathbf{M}_1^+$ the function $[s,t] \ni u \mapsto \sup_{x\in K'_\phi} |b_u(x,\lambda) - b_u(x,\lambda')| \in \mathbb{R}$ is analytically measurable; see Bertsekas and Shreve [7, Sec. 7.6], and thus universally measurable; see [7, Cor. 7.42.1]. For universally measurable functions, the Lebesgue integral relative to $\mathrm{d}u$ is well defined over sets of finite measure, thus in particular for our compact set $[s,t]$; see [7, Def. 7.19 and Sec. 7.7]. The next definition therefore makes sense.

Fix $\varepsilon_0 > 0$, to be chosen in Step 2.6 below. Recall that $\mathcal{P}_{\mathrm{wk}^*}(\mathbb{R}^d)$ denotes the space of probability measures on $\mathcal{B}(\mathbb{R}^d)$ endowed with narrow convergence. Write $\mathrm{P} := \mathbf{C}([0,T]; \mathcal{P}_{\mathrm{wk}^*}(\mathbb{R}^d))$ and $\mathrm{P}_{\mathrm{ac}}$ for the set of $\lambda \in \mathrm{P}$ which satisfy $\lambda_t \ll \mathrm{d}x$ for a.e. $t \in [0,T]$. For $\lambda \in \mathrm{P}_{\mathrm{ac}}$, we define the set

$$U_\lambda := \left\{\lambda' \in \mathrm{P} \;:\; \int_s^t \sup_{x\in K'_\phi} |b_u(x,\lambda_u) - b_u(x,\lambda'_u)|\,\mathrm{d}u < \varepsilon_0\right\}. \tag{4.28}$$

For $\lambda \in \mathrm{P}_{\mathrm{ac}}$, we get for a.e. $u \in [s,t]$ from Assumption 1.5 that $\mathrm{P} \ni \lambda' \mapsto b_u(x,\lambda'_u) \in \mathbb{R}$ is narrowly continuous at $\lambda$ uniformly with respect to $x \in K_\phi$. It follows that the set $U_\lambda$ is open in P.

We next define

$$O := \bigcup_{\lambda\in\mathrm{P}_{\mathrm{ac}}} U_\lambda \subseteq \mathrm{P}\,.$$

Since unions of open sets are open, $O$ is open and thus a Borel set of P.

***Step 2.2*** Next let $\varepsilon_1 > 0$, which like $\varepsilon_0$ from Step 2.1 will be chosen in Step 2.6 below. We claim that we can find a compact set $K \subseteq O$ such that

$$\tilde{\mathbf{P}}[\tilde{\Lambda}^\infty \in K^c] < \varepsilon_1\,. \tag{4.29}$$

Indeed, $\mathrm{P} = \mathbf{C}([0,T]; \mathcal{P}_{\mathrm{wk}^*}(\mathbb{R}^d))$ is Polish because $[0,T]$ is compact and $\mathcal{P}_{\mathrm{wk}^*}(\mathbb{R}^d)$ is separable, and so every finite Borel measure on P is inner regular. Hence to establish (4.29), it is enough to show that

$$\tilde{\mathbf{P}}[\tilde{\Lambda}^\infty \in O] = 1\,. \tag{4.30}$$

To prove (4.30), we first show that the measure $\tilde{\mathbf{P}} \circ (\tilde{\Lambda}^\infty)^{-1}$ is concentrated on P, i.e., that $\tilde{\Lambda}^\infty \in \mathrm{P}$ $\tilde{\mathbf{P}}$-a.s. To prove this, fix a measurable set $\tilde{\Omega}_0 \subseteq \tilde{\Omega}$ with $\tilde{\mathbf{P}}[\tilde{\Omega}_0] = 1$ and such that $\tilde{\Lambda}^\infty(\tilde{\omega})$ is $\mathbf{M}_1^+$-valued for all $\tilde{\omega} \in \tilde{\Omega}_0$. This is possible by (1.6) and (4.2). Since $\tilde{\Lambda}^\infty \in \mathbf{C}([0,T]; \mathscr{H}_{-p})$ by the Skorohod construction, we have for any sequence $(t_j)_{j\in\mathbb{N}}$ in $[0,T]$ with $t_j \to t_\infty$ as $j \to \infty$ that $\tilde{\Lambda}^\infty_{t_j}(\tilde{\omega}) \to \tilde{\Lambda}^\infty_{t_\infty}(\tilde{\omega})$ in $\mathscr{H}_{-p}$ as $j \to \infty$, and since $\tilde{\Lambda}^\infty_{t_\infty}(\tilde{\omega}) \in \mathbf{M}_1^+$, Lemma A.3 implies that $\tilde{\Lambda}^\infty_{t_j}(\tilde{\omega}) \to \tilde{\Lambda}^\infty_{t_\infty}(\tilde{\omega})$ in $\mathcal{P}_{\mathrm{wk}^*}$ as $j \to \infty$. This shows for all $\tilde{\omega} \in \tilde{\Omega}_0$ that $\tilde{\Lambda}^\infty(\tilde{\omega}) \in \mathrm{P}$. With this observation and in view of the Skorohod construction (4.3), we obtain (4.30) from (3.2) in Theorem 3.1 from which we get that

$$\tilde{\mathbf{P}}\Big[\tilde{\Lambda}^\infty_t \ll \mathrm{d}x \text{ for a.e. } t \in (0,T]\Big] = 1\,.$$

Note now that by construction the collection $(U_\lambda)_{\lambda\in\mathrm{P_{ac}}}$ forms an open cover of the compact set $K$, and so we can choose $\lambda^1, \dots, \lambda^m \in \mathrm{P_{ac}}$ such that $U_{\lambda^1}, \dots, U_{\lambda^m}$ cover $K$. Note that $m$ only depends on our fixed $\phi$, $\varepsilon_0$ from (4.28) and $\varepsilon_1$ in (4.29). We now inductively define

$$\begin{aligned} P_1 &\coloneqq U_{\lambda^1} \cap K\,, \\ P_\ell &\coloneqq (U_{\lambda^\ell} \cap K) \setminus \bigg( \bigcup_{1\le j<\ell} P_j \bigg) \qquad \text{for } 1 < \ell \le m\,. \end{aligned}$$

Since the sets $U_{\lambda^\ell}$ for $\ell \in [m]$ are open and $K$ is measurable, $(P_\ell)_{\ell\in[m]}$ is a partition of $K$ into pairwise disjoint sets in $\mathcal{B}(\mathrm{P})$; so

$$K = \bigcup_{\ell=1}^m P_\ell\,. \tag{4.31}$$

***Step 2.3*** We now consider $\tilde{C}_{\delta,k,2}$ in (4.27) on the set $K$. Using that $K$ equals by (4.31) the disjoint union of $(P_\ell)_{\ell\in[m]}$, we can write

$$\begin{aligned} &\mathbb{1}_K(\tilde{\Lambda}^\infty) \int_s^t \tilde{\Lambda}^k_u \Big[\Big|\Big(b^\delta_u(\tilde{\Lambda}^\infty_u) - b_u(\tilde{\Lambda}^\infty_u)\Big) \cdot \nabla\phi\Big|\Big]\,\mathrm{d}u \\ &= \sum_{\ell=1}^m \mathbb{1}_{P_\ell}(\tilde{\Lambda}^\infty) \int_s^t \tilde{\Lambda}^k_u \Big[\Big|\Big(b^\delta_u(\tilde{\Lambda}^\infty_u) - b_u(\tilde{\Lambda}^\infty_u)\Big) \cdot \nabla\phi\Big|\Big]\,\mathrm{d}u \\ &\le \sum_{\ell=1}^m \bigg( \mathbb{1}_{P_\ell}(\tilde{\Lambda}^\infty) \int_s^t \tilde{\Lambda}^k_u \Big[\Big|\Big(b^\delta_u(\tilde{\Lambda}^\infty_u) - b^\delta_u(\lambda^\ell_u)\Big) \cdot \nabla\phi\Big|\Big]\,\mathrm{d}u \\ &\qquad + \mathbb{1}_{P_\ell}(\tilde{\Lambda}^\infty) \int_s^t \tilde{\Lambda}^k_u \Big[\Big|\Big(b^\delta_u(\lambda^\ell_u) - b_u(\lambda^\ell_u)\Big) \cdot \nabla\phi\Big|\Big]\,\mathrm{d}u \\ &\qquad + \mathbb{1}_{P_\ell}(\tilde{\Lambda}^\infty) \int_s^t \tilde{\Lambda}^k_u \Big[\Big|\Big(b_u(\lambda^\ell_u) - b_u(\tilde{\Lambda}^\infty_u)\Big) \cdot \nabla\phi\Big|\Big]\,\mathrm{d}u \bigg)\,. \end{aligned} \tag{4.32}$$

For each of the terms indexed by $\ell \in [m]$ in the third summand in (4.32), we have from $P_\ell \subseteq U_{\lambda^\ell}$ and the definition of $U_{\lambda^\ell} \subseteq K$ in (4.31) that $\tilde{\mathbf{P}}$-a.s.,

$$\mathbb{1}_{P_\ell}(\tilde{\Lambda}^\infty) \int_s^t \tilde{\Lambda}^k_u \Big[\Big|\Big(b_u(\lambda^\ell_u) - b_u(\tilde{\Lambda}^\infty_u)\Big) \cdot \nabla\phi\Big|\Big]\,\mathrm{d}u \le \mathbb{1}_{P_\ell}(\tilde{\Lambda}^\infty)\|\phi\|_1^*\varepsilon_0\,.$$

This bound can be recycled to also control each of the terms in the first summand in (4.32), as follows. Since convolving with an approximate identity as in (4.5) is a linear operation, we have for $\lambda, \lambda' \in K$ that

$$b^\delta_u(x,\lambda_u) - b^\delta_u(x,\tilde{\lambda}'_u) = \Big(b_u(x,\lambda_u) - b_u(x,\lambda'_u)\Big)^\delta\,.$$

Moreover, our approximate identity is assumed to be compactly supported; so parts 1) and 3) of Lemma 4.3 show for all $u \in [s,t]$ that

$$\begin{aligned}\sup_{x\in K_\phi} |b^\delta_u(x,\lambda_u) - b^\delta_u(x,\lambda'_u)| &= \|b^\delta_u(\lambda_u) - b^\delta_u(\lambda'_u)\|_{\mathbb{L}^\infty(K_\phi,\mathrm{d}x)}\\ &\le \|b_u(\lambda_u) - b_u(\lambda'_u)\|_{\mathbb{L}^\infty(K'_\phi,\mathrm{d}x)}\\ &\le \sup_{x\in K'_\phi} |b_u(x,\lambda_u) - b_u(x,\lambda'_u)|\,.\end{aligned}$$

For each $\delta > 0$, we thus obtain for all $\ell \in [m]$ the $\tilde{\mathbf{P}}$-a.s. bound

$$\mathbb{1}_{P_\ell}(\tilde\Lambda^\infty)\int_s^t \tilde\Lambda^k_u\Big[\Big|\Big(b^\delta_u(\lambda^\ell_u) - b^\delta_u(\tilde\Lambda^\infty_u)\Big)\cdot\nabla\phi\Big|\Big]\,\mathrm{d}u \le \mathbb{1}_{P_\ell}(\tilde\Lambda^\infty)\|\phi\|^*_1\varepsilon_0\,. \tag{4.33}$$

Finally, using $\mathbb{1}_{P_\ell}(\tilde\Lambda^\infty) \le \mathbb{1}_K(\tilde\Lambda^\infty)$, each of the terms in the second summand in (4.32) can be trivially majorized by

$$\mathbb{1}_{P_\ell}(\tilde\Lambda^\infty)\int_s^t \tilde\Lambda^k_u\Big[\Big|\Big(b^\delta_u(\lambda^\ell_u) - b_u(\lambda^\ell_u)\Big)\cdot\nabla\phi\Big|\Big]\,\mathrm{d}u \le \mathbb{1}_K(\tilde\Lambda^\infty)\int_s^t \tilde\Lambda^k_u\Big[\Big|\Big(b^\delta_u(\lambda^\ell_u) - b_u(\lambda^\ell_u)\Big)\cdot\nabla\phi\Big|\Big]\,\mathrm{d}u.$$

***Step 2.4*** On the set $K$, we can thus control $\tilde C_{\delta,k,2}$ from (4.27) by combining the estimates from Step 2.3 via (4.32), and by using that $\mathbb{1}_K = \sum_{\ell=1}^m \mathbb{1}_{P_\ell}$ and $\mathbb{1}_{P_\ell} \le \mathbb{1}_K$, which gives the bound

$$\begin{aligned}|\tilde C_{\delta,k,2}\mathbb{1}_K(\tilde\Lambda^\infty)| &= \Big|\mathbb{1}_K(\tilde\Lambda^\infty)\int_s^t g_u(\tilde\Lambda^\infty)\tilde\Lambda^k_u\Big[\Big|\Big(b^\delta_u(\tilde\Lambda^\infty_u) - b_u(\tilde\Lambda^\infty_u)\Big)\cdot\nabla\phi\Big|\Big]\,\mathrm{d}u\Big|\\ &\le \|g\|_\infty\mathbb{1}_K(\tilde\Lambda^\infty)\|\phi\|^*_1\varepsilon_0\\ &\quad + \|g\|_\infty\mathbb{1}_K(\tilde\Lambda^\infty)\bigg(\sum_{\ell=1}^m\int_s^t \tilde\Lambda^k_u\Big[\Big|\Big(b^\delta_u(\lambda^\ell_u) - b_u(\lambda^\ell_u)\Big)\cdot\nabla\phi\Big|\Big]\,\mathrm{d}u\bigg)\\ &\quad + \|g\|_\infty\mathbb{1}_K(\tilde\Lambda^\infty)\|\phi\|^*_1\varepsilon_0\,.\end{aligned}$$

On the set $K^c$, on the other hand, using that $b$ and $b^\delta$ are bounded by Assumption 1.5 and appealing to part 3) of Lemma 4.3 gives

$$|\tilde C_{\delta,k,2}\mathbb{1}_{K^c}(\tilde\Lambda^\infty)| \le 2\|g\|_\infty\|b\|_\infty\|\phi\|^*_1(t-s)\mathbb{1}_{K^c}(\tilde\Lambda^\infty).$$

Upon taking expectations, we can combine the bounds on the sets $K$ and $K^c$ by using on $K$ the trivial bound $\mathbb{1}_K \le 1$ and on $K^c$ the construction (4.29) by which $\mathbf{E}^{\tilde{\mathbf{P}}}[\mathbb{1}_{K^c}] \le \varepsilon_1$, to find

$$\begin{aligned}C_{\delta,k,2} := |\,\mathbf{E}^{\tilde{\mathbf{P}}}[\tilde C_{\delta,k,2}]| \le \|g\|_\infty\sum_{\ell=1}^m\int_s^t \mathbf{E}^{\tilde{\mathbf{P}}}\Big[\tilde\Lambda^k_u\Big[\Big|\Big(b^\delta_u(\lambda^\ell_u) - b_u(\lambda^\ell_u)\Big)\cdot\nabla\phi\Big|\Big]\Big]\,\mathrm{d}u\\ + 2\|g\|_\infty\|\phi\|^*_1\varepsilon_0 + 2\|g\|_\infty\|b\|_\infty\|\phi\|^*_1(t-s)\varepsilon_1\,.\end{aligned} \tag{4.34}$$

The right-hand side of (4.34) is independent of $\tilde\Lambda^\infty$. This is the bound that we were after, and its usefulness is seen as follows.

***Step 2.5*** We consider first the case $s > 0$, and defer the case $s = 0$ to Step 2.7 below. Since $\mathrm{Law}_{\tilde{\mathbf{P}}}(\tilde{\Lambda}^k) = \mathrm{Law}_{\mathbf{P}^k_{x_0}}(\Lambda)$ by the Skorohod construction in (4.3) and $\mathbf{P}^k_{x_0} = \mathbb{P}^{n_k}_{x_0} \circ (\mu^{k_n})^{-1}$ by the definition (1.5), we see that each of the $m$ summands involving the expectation with respect to $\tilde{\mathbf{P}}$ in (4.34) satisfies

$$\begin{aligned} &\int_s^t \mathbf{E}^{\tilde{\mathbf{P}}}\left[\tilde{\Lambda}^k_u\Big[\Big|\Big(b_u(\lambda^\ell_u) - b^\delta_u(\lambda^\ell_u)\Big)\cdot\nabla\phi\Big|\Big]\right]\mathrm{d}u \\ &= \int_s^t \mathbf{E}^{\mathbf{P}^k_{x_0}}\left[\Lambda_u\Big[\Big|\Big(b_u(\lambda^\ell_u) - b^\delta_u(\lambda^\ell_u)\Big)\cdot\nabla\phi\Big|\Big]\right]\mathrm{d}u \\ &= \int_s^t \mathbf{E}^{\mathbb{P}^{n_k}_{x_0}}\left[\mu^{n_k}_u\Big[\Big|\Big(b_u(\lambda^\ell_u) - b^\delta_u(\lambda^\ell_u)\Big)\cdot\nabla\phi\Big|\Big]\right]\mathrm{d}u\,. \end{aligned} \tag{4.35}$$

The exchangeability of the particle system shown in Lemma 1.3 implies for the last equality in (4.35) that

$$\begin{aligned} &\int_s^t \mathbf{E}^{\mathbb{P}^{n_k}_{x_0}}\left[\mu^{n_k}_u\Big[\Big|\Big(b_u(\lambda^\ell_u) - b^\delta_u(\lambda^\ell_u)\Big)\cdot\nabla\phi\Big|\Big]\right]\mathrm{d}u \\ &= \frac{1}{n_k}\sum_{i=1}^{n_k}\int_s^t \mathbf{E}^{\mathbb{P}^{n_k}_{x_0}}\Big[\Big|\Big(b_u(X^{i,n_k}_u,\lambda^\ell_u) - b^\delta_u(X^{i,n_k}_u,\lambda^\ell_u)\Big)\cdot\nabla\phi(X^{i,n_k}_u)\Big|\Big]\mathrm{d}u \\ &= \int_s^t \mathbf{E}^{\mathbb{P}^{n_k}_{x_0}}\Big[\Big|\Big(b_u(X^{1,n_k}_u,\lambda^\ell_u) - b^\delta_u(X^{1,n_k}_u,\lambda^\ell_u)\Big)\cdot\nabla\phi(X^{1,n_k}_u)\Big|\Big]\mathrm{d}u\,. \end{aligned}$$

In the last line, we can now use the fact that by Lemma 3.3, the law of $X^{1,n_k}_u$ is absolutely continuous relative to Lebesgue measure. Writing $p^{1,n_k}_t \coloneqq p^{1,n_k}(t)$ for the function in (3.8), Hölder's inequality with conjugate exponents $1/r + 1/r' = 1$ and the integrability bound (3.9) show that

$$\begin{aligned} &\int_s^t \mathbf{E}^{\mathbb{P}^{n_k}_{x_0}}\Big[\Big|\Big(b_u(X^{1,n_k}_u,\lambda^\ell_u) - b^\delta_u(X^{1,n_k}_u,\lambda^\ell_u)\Big)\cdot\nabla\phi(X^{1,n_k}_u)\Big|\Big]\mathrm{d}u \\ &= \int_s^t\int_{\mathbb{R}^d}\Big|\Big(b_u(x,\lambda^\ell_u) - b^\delta_u(x,\lambda^\ell_u)\Big)\cdot\nabla\phi(x)\Big|\,p^{1,n_k}_t(x)\mathrm{d}x\,\mathrm{d}u \\ &\le c_{r,\gamma}\|\phi\|^*_1\int_s^t(1\wedge u)^{-\gamma}\|b_u(\lambda^\ell_u) - b^\delta_u(\lambda^\ell_u)\|_{\mathbb{L}^{r'}(K_\phi,\mathrm{d}x)}\,\mathrm{d}u\,, \end{aligned} \tag{4.36}$$

where $c_{r,\gamma} \coloneqq c_{\text{Lem. 3.3}} < \infty$. Since $\gamma < 1$, we can choose $w > 1$ to have $w\gamma < 1$. Then $c_{\gamma,w} \coloneqq (\int_0^T(1\wedge u)^{-\gamma w}\,\mathrm{d}u)^{1/w} < \infty$, and Hölder's inequality with conjugate exponents $1/w + 1/w' = 1$ applied to the right-hand side of (4.36) shows for (4.35) that

$$\begin{aligned} &\int_s^t \mathbf{E}^{\tilde{\mathbf{P}}}\left[\tilde{\Lambda}^k_u\Big[\Big|\Big(b_u(\lambda^\ell_u) - b^\delta_u(\lambda^\ell_u)\Big)\cdot\nabla\phi\Big|\Big]\right]\mathrm{d}u \\ &\le c_{r,\gamma,w,\phi}\Big(\int_s^t\|b_u(\lambda^\ell_u) - b^\delta_u(\lambda_\ell)\|^{w'}_{\mathbb{L}^{r'}(K_\phi,\mathrm{d}x)}\,\mathrm{d}u\Big)^{\frac{1}{w'}} \end{aligned} \tag{4.37}$$

for each $\ell \in [m]$, with $c_{r,\gamma,w} \coloneqq c_{r,\gamma}c_{\gamma,w}\|\phi\|^*_1$. Returning to (4.34) and using (4.37) gives

$$\begin{aligned} C_{\delta,k,2} \le\ & c_{r,\gamma,w,\phi}\|g\|_\infty\sum_{\ell=1}^m\Big(\int_s^t\|b_u(\lambda^\ell_u) - b^\delta_u(\lambda_\ell)\|^{w'}_{\mathbb{L}^{r'}(K_\phi,\mathrm{d}x)}\,\mathrm{d}u\Big)^{\frac{1}{w'}} \\ &+ 2\|g\|_\infty\|\phi\|^*_1\varepsilon_0 + 2\|g\|_\infty\|b\|_\infty\|\phi\|^*_1(t-s)\varepsilon_1\,. \end{aligned} \tag{4.38}$$

***Step 2.6*** We now choose our initially arbitrary $\varepsilon_0 > 0$ and $\varepsilon_1 > 0$ from Steps 2.1 and 2.2, respectively, such that

$$2\|g\|_\infty \|\phi\|_1^* \varepsilon_0 + 2\|g\|_\infty \|b\|_\infty \|\phi\|_1^* T\varepsilon_1 < \varepsilon/4\,.$$

Note that this choice is made independently of $s$ and $t$. Having chosen $\varepsilon_0$ and $\varepsilon_1$, we have also fixed our compact set $K$ via (4.29), the open sets $U_\lambda$ in (4.31) covering it, and thus the finite number $m$ appearing in the sum in (4.38) in Step 2.3. By Lemma 4.3, part 2), we can now choose $\delta_C(\varepsilon)$ to have

$$c_{\phi,r,\gamma,\phi}\|g\|_\infty \sum_{\ell=1}^{m} \int_s^t \|b_u(\lambda_u^\ell) - b_u^\delta(\lambda_u^\ell)\|_{\mathbb{L}^{r'}(K_\phi,\mathrm{d}x)} < \varepsilon/4$$

whenever $0 < \delta < \delta_C(\varepsilon)$. Provided $s > 0$, the choices of $\varepsilon_0$, $\varepsilon_1$ and $\delta$ show via (4.38) in Step 2.5 that for $k \in \mathbb{N}$ and $\delta < \delta_C(\varepsilon)$, we have

$$|C_{\delta,k,2}| = |\mathbf{E}^{\tilde{\mathbf{P}}}[\tilde{C}_{\delta,k,2}]| < \varepsilon/2\,. \tag{4.39}$$

***Step 2.7*** Now consider the case $s = 0$. While the general strategy is familiar from the proof of Lemma 4.4, implementing it requires slightly more care. Indeed, for $s' \in (0,t)$ to be chosen momentarily, we split in (4.34) the integral $\int_0^t = \int_0^{s'} + \int_{s'}^t$ to get

$$\begin{aligned} C_{\delta,k,2} \le{}& \|g\|_\infty \sum_{\ell=1}^{m} \int_0^{s'} \mathbf{E}^{\tilde{\mathbf{P}}}\left[\tilde{\Lambda}_u^k\Big[\Big|\Big(b_u^\delta(\lambda_u^\ell) - b_u(\lambda_u^\ell)\Big)\cdot\nabla\phi\Big|\Big]\right]\mathrm{d}u \\ &+ \|g\|_\infty \sum_{\ell=1}^{m} \int_{s'}^{t} \mathbf{E}^{\tilde{\mathbf{P}}}\left[\tilde{\Lambda}_u^k\Big[\Big|\Big(b_u^\delta(\lambda_u^\ell) - b_u(\lambda_u^\ell)\Big)\cdot\nabla\phi\Big|\Big]\right]\mathrm{d}u \\ &+ 2\|g\|_\infty\|\phi\|_1^*\varepsilon_0 + 2\|g\|_\infty\|b\|_\infty\|\phi\|_1^*(t-s)\varepsilon_1\,. \end{aligned}$$

We first choose $\varepsilon_0$, $\varepsilon_1$ as in Step 2.6, replacing $\varepsilon/4$ by $\varepsilon/8$, to have

$$2\|g\|_\infty \|\phi\|_1^* \varepsilon_0 + 2\|g\|_\infty \|b\|_\infty \|\phi\|_1^* T\varepsilon_1 < \varepsilon/8\,.$$

Note that this choice is independent of $s'$, and recall from the comment we made in Step 2.6 that fixing $\varepsilon_0$ and $\varepsilon_1$ also fixes $m$ in (4.34). Therefore, also $m$ is independent of $s'$. Since $\tilde{\Lambda}^k$ is $\mathbf{M}_1^+$-valued for all $k \in \mathbb{N}$, $b$ is bounded by Assumption 1.5 and $b^\delta$ is bounded uniformly over $\delta > 0$ by part 3) of Lemma 4.3, we can choose $s' > 0$ to have

$$\|g\|_\infty \sum_{\ell=1}^{m} \int_0^{s'} \mathbf{E}^{\tilde{\mathbf{P}}}\left[\tilde{\Lambda}_u^k\Big[\Big|\Big(b_u^\delta(\lambda_u^\ell) - b_u(\lambda_u^\ell)\Big)\cdot\nabla\phi\Big|\Big]\right]\mathrm{d}u < \varepsilon/4$$

uniformly for $\delta > 0$ and $k \in \mathbb{N}$, and without changing $m$. For the remaining term involving the integral from $s'$ to $t$, we can now follow Steps 2.5 and 2.6 and choose $\delta_C(\varepsilon) > 0$ to have for any $k \in \mathbb{N}$ and $0 < \delta < \delta_C(\varepsilon)$ that

$$\begin{aligned} &\|g\|_\infty \sum_{\ell=1}^{m} \int_{s'}^{t} \mathbf{E}^{\tilde{\mathbf{P}}}\left[\tilde{\Lambda}_u^k\Big[\Big|\Big(b_u^\delta(\lambda_u^\ell) - b_u(\lambda_u^\ell)\Big)\cdot\nabla\phi\Big|\Big]\right]\mathrm{d}u \\ &\le 2\|g\|_\infty\|\phi\|_1^*\varepsilon_0 + 2\|g\|_\infty\|b\|_\infty\|\phi\|_1^* T\varepsilon_1 < \varepsilon/8\,. \end{aligned}$$

Combining the above shows just as in (4.39) that also for $s = 0$, we have $|C_{\delta,k,2}| < \varepsilon/2$.

***Step 3*** We now combine the estimates from Steps 1 and 2. Our choices of $k_C(\varepsilon)$ from (4.26) and $\delta_C(\varepsilon)$ from (4.39) imply via (4.24) that for any $k > k_C(\varepsilon)$ and any $0 < \delta < \delta_C(\varepsilon)$, we have

$$|C_{k,\delta}| = |\mathbf{E}^{\tilde{\mathbf{P}}}[\tilde{C}_{k,\delta}]| \le |\mathbf{E}^{\tilde{\mathbf{P}}}[\tilde{C}_{k,\delta,1}]| + |\mathbf{E}^{\tilde{\mathbf{P}}}[\tilde{C}_{\delta,k,2}]| < \varepsilon\,.$$

This is the bound we sought, and so the proof is complete. □

We are left with the term $D_n$, which requires Assumption 1.5.

**Lemma 4.7** | *Under Assumptions 1.2 and 1.5, fix $\phi \in \mathbf{C}_c^\infty(\mathbb{R}^d)$ and a bounded measurable $g : [0,T] \times \Omega \to \mathbb{R}$. Then for any $\varepsilon > 0$, there exists $k_D(\varepsilon) > 0$ such that $|D_k| < \varepsilon/4$ for all $k > k_D(\varepsilon)$. In other words, we have $\lim_{k\to\infty} |D_k| = 0$.*

***Proof*** As earlier, $K_\phi := \mathrm{cl}_{(\mathbb{R}^d,|\,\cdot\,|)}\{x \in \mathbb{R}^d : |\nabla\phi(x)| > 0\}$ denotes the support of the gradient of $\phi$. In view of (4.17) we then have $\tilde{\mathbf{P}}$-a.s. simultaneously for $u \in [0,T]$ by Assumption 1.5 that

$$\lim_{k\to\infty} \sup_{x\in K_\phi} |b_u(x,\tilde{\Lambda}^\infty_u) - b_u(x,\tilde{\Lambda}^k_u)| = 0\,,$$

and by the Hölder-continuity from Assumption 1.2(H) that

$$\lim_{k\to\infty} \sup_{x\in K_\phi} |\sigma_u(x,\tilde{\Lambda}^\infty_u) - \sigma_u(x,\tilde{\Lambda}^k_u)| = 0\,,$$
$$\lim_{k\to\infty} \sup_{x\in K_\phi} |\bar{\sigma}_u(x,\tilde{\Lambda}^\infty_u) - \bar{\sigma}_u(x,\tilde{\Lambda}^k_u)| = 0\,.$$

With this, the definition of $\mathscr{L}_u(\tilde{\Lambda}^\infty_u)$ in (2.8) shows that $\tilde{\mathbf{P}}$-a.s., we have simultaneously for all $u \in [0,T]$ that

$$\lim_{k\to\infty} \sup_{x\in K_\phi} \Big|\Big(\Big(\mathscr{L}_u(\tilde{\Lambda}^\infty_u) - \mathscr{L}_u(\tilde{\Lambda}^k_u)\Big)\phi\Big)(x)\Big| = 0\,. \tag{4.40}$$

To make use of (4.40), note that the definition of $\tilde{D}_k$ in (4.12) gives

$$\begin{aligned}|\tilde{D}_k| &= \Big|\int_s^t g_u(\tilde{\Lambda}^k)\Big(\tilde{\Lambda}^k_u[\mathscr{L}_u(\tilde{\Lambda}^\infty_u)\phi] - \Lambda^k_u[\mathscr{L}_u(\tilde{\Lambda}^k_u)\phi]\Big)\,\mathrm{d}u\Big| \\ &\le \|g\|_\infty \int_s^t \Big|\tilde{\Lambda}^k_u\Big[\Big(\mathscr{L}_u(\tilde{\Lambda}^\infty_u) - \mathscr{L}_u(\tilde{\Lambda}^k_u)\Big)\phi\Big]\Big|\,\mathrm{d}u\,.\end{aligned}$$

Using the estimate (2.32) from Remark 2.7, we see that

$$\Big|\tilde{\Lambda}^k_u\Big[\Big(\mathscr{L}_u(\tilde{\Lambda}^\infty_u) - \mathscr{L}_u(\tilde{\Lambda}^k_u)\Big)\phi\Big]\Big| \le \sup_{x\in K_\phi}\Big|\Big(\Big(\mathscr{L}_u(\tilde{\Lambda}^\infty_u) - \mathscr{L}_u(\tilde{\Lambda}^k_u)\Big)\phi\Big)(x)\Big| \le 2c_{b,\sigma,\bar{\sigma}}\|\phi\|_2^*\,, \tag{4.41}$$

simultaneously for all $u \in [s,t]$. So dominated convergence along with (4.40) lets us conclude that $D_k = \mathbf{E}^{\tilde{\mathbf{P}}}[\tilde{D}_k] \to 0$ as $k \to \infty$. Hence there is $k_D(\varepsilon)$ such that $k > k_D(\varepsilon)$ implies $|D_k| \le \varepsilon/4$. □

## 4.2 Proof of Proposition 4.2 and Theorem 1.6

We are ready to prove Proposition 4.2 which, to recall, claims that

$$\lim_{k\to\infty} \mathbf{E}^{\mathbf{P}^{n_k}_{x_0}}\Big[\int_s^t g_u(\Lambda)\mathscr{A}_u(\Lambda_u)[\phi]\,\mathrm{d}u\Big] = \mathbf{E}^{\mathbf{P}^{\infty}_{x_0}}\Big[\int_s^t g_u(\Lambda)\mathscr{A}_u(\Lambda_u)[\phi]\,\mathrm{d}u\Big] \tag{4.1}$$

for any $\phi \in \mathcal{S}$, all $s,t \in [0,T]$ with $s \le t$, and all bounded measurable functions $g : [0,T]\times\Omega \to \mathbb{R}$ such that for each $u \in [0,T]$, the real-valued function $g_u(\,\cdot\,) \coloneqq g(u,\,\cdot\,)$ on $\Omega$ is continuous.

***Proof of Proposition 4.2*** Fix $\phi \in \mathcal{S}$ and $\varepsilon > 0$, and recall from the definition (2.26) that $\mathscr{A}_u(\Lambda_u)[\phi] = \Lambda_u[\mathscr{L}_u(\Lambda_u)\phi]$. As seen after Proposition 4.2, to establish (4.1), we can instead prove (4.4), i.e., that there exists $k_0(\varepsilon) \in \mathbb{N}$ such that for all $k > k_0(\varepsilon)$, we have

$$\bigg|\mathbf{E}^{\tilde{\mathbf{P}}}\Big[\int_s^t g_u(\tilde\Lambda^k)\tilde\Lambda^k_u[\mathscr{L}_u(\tilde\Lambda^k_u)\phi]\,\mathrm{d}u - \int_s^t g_u(\tilde\Lambda^\infty)\tilde\Lambda^\infty_u[\mathscr{L}_u(\tilde\Lambda^\infty_u)\phi]\,\mathrm{d}u\Big]\bigg| < \varepsilon\,. \tag{4.42}$$

We first claim that there exists a function $\psi \in \mathbf{C}^\infty_{\mathrm{c}}$ such that for all $k \in \mathbb{N}\cup\{\infty\}$, we have

$$\bigg|\mathbf{E}^{\tilde{\mathbf{P}}}\Big[\int_s^t g_u(\tilde\Lambda^k)\tilde\Lambda^k_u[\mathscr{L}_u(\tilde\Lambda^k_u)(\phi-\psi)]\,\mathrm{d}u\Big]\bigg| < \varepsilon/3\,. \tag{4.43}$$

Indeed, Assumptions 1.2 and 1.5 allow us to apply the estimates from (2.32) in Remark 2.7 to find that

$$\bigg|\mathbf{E}^{\tilde{\mathbf{P}}}\Big[\int_s^t g_u(\tilde\Lambda^k)\tilde\Lambda^k_u[\mathscr{L}_u(\tilde\Lambda^k_u)(\phi-\psi)]\,\mathrm{d}u\Big]\bigg| \le c_{b,\sigma,\bar\sigma}\|g\|_\infty(t-s)\|\phi-\psi\|_2^*\,. \tag{4.44}$$

Since $\mathbf{C}^\infty_{\mathrm{c}}$ is dense in $\mathcal{S}$, we can choose $\psi \in \mathbf{C}^\infty_{\mathrm{c}}$ to make $\|\phi-\psi\|_2^*$ and hence the right-hand side above arbitrarily small.

Now fix $\psi \in \mathbf{C}^\infty_{\mathrm{c}}$ satisfying (4.43) for all $k \in \mathbb{N}\cup\{\infty\}$. Then linearity and twice the triangle inequality show that to prove (4.42), it is sufficient to establish that

$$\begin{aligned}&\bigg|\mathbf{E}^{\tilde{\mathbf{P}}}\Big[\int_s^t g_u(\tilde\Lambda^k)\tilde\Lambda^k_u[\mathscr{L}_u(\tilde\Lambda^k_u)\psi]\,\mathrm{d}u - \int_s^t g_u(\tilde\Lambda^\infty)\tilde\Lambda^\infty_u[\mathscr{L}_u(\tilde\Lambda^\infty_u)\psi]\,\mathrm{d}u\Big]\bigg|\\ &\le A_\delta + B_{\delta,k} + C_{\delta,k} + D_k < \varepsilon/3\,,\end{aligned}$$

where $A_\delta$, $B_{\delta,k}$, $C_{\delta,k}$, $D_k$ are defined as in (4.13) and (4.9)–(4.12), with $\phi$ there replaced by $\psi$ here. But this is immediate from Lemmas 4.4, 4.5, 4.6 and 4.7 by which we can choose for any $\delta < \min\{\delta_A(\varepsilon),\delta_B(\varepsilon),\delta_C(\varepsilon)\}$ a $k > \max\{k_B(\varepsilon,\delta),k_C(\varepsilon),k_D(\varepsilon)\}$ with $\max\{A_\delta,B_{\delta,k},C_{\delta,k},D_k\} < \varepsilon/12$. This shows (4.42) and concludes the proof. □

We can now proceed to the proof of our main result in Theorem 1.6. To make its proof more transparent, we start with an easy auxiliary lemma.

**Lemma 4.8** | *Impose Assumption 1.2. For all $\phi,\psi \in \mathcal{S}$ and $t \in [0,T]$, the mapping*

$$\lambda \mapsto \mathscr{Q}_t(\lambda)[\phi,\psi] = \lambda[\mathscr{R}_t(\lambda)\phi]\cdot\lambda[\mathscr{R}_t(\lambda)\psi] = \lambda[\bar\sigma_t^\top(\lambda)\nabla\phi]\cdot\lambda[\bar\sigma_t^\top(\lambda)\nabla\psi]$$

*defined in* (2.27) *is continuous on $\mathcal{P}_{\mathrm{wk}^*}$, the space of probability measures on $\mathcal{B}(\mathbb{R}^d)$ equipped with the narrow topology.*

***Proof*** Let $\lambda^n \to \lambda$ in $\mathcal{P}_{\mathrm{wk}^*}$ as $n \to \infty$ and write

$$|\lambda^n[\bar\sigma_t^{\mathsf{T}}(\lambda^n)\nabla\phi] - \lambda[\bar\sigma_t^{\mathsf{T}}(\lambda)\nabla\phi]| \le \left|\lambda^n\Big[\Big(\bar\sigma_t^{\mathsf{T}}(\lambda^n) - \bar\sigma_t^{\mathsf{T}}(\lambda)\Big)\nabla\phi\Big]\right| + |(\lambda^n - \lambda)[\bar\sigma_t^{\mathsf{T}}(\lambda)\nabla\phi]|\,.$$

We consider the two terms on the right-hand side in turn.

For the first, since Hölder-continuity implies uniform continuity, Assumption 1.2 implies $\lim_{n\to\infty}\|\bar\sigma_t(\,\cdot\,,\lambda^n) - \bar\sigma_t(\,\cdot\,,\lambda)\|_\infty = 0$ for all $t \in [0,T]$. Since $\lambda^n$ is a probability measure, we moreover have that $|\lambda^n[\bar\sigma_t^{\mathsf{T}}(\lambda^n)\nabla\phi] - \lambda^n[\bar\sigma_t^{\mathsf{T}}(\lambda)\nabla\phi]| \le \|\phi\|_1^*\|\bar\sigma_t(\,\cdot\,,\lambda^n) - \bar\sigma_t(\,\cdot\,,\lambda)\|_\infty$. In combination, this shows that for all $t \in [0,T]$, we have

$$\lim_{n\to\infty}|\lambda^n[\bar\sigma_t^{\mathsf{T}}(\lambda^n)\nabla\phi] - \lambda^n[\bar\sigma_t^{\mathsf{T}}(\lambda)\nabla\phi]| = 0\,.$$

For the second term, since $\lambda^n \to \lambda$ narrowly and $x \mapsto \bar\sigma_t(x,\lambda)$ is continuous and bounded by Assumption 1.2, we get

$$\lim_{n\to\infty}|\lambda^n[\bar\sigma_t^{\mathsf{T}}(\lambda)\nabla\phi] - \lambda[\bar\sigma_t^{\mathsf{T}}(\lambda)\nabla\phi]| = 0$$

for each $t \in [0,T]$. The assertion now follows. □

To prove Theorem 1.6, we must show that a cluster point $\mathbf{P}_{x_0}^\infty$ solves the martingale problem $\mathrm{MP}(\delta_{x_0}, \mathscr{A}, \mathscr{Q})$, i.e., satisfies the requirements of Definition 2.3 with $\nu_0 = \delta_{x_0}$, $A = \mathscr{A}$ and $Q = \mathscr{Q}$. More specifically, we must show that under the measure $\mathbf{P}_{x_0}^\infty$, the decomposition

$$\Lambda_t = \delta_{x_0} + \int_0^t \mathscr{A}_s(\Lambda_s)\,\mathrm{d}s + M_t\,, \qquad t \in [0,T]\,, \tag{4.45}$$

from (1.13) is such that

1) under $\mathbf{P}_{x_0}^\infty$ we have that $\Lambda_0 \overset{d}{=} \nu$;

2) the absolutely continuous part is almost surely finite-valued, i.e. for all $\phi \in \mathcal{S}$ we have that $\mathbf{P}_{x_0}^\infty[\int_0^T |\mathscr{A}_t(\Lambda_t)[\phi]|\,\mathrm{d}t < \infty] = 1$;

3) the process $M$ defined by (4.45) is an $\mathcal{S}'$-valued local $(\mathbb{F}, \mathbf{P}_{x_0}^\infty)$-martingale;

4) $M$ has the quadratic variation process $\langle\!\langle M\rangle\!\rangle = \int \mathscr{Q}_t\,\mathrm{d}t$.

With Proposition 4.2 in our hands, the proof is easy and follows classical strategies; see e.g. Karatzas and Shreve [38, Ch. 5.4] or Vaillancourt [65] for similar proofs.

***Proof of Theorem 1.6*** By Proposition 2.6, the canonical process $\Lambda$ is for each $n \in \mathbb{N}$ an $\mathcal{S}'$-valued $(\mathbb{F}, \mathbf{P}_{x_0}^n)$-semimartingale with continuous trajectories, with canonical decomposition given by (4.45), and with martingale part $M \in \mathscr{M}_2^c(\mathbf{P}_{x_0}^n; \mathcal{S}')$. Hence for $\phi \in \mathcal{S}$, the process $\Lambda[\phi] = (\Lambda_t[\phi])_{t\in[0,T]}$ is a continuous real-valued $(\mathbb{F}, \mathbf{P}_{x_0}^n)$-semimartingale with canonical decomposition

$$\Lambda_t[\phi] = \delta_{x_0}[\phi] + \Big(\int_0^t \mathscr{A}_u(\Lambda_u)\,\mathrm{d}u\Big)[\phi] + M_t[\phi] = \delta_{x_0}[\phi] + \int_0^t \mathscr{A}_u(\Lambda_u)[\phi]\,\mathrm{d}u + M_t[\phi] \tag{4.46}$$

for $t \in [0,T]$, where $M[\phi] = (M_t[\phi])_{t\in[0,T]} \in \mathscr{M}_2^c(\mathbf{P}_{x_0}^n; \mathbb{R})$ has by (2.30) the quadratic variation process

$$\begin{aligned}\langle M[\phi]\rangle_t &= \Big(\int_0^t \mathscr{Q}_u(\Lambda_u)\,\mathrm{d}u + \frac{1}{n}\int_0^t \mathscr{C}_u(\Lambda_u)\,\mathrm{d}u\Big)[\phi,\phi] \\ &= \int_0^t \mathscr{Q}_u(\Lambda_u)[\phi,\phi]\,\mathrm{d}u + \frac{1}{n}\int_0^t \mathscr{C}_u(\Lambda_u)[\phi,\phi]\,\mathrm{d}u\,. \end{aligned}\tag{4.47}$$

This uses that $\mathbf{P}^n_{x_0}$ solves by Proposition 2.6 the martingale problem $\mathrm{MP}(\delta_{x_0}, \mathscr{A}, \mathscr{Q} + \frac{1}{n}\mathscr{C})$.

***Step 1*** Let $\mathbf{P}^\infty_{x_0}$ be a cluster point of $(\mathbf{P}^n_{x_0})_{n\in\mathbb{N}}$ and $(\mathbf{P}^{n_k}_{x_0})_{k\in\mathbb{N}}$ a subsequence with $\mathbf{P}^{n_k}_{x_0} \to \mathbf{P}^\infty_{x_0}$ narrowly in $\mathbf{M}^+_1(\Omega)$ as $k \to \infty$, where $\Omega = \mathbf{C}([0,T]; \mathcal{S}')$. We verify that $\mathbf{P}^\infty_{x_0}$ is a solution to the martingale problem $\mathrm{MP}(\delta_{x_0}, \mathscr{A}, \mathscr{Q})$. First note that since $X^{1,n}_0 \overset{d}{=} \delta_{x_0}$ and $(X^{i,n}_0)_{i=1,\dots,n}$ are i.d.d. for all $n \in \mathbb{N}$, the law of large numbers gives that $\Lambda_0[\phi] \overset{d}{=} \delta_{x_0}[\phi]$ under $\mathbf{P}^\infty_{x_0}$ for all $\phi \in \mathcal{S}$. We therefore have part 1) of Definition 2.3.

Next note that $\mathscr{A}_s(\Lambda_s)[\phi] = \Lambda_s[\mathscr{L}_s(\Lambda_s)\phi]$ by (2.26), from which the estimate (2.32) in Remark 2.7 gives the bound

$$\int_0^T |\mathscr{A}_t(\Lambda_t)[\phi]|\,\mathrm{d}t = \int_0^T |\Lambda_t[\mathscr{L}_t(\Lambda_t)\phi]|\,\mathrm{d}t \le T c_{b,\sigma,\bar\sigma} \|\phi\|^*_2 < \infty\,;$$

so part 2) of Definition 2.3 is satisfied. For parts 3) and 4), let $F \in \mathcal{S}$ and $0 \le s \le t \le T$. Then Itô's formula, (4.46) and (4.47) give

$$\begin{aligned} &F(\Lambda_t[\phi]) - F(\Lambda_s[\phi]) - \int_s^t F'(\Lambda_u[\phi])\mathscr{A}_u(\Lambda_u)[\phi]\,\mathrm{d}u \\ &\quad - \frac{1}{2}\int_s^t F''(\Lambda_u[\phi])\mathscr{Q}_u(\Lambda_u)[\phi,\phi]\,\mathrm{d}u - \frac{1}{2n}\int_s^t F''(\Lambda_u[\phi])\mathscr{C}_u(\Lambda_u)[\phi,\phi]\,\mathrm{d}u \\ &= \int_s^t F'(\Lambda_u[\phi])\,\mathrm{d}M_u[\phi]\,. \end{aligned} \tag{4.48}$$

Since $F \in \mathcal{S}$, $F'$ is bounded so that the process $\int F'(\Lambda_u[\phi])\mathrm{d}M_u[\phi]$ is an $(\mathbb{F}, \mathbf{P}^n_{x_0})$-martingale. Conditioning (4.48) on $\mathcal{F}_s$ thus gives for all continuous, bounded, $\mathcal{F}_s$-measurable functions $G : \Omega \to \mathbb{R}$ that

$$\mathbf{E}^{\mathbf{P}^{n_k}_{x_0}}\Big[\Big(\int_s^t F'(\Lambda_u[\phi])\mathrm{d}M_u[\phi]\Big)G(\Lambda)\Big] = 0 \qquad \text{for } k \in \mathbb{N}\,. \tag{4.49}$$

For the individual terms on the left-hand side of (4.48), we note the following. Clearly, by the narrow convergence $\mathbf{P}^{n_k}_{x_0} \to \mathbf{P}^\infty_{x_0}$ and the continuity of $G$ and of the evaluation map $\phi \mapsto \Lambda[\phi]$, we get

$$\lim_{k\to\infty} \mathbf{E}^{\mathbf{P}^{n_k}_{x_0}}\Big[\Big(F(\Lambda_t[\phi]) - F(\Lambda_s[\phi])\Big)G(\Lambda)\Big] = \mathbf{E}^{\mathbf{P}^\infty_{x_0}}\Big[\Big(F(\Lambda_t[\phi]) - F(\Lambda_s[\phi])\Big)G(\Lambda)\Big]\,. \tag{4.50}$$

Consider now the map $g : [0,T] \times \Omega \to \mathbb{R}$ defined by $g_u(\Lambda) := g(u,\Lambda) := G(\Lambda)F'(\Lambda_u[\phi])$. It is clearly bounded and measurable. In addition, $\lambda \mapsto g_u(\lambda)$ is continuous on $\Omega$ for each $u \in [0,T]$. So $g$ satisfies the hypotheses of Proposition 4.2, and we get that

$$\begin{aligned} &\lim_{k\to\infty} \mathbf{E}^{\mathbf{P}^{n_k}_{x_0}}\Big[\Big(\int_s^t F'(\Lambda_u[\phi])\mathscr{A}_u(\Lambda_u)[\phi]\,\mathrm{d}u\Big)G(\Lambda)\Big] \\ &= \lim_{k\to\infty} \mathbf{E}^{\mathbf{P}^{n_k}_{x_0}}\Big[\int_s^t g_u(\Lambda)\mathscr{A}_u(\Lambda_u)[\phi]\,\mathrm{d}u\Big] \\ &= \mathbf{E}^{\mathbf{P}^\infty_{x_0}}\Big[\int_s^t g_u(\Lambda)\mathscr{A}_u(\Lambda_u)[\phi]\,\mathrm{d}u\Big] \\ &= \mathbf{E}^{\mathbf{P}^\infty_{x_0}}\Big[\Big(\int_s^t F'(\Lambda_u[\phi])\,\mathscr{A}_u(\Lambda_u)[\phi]\,\mathrm{d}u\Big)G(\Lambda)\Big]\,. \end{aligned} \tag{4.51}$$

In addition, the function $\lambda \mapsto F''(\lambda_u[\phi])G(\lambda)$ is bounded, measurable and continuous on $\Omega$ for each $u \in [0,T]$. Hence applying Lemma 4.8 shows that

$$\begin{aligned}\lim_{k\to\infty} \mathbf{E}^{\mathbf{P}^{n_k}_{x_0}} &\left[\left(\frac{1}{2}\int_s^t F''(\Lambda_u[\phi])\mathscr{Q}_u(\Lambda_u)[\phi,\phi]\,\mathrm{d}u\right)G(\Lambda)\right] \\ &= \mathbf{E}^{\mathbf{P}^{\infty}_{x_0}} \left[\left(\frac{1}{2}\int_s^t F''(\Lambda_u[\phi])\mathscr{Q}_u(\Lambda_u)[\phi,\phi]\,\mathrm{d}u\right)G(\Lambda)\right]. \end{aligned} \tag{4.52}$$

Finally, the boundedness of $\sigma$ from Assumption 1.2, the smoothness of $\phi$ and the fact that $\Lambda$ is $\mathbf{M}_1^+$-valued $\mathbf{P}^{n_k}_{x_0}$-a.s. for all $k \in \mathbb{N}$ imply that $\mathscr{C}_u(\Lambda_u)[\phi,\phi] = \Lambda_u[(\sigma_t^{\mathsf{T}}(\Lambda_u)\nabla\phi)\cdot(\sigma_t^{\mathsf{T}}(\Lambda_u)\nabla\phi)]$ in (2.28) is uniformly bounded on $[0,T]\times\Omega$, so that

$$\lim_{k\to\infty} \frac{1}{n_k}\mathbf{E}^{\mathbf{P}^{n_k}_{x_0}}\left[\left(\int_s^t F''(\Lambda_u[\phi])\mathscr{C}_u(\Lambda_u)[\phi,\phi]\,\mathrm{d}u\right)G(\Lambda)\right] = 0\,. \tag{4.53}$$

Starting from (4.48), multiplying both sides by $G(\Lambda)$, taking expectations with respect to $\mathbf{P}^{n_k}_{x_0}$, passing to the limit $k\to\infty$ and using (4.49)–(4.53), we therefore see that for all $F \in \mathcal{S}$ and all continuous, bounded and $\mathcal{F}_s$-measurable $G : \Omega \to \mathbb{R}$, we have that

$$\begin{aligned}\mathbf{E}^{\mathbf{P}^{\infty}_{x_0}} \Bigg[\Bigg(F(\Lambda_t[\phi]) - F(\Lambda_s[\phi]) &- \int_s^t F'(\Lambda_u[\phi])\mathscr{A}_u(\Lambda_u)[\phi]\,\mathrm{d}u \\ &- \frac{1}{2}\int_s^t F''(\Lambda_u[\phi])\mathscr{Q}_u(\Lambda_u)[\phi,\phi]\,\mathrm{d}u\Bigg)G(\Lambda)\Bigg] = 0\,. \end{aligned} \tag{4.54}$$

By a standard functional monotone class argument, see for instance Ethier and Kurtz [20, Sec. 3, Eqn. (3.4) and Thm. A.4.3], (4.54) can be extended from continuous, bounded and $\mathcal{F}_s$-measurable $G$ to general bounded and $\mathcal{F}_s$-measurable $G$, and this shows that the process $N^{F,\phi} := (N^{F,\phi}_t)_{t\in[0,T]}$ defined by

$$\begin{aligned} N^{F,\phi}_t := F(\Lambda_t[\phi]) - F(\delta_{x_0}[\phi]) &- \int_0^t F'(\Lambda_u[\phi])\mathscr{A}_u(\Lambda_u)[\phi]\,\mathrm{d}u \\ &- \frac{1}{2}\int_0^t F''(\Lambda_u[\phi])\mathscr{Q}_u(\Lambda_u)[\phi,\phi]\,\mathrm{d}u\,, \qquad t\in[0,T], \end{aligned} \tag{4.55}$$

is a continuous real-valued $(\mathbb{F},\mathbf{P}^{\infty}_{x_0})$-martingale. Since $F \in \mathcal{S}$ is arbitrary, standard results on the characterization of diffusion processes, see e.g. Revuz and Yor [52, Prop. VII.2.4], show that for each $\phi \in \mathcal{S}$, the process $M^\phi = (M^\phi_t)_{t\in[0,T]}$ defined by

$$M^\phi_t := \Lambda_t[\phi] - \delta_{x_0}[\phi] - \int_0^t \mathscr{A}_u(\Lambda_u)[\phi]\,\mathrm{d}u\,, \qquad t\in[0,T], \tag{4.56}$$

is a continuous $(\mathbb{F},\mathbf{P}^{\infty}_{x_0})$-martingale with quadratic variation process

$$\langle M^\phi\rangle_t = \int_0^t \mathscr{Q}_u(\Lambda_u)[\phi,\phi]\,\mathrm{d}u\,, \qquad t\in[0,T]\,. \tag{4.57}$$

So Lemma 2.2 implies that $\Lambda$ is an $\mathcal{S}'$-valued $(\mathbb{F},\mathbf{P}^{\infty}_{x_0})$-semimartingale with continuous (in $\mathcal{S}'$) trajectories and decomposition given by (1.13) or equivalently (4.45). Moreover, from Equation (2.3), Remark 2.4 and polarization, we obtain the quadratic variation process $\mathscr{Q}$ in (1.14) directly from (4.57). This establishes parts 3) and 4) of Definition 2.3 and shows that $\mathbf{P}^{\infty}_{x_0}$ solves $\mathrm{MP}(\delta_{x_0},\mathscr{A},\mathscr{Q})$.

***Step 2*** We now verify that the empirical propagation of chaos property (1.15) holds. Since the sequence $(\mathbf{P}^n_{x_0})_{n\in\mathbb{N}}$ is tight by Lemma 1.4, there is a subsequence $(\mathbf{P}^{n_k}_{x_0})_{k\in\mathbb{N}}$ converging to $\mathbf{P}^\infty_{x_0}$ narrowly in $\mathbf{M}^+_1(\Omega)$. But by the definition (1.5), we have $\mathrm{Law}_{\mathbb{P}^n_{x_0}}(\mu^n) = \mathbf{P}^n_{x_0}$ for all $n \in \mathbb{N}$, and so (1.15) holds trivially. This concludes the proof. □

## 5 AN EXAMPLE

Let us now give a simple example of a drift coefficient that falls within our framework.

Consider for some $R > 0$ the kernel

$$k_R(y) \coloneqq \mathbb{1}_{B_R}(y)\,, \tag{5.1}$$

where $B_R \coloneqq B_R(0) \coloneqq \{x \in \mathbb{R}^d \,:\, |x| < 1\}$ and set

$$b(x,\mu) \coloneqq (k_R * \mu)(x) = \int_{\mathbb{R}^d} k_R(y-x)\,\mu(\mathrm{d}y) = \mu[B_R(x)]\,. \tag{5.2}$$

We sketch an argument to prove that $b : \mathbb{R}^d \times \mathbf{M}^+_1 \to \mathbb{R}$ satisfies Assumption 1.5. Let $\mu \in \mathbf{M}^+_1 \cap \mathbb{L}^r(\mathrm{d}x)$ with $r > 1$ and $K \subseteq \mathbb{R}^d$ be compact. Let $\mu_k \to \mu$ narrowly in $\mathbf{M}^+_1$ as $k \to \infty$. We must verify that

$$\lim_{k\to\infty} \sup_{x\in K} |b(x,\mu_k) - b(x,\mu)| = 0\,. \tag{5.3}$$

We show in i) below that $x \mapsto b(x,\mu)$ is continuous. The difficulty is to prove that in (5.3), the potentially discontinuous function $x \mapsto b(x,\mu_k)$ converges to the continuous function $x \mapsto b(x,\mu)$ uniformly on $K$ as $k \to \infty$. For this, we use a sandwich argument.

i) We argue that $K \ni x \mapsto b(x,\mu) = \mu[B_R(x)] \in \mathbb{R}$ is continuous and thus also uniformly continuous. Indeed, since $\mu \in \mathbb{L}^r(\mathrm{d}x)$, we get with Hölder's inequality that

$$|\mu[B_R(x)] - \mu[B_R(x')]| \le \|\mu\|_{\mathbb{L}^r(\mathrm{d}x)} \|\mathbb{1}_{B_R}(x) - \mathbb{1}_{B_R}(x')\|_{\mathbb{L}^{r'}(\mathrm{d}x)}\,.$$

Since $r' < \infty$, the right-hand side vanishes as $x \to x'$. This shows continuity. Uniform continuity follows from the compactness of $K$.

ii) Consider functions $\phi_j, \psi_j \in \mathbf{C}^\infty_c(\mathbb{R}^d)$ such that $\phi_j \uparrow \mathbb{1}_{B_R(0)}$ and $\psi_j \downarrow \mathbb{1}_{B_R(0)}$ as $j \to \infty$. Denote by $\tau_x f$ the translation of $f$ by $x$, i.e. $\tau_x f(y) = f(y-x)$. Note that for each fixed $k \in \mathbb{N}$, the family $\{\mathbb{R}^d \ni y \mapsto \tau_x\phi_j(y) \in \mathbb{R}\}_{x\in K} \subseteq \mathbf{C}^\infty_c(\mathbb{R}^d)$ is at each point $y \in \mathbb{R}^d$ equicontinuous. Since $\mu_k \to \mu$ narrowly, we get from Stroock and Varadhan [63, Cor. 1.1.5] that

$$\lim_{n\to\infty} \sup_{x\in K} |\mu_k[\tau_x\phi_j] - \mu[\tau_x\phi_j]| = 0\,,$$
$$\lim_{n\to\infty} \sup_{x\in K} |\mu_k[\tau_x\psi_j] - \mu[\tau_x\psi_j]| = 0\,.$$

iii) Similarly to step i), Hölder's inequality yields the bound

$$|\mu[B_R(x)] - \mu[\tau_x\phi_j]| \le \|\mu\|_{\mathbb{L}^r(\mathrm{d}x)} \|\tau_x(\mathbb{1}_{B_R(0)} - \phi_j)\|_{\mathbb{L}^{r'}(\mathrm{d}x)} \le \|\mu\|_{\mathbb{L}^r(\mathrm{d}x)} \|\mathbb{1}_{B_R(0)} - \phi_j\|_{\mathbb{L}^{r'}(\mathrm{d}x)}\,.$$

The same estimate holds with $\psi_j$ instead of $\phi_j$; so noting that $r' < \infty$, we get from the choice of the sequences $(\phi_j)_{j\in\mathbb{N}}$, $(\psi_j)_{j\in\mathbb{N}}$ that

$$\lim_{j\to\infty}\sup_{x\in K}|\mu[B_R(x)] - \mu[\tau_x\phi_j]| = 0\,,$$
$$\lim_{j\to\infty}\sup_{x\in K}|\mu[B_R(x)] - \mu[\tau_x\psi_j]| = 0\,.$$

iv) We now make ourselves a sandwich. Let $\varepsilon > 0$. Since $\mu$ is a probability measure and since $\phi_j \le \psi_j$ for all $j \in \mathbb{N}$, we can find by step iii) some $j_0 \in \mathbb{N}$ such that for all $x \in K$, we have

$$\mu[B_R(x)] - \varepsilon/2 < \mu[\tau_x\phi_{j_0}] \le \mu[\tau_x\psi_{j_0}] < \mu[B_R(x)] + \varepsilon/2\,.$$

By step ii), we can find $k_0 \in \mathbb{N}$ such that for all $k \ge k_0$ and $x \in K$,

$$\mu[B_R(x)] - \varepsilon < \mu_k[\tau_x\phi_{j_0}] \le \mu_k[\tau_x\psi_{j_0}] < \mu[B_R(x)] + \varepsilon\,.$$

But since $\tau_x\phi_j \le \mathbb{1}_{B_R(x)} \le \tau_x\psi_j$ for all $x \in \mathbb{R}^d$ and $\mu_n \in \mathbf{M}_1^+$, we get for all $k \in \mathbb{N}$ and $x \in K$ that $\mu_k[\tau_x\phi_{k_0}] \le \mu_k[B_R(x)] \le \mu_k[\tau_x\psi_{k_0}]$. We therefore have

$$\mu[B_R(x)] - \varepsilon < \mu_k[B_R(x)] < \mu[B_R(x)] + \varepsilon$$

for all $x \in K$ and $k \ge k_0$. This shows that (5.3) holds.

Since $r > 1$ was arbitrary, Assumption 1.5 holds for a drift $b$ of the form in (5.2). If we in addition consider $\sigma, \bar\sigma : [0,T]\times\mathbb{R}^d\times\mathbf{M}_1^+ \to \mathbb{R}^{d\times d}$ that satisfy Assumption 1.2, then we can use Theorem 1.6 to show that cluster points $\mathbf{P}^\infty_{x_0}$ of the sequence of empirical measure flow laws $(\mathbf{P}^n_{x_0})_{n\in\mathbb{N}}$ associated with the particle system (1.1), (1.2) with coefficients $b$, $\sigma$, $\bar\sigma$ solve the formally limiting martingale problem $\mathrm{MP}(\delta_{x_0}, \mathscr{A}, \mathscr{Q})$.

A drift of a comparable type is studied in the recent paper by Chen et al. [10]. There the authors consider a real-valued setting, i.e. $d = 1$, and study the Hegselmann–Krause model, which is given by the kernel

$$k_{\mathrm{HK}}(x) \coloneqq \mathbb{1}_{[0,R]}(|x|)x = \mathbb{1}_{B_R}(x)x$$

that induces the drift

$$b_{\mathrm{HK}}(x,\mu) \coloneqq -k_{\mathrm{HK}} * \mu(x) = -\int_{\mathbb{R}^d} k_{\mathrm{HK}}(x-y)\mu(\mathrm{d}x)\,.$$

The authors then fix $\kappa > 0$ and consider as diffusion coefficients a map $\sigma : [0,T]\times\mathbb{R}\to\mathbb{R}$ of class $\mathbf{C}^3_b(\mathbb{R})$ with $\sigma_t(x) \ge \kappa$ for all $(t,x) \in [0,T]\times\mathbb{R}$ and a constant $\bar\sigma \ge \kappa$. This gives rise to the $n$-particle systems

$$X^{i,n}_t = x_0 - \int_0^t (k_{\mathrm{HK}} * \mu^n_s)(X^{i,n}_s)\,\mathrm{d}s + \int_0^t \sigma_s(X^{i,n}_s)\,\mathrm{d}B^{i,n}_s + \int_0^t \bar\sigma_s\,\mathrm{d}Z_s\,, \tag{5.4}$$

$$\mu^n_t = \frac{1}{n}\sum_{i=1}^n \delta_{X^{i,n}_t}\,. \tag{5.5}$$

The authors next define a *regularized particle system* given for $\delta > 0$ in terms of a mollified kernel $k^\delta_{\mathrm{HK}} := k_{\mathrm{HK}} * h^\delta$ by

$$X^{i,n,\delta}_t = x_0 - \int_0^t (k^\delta_{\mathrm{HK}} * \mu^n_s)(X^{i,n,\delta}_s)\,\mathrm{d}s + \int_0^t \sigma_s(X^{i,n,\delta}_s)\,\mathrm{d}B^{i,n,\delta}_s + \bar\sigma Z^{n,\delta}_t, \tag{5.6}$$

$$\mu^{n,\delta}_t = \frac{1}{n}\sum_{i=1}^n \delta_{X^{i,n,\delta}_t}\,. \tag{5.7}$$

Just as in Lemma 1.4 we can define the measure $\mathbf{P}^{n,\delta}_{x_0} := \mathrm{Law}(\mu^{n,\delta}_t)$. In [10], the authors show that for a certain sequence $(\delta(n))_{n\in\mathbb{N}}$ with $\delta(n) \to 0$ as $n \to \infty$, the sequence $(\mathbf{P}^{\delta(n),n}_{x_0})_{n\in N}$ possess a unique limit, and that this limit solves the martingale problem $\mathrm{MP}(\delta_{x_0}, \mathscr{A}, \mathscr{Q})$ associated with $\delta_{x_0}$, $\mathscr{A}$, $\mathscr{Q}$ defined as in Sections 2.2 and 2.3 in terms of the coefficients $b_{\mathrm{HK}}$, $\sigma$, $\bar\sigma$ appearing in (5.4).

In our work, we do not consider the question of uniqueness. The Hegselmann–Krause drift bears strong similarity with the drift $b$ in (5.2) for which we sketched an argument to verify that it satisfies Assumption 1.5. We therefore expect that also $b_{\mathrm{HK}}$ satisfies this assumption. With Theorem 1.6, we could then deduce that the cluster points of the empirical measure flow laws $(\mathbf{P}^n_{x_0})_{n\in\mathbb{N}}$ of the *unregularized* $n$-particle systems (5.4) and (5.5) also solve the associated limiting martingale problem $\mathrm{MP}(\delta_{x_0}, \mathscr{A}, \mathscr{Q})$. In fact, with the uniqueness result of [10], the solution of $\mathrm{MP}(\delta_{x_0}, \mathscr{A}, \mathscr{Q})$ we find with Theorem 1.6 as the limit of the sequence of empirical measure flow laws of the unregularized particle systems (5.4) and (5.5) coincides with the solution found in [10], that is obtained via the modified, i.e., the regularized particle systems (5.6) and (5.7). We should note, however, that our existence result holds even in the case when $\sigma$, $\bar\sigma$ are more general functions than those appearing in (5.4) as long as they satisfy Assumption 1.2.

Thus, via Theorem 1.6 we should obtain the empirical propagation of chaos for the original particle systems. As we mentioned in the introduction, this can be used to construct an appropriate version of the classical notion of propagation of chaos for the associated McKean–Vlasov dynamics; see Crowell [14].

## 6 DISCUSSION AND CONCLUDING REMARKS

In our main result in Theorem 1.6, we showed that under Assumptions 1.2 and 1.5, any cluster point $\mathbf{P}^\infty_{x_0}$ of the tight sequence $(\mathbf{P}^n_{x_0})_{n\in\mathbb{N}}$ solves the martingale problem $\mathrm{MP}(\delta_{x_0}, \mathscr{A}, \mathscr{Q})$. Theorem 1.6 thus contains both an existence and a convergence result.

In aiming for such a result, we face certain barriers. Indeed, it can be shown that without continuity of the drift $b$ in the measure argument $\mu$, we cannot expect a solution of $\mathrm{MP}(\delta_{x_0}, \mathscr{A}, \mathscr{Q})$ to exist in general even for nice $\sigma$, $\bar\sigma$; see [12, Ex. 7.1]. Accordingly, Assumption 1.5 requires a certain type of *local* continuity in the measure argument $\mu$ but allows for irregularity in the space variable $x$. As we mention below, the specific type of continuity in the measure argument required to obtain only the existence of solutions or, in addition, convergence may vary depending on further structural properties of the drift $b$. In this paper, we aim to keep Assumption 1.5 simple and transparent while being general enough to accommodate relevant examples, such as those mentioned in Section 5. We want to emphasize, however, that our approach is flexible enough to be adapted to other settings.

With coefficients $b$, $\sigma$, $\bar{\sigma}$ that are jointly continuous in $(x, \mu)$, a result as in Theorem 1.6 is classical; see, for instance, Vaillancourt [65], Dawson and Vaillancourt [17], and Dawson [16]. These works build on studies that examine the case without common noise, i.e., when $\bar{\sigma} \equiv 0$; see [65] for references. In all these works, and more generally in the established solution theory for martingale problems, continuity of the coefficients plays a central role. Without it, standard weak-convergence techniques fail, creating obstacles that are not easily overcome. Thus, it is unsurprising that a results such as the one in Theorem 1.6, where some coefficients can be discontinuous, are rare; we mention a few notable examples below.

To clarify this point, let us draw a parallel to $\mathbb{R}$-valued martingale problems, which are commonly used to obtain weak solutions of one-dimensional SDEs; see, e.g., Karatzas and Shreve [38, Sec. 5.4]. Consider such an SDE with a bounded and measurable drift $\tilde{b}$ and a bounded, continuous, and uniformly elliptic (i.e., strictly positive and lower-bounded) diffusion coefficient $\tilde{\sigma}$. In this case, it is well-understood how to obtain a solution for the associated martingale problem $\mathrm{MP}(x_0, \tilde{b}, \tilde{\sigma}^2)$, formulated on the canonical space $\mathbf{C}([0,T];\mathbb{R})$ with canonical coordinate process $Y = (Y_t)_{t\in[0,T]}$. The orthodox strategy is to *first solve the martingale problem without the irregular drift*, i.e., $\mathrm{MP}(x_0, 0, \tilde{\sigma}^2)$, and then *transform the solution* using the Girsanov theorem. In this way, a solution to the original problem $\mathrm{MP}(x_0, \tilde{b}, \tilde{\sigma}^2)$ is obtained. Note that in this approach, the ellipticity of $\tilde{\sigma}$ guarantees a sufficiently rich set of transformations to induce the drift $\tilde{b}$. This allows us to avoid difficulties posed by the irregular drift in the first step and address them through a probabilistic argument only in the second step. In fact, the solution $\boldsymbol{X}^n$ in Lemma 1.3 of the $n$-particle system (1.1) and (1.2) is constructed precisely in this way.

However, adapting this strategy to solve the $\mathcal{S}'$-valued problem $\mathrm{MP}(\delta_{x_0}, \mathscr{A}, \mathscr{Q})$ fails in general. A notable exception is the case where $\sigma$ and $\bar{\sigma}$ in (1.1) and (1.2) do not depend on the measure argument; see for example [43]. Informally, this strategy fails in the general case because the tensor variation $\langle\!\langle M \rangle\!\rangle_t = \mathscr{Q}_t \,\mathrm{d}t$ is a finite-rank operator, and this limits the set of transformations for the infinite-dimensional process $\Lambda$ that can be induced by an absolutely continuous change of measure via the Girsanov theorem. Consequently, we cannot transform solutions of the martingale problem $\mathrm{MP}(\delta_{x_0}, 0, \mathscr{Q})$ without drift to include a general drift $\mathscr{A}$.

Thus, we must approach the full martingale problem $\mathrm{MP}(\delta_{x_0}, \mathscr{A}, \mathscr{Q})$ directly, and the irregular part $\mathscr{A}$ complicates this substantially. Indeed, even in the case of classical SDEs, solving $\mathrm{MP}(x_0, \tilde{b}, \tilde{\sigma}^2)$ with a discontinuous drift $\tilde{b}$ directly is not straightforward. We are unaware of a reference that carries out this analysis explicitly. Instead, we refer to Figalli [21] or Stroock and Varadhan [63, Sec. 11.3], where the necessary ideas appear more or less explicitly. The key is to exploit the regularity of the density of $\mathrm{Law}(Y_t)$ for $t \in (0,T]$ in an approximation argument for the irregular drift $\tilde{b}$. Note that in this approach, the ellipticity of $\tilde{\sigma}$ plays a different role compared to the argument involving the Girsanov transform. Ellipticity is the source of the regularity of $\mathrm{Law}(Y_t)$.

The same fundamental principle is used in our work, with the emergence of regularity being the key novel conceptual and technical ingredient that paves the way. It is paired with new arguments to handle the randomness stemming from the common noise. The combination of these ideas lets us solve $\mathrm{MP}(\delta_{x_0}, \mathscr{A}, \mathscr{Q})$.

Let us give a high-level view of the central steps. The key technical step is provided by Proposition 4.2. Its proof is efficiently approached via the Skorohod representation. The benefit of using the Skorohod construction is that it allows us to consider only those convergent sequences $(\tilde{\Lambda}^k)_{k\in\mathbb{N}\cup\{\infty\}}$ that are charged under $\tilde{\mathbf{P}}$ and for which we need to establish continuity in mean along; see (4.42). The representation in terms of the random variables $(\tilde{\Lambda}^k)_{k\in\mathbb{N}\cup\{\infty\}}$ lets us define the terms $A_\delta$, $B_{\delta,k}$, $C_{\delta,k}$, and $D_k$ in (4.8)–(4.13), and thus transparently formulate the approximation argument that we carry out via Lemmas 4.4–4.7. To control $A_\delta$, the emergence of regularity from Theorem 3.1 is indispensable. Bounding the term $C_{\delta,k}$ is difficult because the Skorohod construction, while useful, necessarily couples the random variables $(\tilde{\Lambda}^k)_{k\in\mathbb{N}\cup\{\infty\}}$. To control the dependence this coupling introduces, we require the intricate compactness argument given in the proof of Lemma 4.6. It lets us access, via the exchangeability of the $n_k$-particle system in (1.1) and (1.2), the regularity of $\mathrm{Law}_{\mathbb{P}^{n_k}_{x_0}}(X^{1,n_k}_t)$, which Lemma 3.3 shows to be uniform in $k \in \mathbb{N}$.

Key for our argument is the emergence of regularity. Only with it can we show that cluster points $\mathbf{P}^\infty_{x_0}$ of the sequence $(\mathbf{P}^n_{x_0})_{n\in\mathbb{N}}$ solve the formally limiting martingale problem $\mathrm{MP}(\delta_{x_0}, \mathscr{A}, \mathscr{Q})$. An approach of a similar spirit, yet followed in a very different manner and exercised in a very different way, is found in Varadhan [66, Thm. 3.2]. There, a result in the same philosophy as the emergence of regularity in Theorem 3.1 is obtained in a setting without common noise and for dynamics of a different type than those in (1.1), (1.2). This is then used to solve a formally limiting martingale problem. However, the techniques used in [66] are based on entropy estimates that cannot be obtained in our setting and require additional structural assumptions and more regularity for the drift than we are willing to impose.

For certain special cases, existence and convergence results for the formally limiting martingale problem $\mathrm{MP}(\delta_{x_0}, \mathscr{A}, \mathscr{Q})$ with an irregular drift or irregular coefficients more generally are available in the literature. To discuss them, we should recall that solutions of the common noise McKean–Vlasov SDE with coefficients $b$, $\sigma$, $\bar{\sigma}$, or of the associated stochastic Fokker–Planck equation yield solutions of $\mathrm{MP}(\delta_{x_0}, \mathscr{A}, \mathscr{Q})$; see, e.g., Lacker et al. [44], Coghi and Gess [11].

When $b$, $\sigma$, $\bar{\sigma}$ depend linearly on $\mu$ for instance, Hammersley et al. [28] and Hammersley [27] obtain the existence of solutions of the McKean–Vlasov SDE and thus also for $\mathrm{MP}(\delta_{x_0}, \mathscr{A}, \mathscr{Q})$. These authors rely on classical regularity estimates for the time-marginals of diffusion processes such as those in Krylov [42] or in Stroock and Varadhan [63], and the assumed linearity of the coefficients in the measure argument $\mu$ allows the authors to use these estimates. More recently, Chen et al. [10] obtained results on existence, uniqueness, and a type of convergence for the Hegesmann–Krause drift mentioned in Section 5. This drift is again linear in the measure argument, and their technique requires restrictive assumptions on $\sigma$ and $\bar{\sigma}$, in particular independence of $\mu$. We should note that their convergence result holds for a sequence of empirical measure flow laws of smoothed particle systems converging to a solution of $\mathrm{MP}(\delta_{x_0}, \mathscr{A}, \mathscr{Q})$ when the smoothing vanishes as the system size increases; see Section 5.

Another set of examples is provided by rank-based models, i.e., models in which the dependence on the measure argument is through the quantile function of $\mu$. Works that obtain convergence results in this case include Kolli and Shkolnikov [40], Pal and Pitman [47], Pal and Shkolnikov [48] without common noise, and Shkolnikov and Yeung [57] in the case with common noise.

**Ramifications and extensions** As we mentioned in Section 1.3, our main interest in Theorem 1.6 is its connection to the existence of weak solutions and the propagation of chaos property for McKean–Vlasov systems and SDEs in the cases where the drift possesses only low regularity properties. Next to our intrinsic interest in Theorem 1.6, this result also serves as a stepping stone for various important extensions. One such instance we highlight is a more general existence results for solutions of $\mathrm{MP}(\nu, \mathscr{A}, \mathscr{Q})$.

**Assumption 6.1** | The probability measure $\nu$ on $\mathcal{B}(\mathbb{R}^d)$ is absolutely continuous with respect to Lebesgue measure.

**Assumption 6.2** | The drift coefficient $b : [0,T] \times \mathbb{R}^d \times \mathbf{M}_1^+(\mathbb{R}^d) \to \mathbb{R}^d$ is bounded and measurable relative to $\mathcal{E}$. In addition, there exists $r_0 > 1$ such that $b$ satisfies:

($\mathrm{C_w}$) If $(\mu_k)_{k\in\mathbb{N}\cup\{\infty\}}$ in $\mathbf{M}_1^+(\mathbb{R}^d)$ consists of measure that are absolutely continuous with respect to Lebesgue measures and $\mathrm{d}\mu_k/\mathrm{d}x \to \mathrm{d}\mu_\infty/\mathrm{d}x$ in $\mathbb{L}^1(\mathrm{d}x)$ as $k \to \infty$, then

$$\lim_{k\to\infty} \|b_t(\mu_k) - b_t(\mu_\infty)\|_{\mathbb{L}^1(K,\mathrm{d}x)} = 0$$

for any compact subset $K \subseteq \mathbb{R}^d$ and $t \in [0,T]$.

Clearly, Assumption 1.5 implies Assumption 6.2. In the companion paper Crowell [13], we obtain the following existence result, which we preview here.

**Theorem 6.3** | *Under Assumptions 1.2, 6.1 and 6.2, there exists a measure* $\mathbf{P}$ *on* $(\Omega, \mathcal{F})$ *solving the martingale problem* $\mathrm{MP}(\nu, \mathscr{A}, \mathscr{Q})$.

This particular extension is useful to treat certain cases of singular coefficients, more specifically, the case when the drift depends $b$ locally on the measure $\mu$, explicitly, if $b_t(x,\mu) = B_t(x, f_\mu(x))$, for $f_\mu := \mathrm{d}\mu/\mathrm{d}x$. Clearly, this is defined only for measures $\mu$ that are absolutely continuous with respect to Lebesgue measure and $\mathrm{d}x$-a.e., with exceptional sets depending on the density $f_\mu$. Note that for such a drift the particle system (1.1), (1.2) is not defined, so the convergence question as we phrased it here does not make sense. Works that consider this type of interaction, yet without common noise, are considered in, e.g., Barbu and Röckner [5] and Grube [25, 26].

Theorem 1.6, and indeed, our entire approach can be generalized in several ways. For one, the boundedness condition on $b$ in Assumption 1.5 (and thus also Assumption 6.2) can be replaced by an appropriate integrability condition. However, this extension introduces additional technicalities that we wish to avoid here. To approach this generalization, it is necessary to extend the emergence of regularity result in [12] to handle an unbounded drift with a suitable integrability condition, as discussed in Section 7 there or in the work by Romito [55]. The precise nature of the integrability condition depends on the Hölder-regularity parameter $\beta$ of $\sigma$, $\bar{\sigma}$, and the dimension $d$. To formulate this condition, it is necessary to specify the values of $r$ and $s$ that appear in Theorem 3.1. This can be achieved by choosing the various parameters explicitly in the proof of Proposition 5.1 and Theorem 1.5 in [12]. This extension is useful for treating another class of singular drifts, such as e.g. $b_t(x,\mu) = (\mu * k_t)(x)$ for a singular kernel $k$. Models of this type and without common noise are considered in, e.g., Krylov and Röckner [41], Röckner and Zhang [53], Röckner et al. [54], Tomaševic [64] or Hao et al. [29].

Finally, it may be possible to replace Assumption 1.5 with other assumptions. However, requires more work.

# A BACKGROUND

## A.1 Functional analysis

**The spaces of Schwartz functions an tempered distributions.** Let $\alpha \in \mathbb{N}_0^d$ be a multi-index, define $|\alpha| = \alpha_1 + \cdots + \alpha_d$, write $\partial^\alpha := \partial_1^{\alpha_1} \cdots \partial_d^{\alpha_d}$ and set

$$\langle x \rangle := (1 + |x|^2)^{1/2} \quad \text{for } x \in \mathbb{R}^d . \tag{A.1}$$

For each $m \in \mathbb{N}_0$ and $\phi \in \mathbf{C}^\infty(\mathbb{R}^d)$, we define the seminorm

$$\|\phi\|_m^* := \max_{|\alpha| \leq m} \sup_{x \in \mathbb{R}^d} |\langle x \rangle^m \partial^\alpha \phi(x)| . \tag{A.2}$$

The *space of rapidly decreasing functions* or Schwartz functions is the vector space

$$\mathcal{S} := \mathcal{S}(\mathbb{R}^d) := \{\phi \in \mathbf{C}^\infty(\mathbb{R}^d) \ : \ \|\phi\|_m^* < \infty \text{ for all } m \in \mathbb{N}_0\}.$$

We endow $\mathcal{S}$ with the Fréchet topology generated by the family $\| \cdot \|_m^*$ for $m \in \mathbb{N}_0$, i.e. $\phi_n \to \phi$ in $\mathcal{S}$ if and only if $\|\phi_n - \phi\|_m^* \to 0$ when $n \to \infty$, for each $m \in \mathbb{N}$. As usual, a set $B \subseteq \mathcal{S}$ is said to be *bounded* if $\sup_{\phi \in B} \|\phi\|_m^* < \infty$ for all $m \in \mathbb{N}_0$.

It is well known that $\mathcal{S}$ is metrizable, complete and separable. Indeed, let $\mathbf{C}_c^\infty = \mathbf{C}_c^\infty(\mathbb{R}^d)$ be the Fréchet space of smooth and compactly supported functions $\mathbb{R}^d \to \mathbb{R}$; see Rudin [56, Sec. 6.2] for details. Then $\mathbf{C}_c^\infty \hookrightarrow \mathcal{S}$, which is to say that $\mathbf{C}_c^\infty$ embeds continuously into $\mathcal{S}$. Moreover, this embedding is dense. As a consequence we also have the continuous embedding $\mathcal{S} \hookrightarrow \mathbb{L}^r(\mathrm{d}x) = \mathbb{L}^r(\mathbb{R}^d, \mathrm{d}x)$ for all $r \in [1, \infty]$, which is dense if $r \in [1, \infty)$.

Let $\lambda : \mathcal{S} \to \mathbb{R}$ be a linear functional. We denote the duality pairing by $\langle \lambda; \phi \rangle_{\mathcal{S}' \times \mathcal{S}}$ for $\phi \in \mathcal{S}$, and usually write $\lambda[\phi] := \langle \lambda; \phi \rangle_{\mathcal{S}' \times \mathcal{S}}$. The functional $\lambda$ is *continuous* if $\phi_n \to \phi$ in $\mathcal{S}$ implies $\lambda[\phi_n] \to \lambda[\phi]$. Equivalently, there exist $C > 0$ and $m \in \mathbb{N}$ such that

$$|\lambda[\phi]| \leq C \|\phi\|_m^* \quad \text{for all } \phi \in \mathcal{S}(\mathbb{R}^d). \tag{A.3}$$

The minimal index $m$ validating this inequality is called the *order* of $\lambda$. Continuous linear functionals are called *(tempered) distributions*; they form the strong topological dual space

$$\mathcal{S}' := \mathcal{S}'(\mathbb{R}^d) := \{\lambda : \mathcal{S}(\mathbb{R}^d) \to \mathbb{R} \ : \ \lambda \text{ is continuous}\}.$$

Its strong topology is generated by uniform convergence on bounded sets, i.e. the topology is generated the family of seminorms $\|\lambda\|_{\mathcal{S}', B} := \sup_{\phi \in B} |\lambda[\phi]|$ as $B$ ranges over the bounded sets of $\mathcal{S}$.

The strong topology generates the Borel-$\sigma$-field on $\mathcal{S}'$, denoted $\mathcal{B}(\mathcal{S}'(\mathbb{R}^d))$ or simply $\mathcal{B}(\mathcal{S}')$. It can be shown that $\mathcal{B}(\mathcal{S}')$ equals the $\sigma$-algebra generated by the collection of sets $\{\lambda \in \mathcal{S}' : \lambda[\phi] < a\}$ for $\phi \in \mathcal{S}$ and $a \in \mathbb{R}$; see e.g. Kallianpur and Xiong [36, Thm. 3.1.1].

Let $f$ be a measurable function. We define

$$\lambda_f[\phi] := \int_{\mathbb{R}^d} \phi(x) f(x)\,\mathrm{d}x\,, \tag{A.4}$$

for $\phi \in \mathcal{S}$, provided the integral exists. Consider $\lambda \in \mathcal{S}'$ and suppose that there exists a $f \in \mathbb{L}^1_{\mathrm{loc}}(\mathrm{d}x)$ such that $\lambda = \lambda_f$ in $\mathcal{S}'$, that is to say $\lambda[\phi] = \lambda_f[\phi]$ for all $\phi \in \mathcal{S}$. Then $f$ is determined uniquely up to $\mathrm{d}x$-a.e. equivalence; see Hörmander [31, Thm. 1.2.5]. In this case it is customary to say that $\lambda$ *is a function*, or that $\lambda$ is representable by $f$. The *canonical embedding* $\mathcal{S} \hookrightarrow \mathcal{S}'$ mapping $f \in \mathcal{S}$ to $\lambda_f \in \mathcal{S}'$ is thus well-defined and, in fact, continuous. This inclusion allows us to identify $\mathcal{S}$ with a subset of $\mathcal{S}'$: If $\lambda \in \mathcal{S}'$ is represented by a function in $\mathcal{S}$, then we freely write $\lambda \in \mathcal{S}$. Similarly, we have the canonical embedding $\mathbb{L}^r(\mathrm{d}x) \hookrightarrow \mathcal{S}'$ for all $r \in [1,\infty]$, and if $\lambda \in \mathcal{S}'$ is represented by a function in $\mathbb{L}^r(\mathrm{d}x)$, we simply write $\lambda \in \mathbb{L}^r(\mathrm{d}x)$. In the same fashion, if $\mu$ is a measure on the Borel sets of $\mathbb{R}^d$, we may define

$$\lambda_\mu[\phi] := \int_{\mathbb{R}^d} \phi(x)\,\mu(\mathrm{d}x) \tag{A.5}$$

for $\phi \in \mathcal{S}$, provided the integral exists. In analogy to the above, we say that $\lambda \in \mathcal{S}'$ is a measure if there is a measure $\mu$ such that $\lambda = \lambda_\mu$ in $\mathcal{S}'$. Let $\mathbf{M}_1^+ = \mathbf{M}_1^+(\mathbb{R}^d)$ denote the set of probability measures on the Borel sets of $\mathbb{R}^d$. It is standard to verify that $\mathbf{M}_1^+ \subseteq \mathcal{S}'$ by virtue of (A.5). In fact, the inclusion is sequentially continuous if $\mathbf{M}_1^+$ carries the narrow topology; see Lemma A.3 below.

We recall some basics from functional analysis, notably the space of tempered distributions. The material in this section is standard and can be found, e.g., in Rudin [56, Ch. 6] and Simon [59, Ch. 6].

**Identification of $\mathcal{S}'$ with a sequence space** The canonical embedding (A.4) lets us write

$$\mathcal{S}(\mathbb{R}^d) \hookrightarrow \mathbb{L}^2(\mathbb{R}^d) \hookrightarrow \mathcal{S}'(\mathbb{R}^d)\,.$$

This allows to give an elegant and useful identification of $\mathcal{S}(\mathbb{R}^d)$ with rapidly decreasing sequences, and dually of $\mathcal{S}'(\mathbb{R}^d)$ with slowly increasing sequences, which we describe by using Hermite functions. An efficient presentation of this classical identification is found in Simon [58], to which we also refer for more background.

Starting from the standard Gaussian density

$$g(x) := \frac{1}{\sqrt{2\pi}} \exp\Big(-\frac{1}{2}x^2\Big) \qquad \text{for } x \in \mathbb{R}\,,$$

we define the *Hermite functions* on $\mathbb{R}$ by

$$h_{k+1}(x) := \frac{(-1)^k}{\sqrt{k!}} g(x)^{-\frac{1}{2}} \Big(\frac{\mathrm{d}}{\mathrm{d}x}\Big)^k g(x) \qquad \text{for } k = 0,1,2,\ldots$$

The above definition extends to $\mathbb{R}^d$ as follows. We define the Hermite functions on $\mathbb{R}^d$ by $h_k(x) := h_{k_1}(x_1)\cdots h_{k_d}(x_d)$ for each multi-index $k \in \mathbb{N}^d$ and $x = (x_1,\ldots,x_d) \in \mathbb{R}^d$. The Hermite functions are elements of $\mathcal{S}$, i.e. $h_k \in \mathcal{S}$ for any $k \in \mathbb{N}^d$. Moreover, they are uniformly bounded, i.e. $\sup_{k\in\mathbb{N}^d} \|h_k\|_\infty \le (2\pi)^{d/4}$; see Indritz [34, p. 981] and Simon [58, Sec. 3, Note 12].

Let $(\,\cdot\,,\,\cdot\,)$ denote the usual inner product on $\mathbb{L}^2(\mathrm{d}x)$, which we use to identify $\mathbb{L}^2(\mathrm{d}x)$ with its dual via the Riesz isomorphism. It is well known that the family $(h_k)_{k\in\mathbb{N}^d}$ is a complete orthonormal system in $\mathbb{L}^2(\mathbb{R}^d, \mathrm{d}x)$ with respect to the usual inner product; see e.g. Simon [59, Thm. 6.4.3]

Since $\mathcal{S} \subseteq \mathbb{L}^2(\mathrm{d}x)$, the basis expansion

$$\phi = \sum_{k\in\mathbb{N}^d} (\phi; h_k)\, h_k \qquad \text{in } \mathbb{L}^2(dx)$$

is well defined for any $\phi \in \mathcal{S}$, and we write $\phi_k^\# := (\phi; h_k)$ with $k \in \mathbb{N}^d$ for the *Hermite coefficients* of $\phi$. Similarly, if $\lambda \in \mathcal{S}'$, we define the Hermite coefficient by $\lambda_k^\# := \lambda[h_k]$, which makes sense since $h_k \in \mathcal{S}$.

We are now ready to discuss the identification of $\mathcal{S}$ and $\mathcal{S}'$ with certain sequences as promised at the beginning of the subsection. We call a sequence $(a_k)_{k\in\mathbb{N}^d}$ of real numbers *rapidly decreasing* if $(\langle k\rangle^m a_k)_{k\in\mathbb{N}^d} \in \ell^2(\mathbb{N}^d)$ for each $m \in \mathbb{N}_0$, with $\langle k\rangle = (1 + |k|^2)^{1/2}$ as in (A.1). We call $(a_k)_{k\in\mathbb{N}^d}$ *slowly increasing* if there exist $m \in \mathbb{N}_0$ and $c_a < \infty$ such that $|a_k| \le c_a \langle k\rangle^m$ for all $k \in \mathbb{N}^d$. Clearly, a rapidly decreasing sequence is slowly increasing.

**Lemma A.1** | 1) *Let $\phi \in \mathcal{S}(\mathbb{R}^d)$. Then $(\phi_k^\#)_{k\in\mathbb{N}^d}$ is rapidly decreasing, i.e. $(\langle k\rangle^m \phi_k^\#)_{k\in\mathbb{N}^d} \in \ell^2(\mathbb{N}^d)$ for any $m \in \mathbb{N}_0$. Conversely, any rapidly decreasing sequence $(a_k)_{k\in\mathbb{N}^d}$ defines a function $\psi \in \mathcal{S}(\mathbb{R}^d)$ by virtue of $\psi = \sum_{k\in\mathbb{N}^d} a_k h_k$ with convergence in the topology of $\mathcal{S}(\mathbb{R}^d)$.*

2) *Let $\lambda \in \mathcal{S}'(\mathbb{R}^d)$. Then $(\lambda_k^\#)_{k\in\mathbb{N}^d}$ is slowly increasing, i.e. there are $m \in \mathbb{N}_0$ and $c_\lambda < \infty$ such that $|\lambda_k^\#| \le c_\lambda \langle k\rangle^m$ for $k \in \mathbb{N}^d$. Conversely, any slowly increasing sequence $(b_k)_{k\in\mathbb{N}^d}$ defines a distribution $\tau \in \mathcal{S}'$ by virtue of $\tau[\phi] := \sum_{k\in\mathbb{N}^d} \phi_k^\# b_k$ in $\mathbb{R}$, for each $\phi \in \mathcal{S}$.*

***Proof*** See e.g. [58, Thms. 1 and 2]. □

By combining parts 1) and 2) of Lemma A.1, it can be shown that for any $\lambda \in \mathcal{S}'$ and $\phi \in \mathcal{S}$, we have

$$\lambda[\phi] = \sum_{k\in\mathbb{N}^d} \lambda_k^\# \phi_k^\# \,. \tag{A.6}$$

Indeed, the sum in (A.6) converges because $(\phi_k^\#)_{k\in\mathbb{N}^d}$ is rapidly decreasing and $(\lambda_k^\#)_{k\in\mathbb{N}^d}$ is slowly increasing; see [58, Thm. 3].

**Intermediate Hilbert spaces** Let $\lambda \in \mathcal{S}'$. By Lemma A.1, the Hermite coefficients $(\lambda_k^\#)_{k\in\mathbb{N}^d}$ are slowly increasing. There thus exists a $p \in \mathbb{R}$ such that

$$\|\lambda\|_{\mathscr{H}_p} := \left( \sum_{k\in\mathbb{N}^d} \left(\langle k\rangle^{p/d} \lambda_k^\#\right)^2 \right)^{\frac{1}{2}} \tag{A.7}$$

is finite. Moreover, if $\phi \in \mathcal{S}$, then $\|\phi\|_{\mathscr{H}_p}$ is finite for any $p \in \mathbb{R}$, since the Hermite coefficients $(\phi_k^\#)_{k\in\mathbb{N}^d}$ are rapidly decreasing.

For each $p \in \mathbb{R}$, we thus define

$$\mathscr{H}_p := \mathscr{H}_p(\mathbb{R}^d) := \{ f \in \mathcal{S}' : \|f\|_{\mathscr{H}_p} < \infty \} \,.$$

Since $\mathcal{S}$ is $\sigma(\mathcal{S}', \mathcal{S})$-dense in $\mathcal{S}'$, the space $(\mathscr{H}_p, \|\cdot\|_{\mathscr{H}_p})$ can equivalently be defined as the Hilbert space completion of $\mathcal{S}$ under the semi-norm $\|\cdot\|_{\mathscr{H}_p}$; so $\mathcal{S}$ is dense in $\mathscr{H}_p$ for any $p \in \mathbb{R}$.

On $\mathscr{H}_p$, (A.7) specifies a Hilbertian semi-norm, that is, a semi-norm satisfying the parallelogram identity. The *Hermite–Fourier space* $(\mathscr{H}_p, \|\cdot\|_{\mathscr{H}_p})$ is thus a Hilbert subspace of $\mathcal{S}'$, with semi-inner product $(\,\cdot\,;\cdot\,)_p$ obtained from $\|\cdot\|_{\mathscr{H}_p}$ by polarization. With this notation, we observe that $(\,\cdot\,;\cdot\,) = (\,\cdot\,;\cdot\,)_0$.

Let $\mathscr{H}_p'$ denote the (strong) topological dual of $\mathscr{H}_p$. By definition, it is equipped with the norm

$$\|\lambda\|_{\mathscr{H}_p'} := \sup\{|\lambda[\phi]| : \|\phi\|_{\mathscr{H}_p} \le 1\}\,. \tag{A.8}$$

The pair $(\mathscr{H}_p', \|\cdot\|_{\mathscr{H}_p'})$ is again a Hilbert space.

It is classical that the canonical inclusion $\mathcal{S}(\mathbb{R}^d) \hookrightarrow \mathscr{H}_p'$ defined by $f \mapsto \lambda_f[\,\cdot\,] = (\,\cdot\,, f)_0$ extends uniquely to a continuous linear map $J : \mathscr{H}_{-p}(\mathbb{R}^d) \to \mathscr{H}_p'(\mathbb{R}^d)$ which is an isometric isomorphism so that $\mathscr{H}_{-p} \simeq \mathscr{H}_p'$; see for instance Kallianpur and Xiong [36, Thm. 1.3.2(c)] and Demengel and Demengel [19, Prop. 4.10].

An important consequence of the above is the chain of continuously embedded spaces

$$\mathcal{S}(\mathbb{R}^d) \hookrightarrow \mathscr{H}_p(\mathbb{R}^d) \hookrightarrow \mathbb{L}^2(\mathbb{R}^d) \hookrightarrow \mathscr{H}_{-p}(\mathbb{R}^d) \hookrightarrow \mathcal{S}'(\mathbb{R}^d)\,. \tag{A.9}$$

We call the parameter $p \ge 0$ the *regularity*, and note that $\mathbb{L}^2(\mathrm{d}x) = \mathscr{H}_0(\mathbb{R}^d)$ with equal norms. The terminology we employ is suggestive: It indicates that objects of negative regularity act on objects of positive regularity by *canonical duality*

From the definition (A.7), we see that $\|\cdot\|_{\mathscr{H}_{p_1}} \le \|\cdot\|_{\mathscr{H}_{p_2}}$ for $p_1 \le p_2$. Thus $\mathscr{H}_{p_2} \hookrightarrow \mathscr{H}_{p_1}$ whenever $p_1 \le p_2$, and dually, $\mathscr{H}_{-p_1} \hookrightarrow \mathscr{H}_{-p_2}$; so the chain (A.9) can be expanded indefinitely. In fact, from the abstract definitions, we can obtain

$$\mathcal{S}(\mathbb{R}^d) = \bigcap_{p\in\mathbb{R}} \mathscr{H}_p(\mathbb{R}^d) \quad \text{and} \quad \mathcal{S}'(\mathbb{R}^d) = \bigcup_{p\in\mathbb{R}} \mathscr{H}_p'(\mathbb{R}^d)\,. \tag{A.10}$$

Both equalities in (A.10) are to be understood not just as equalities of sets, but as equalities of topological spaces. For this, we endow $\bigcap_{p\in\mathbb{R}} \mathscr{H}_p(\mathbb{R}^d)$ with the limit topology and $\bigcup_{p\in\mathbb{R}} \mathscr{H}_p(\mathbb{R}^d)$ with the colimit topology. In other words, if $(\phi_n)_{n\in\mathbb{N}}$ is a sequence in $\mathcal{S}$, then $\phi_n \to \phi$ in $\mathcal{S}$ as $n \to \infty$ if and only if $\phi_n \to \phi$ as $n \to \infty$ in $\mathscr{H}_p$ for any $p \in \mathbb{R}$. Dually, if $(\lambda_n)_{n\in\mathbb{N}}$ is a sequence in $\mathcal{S}'$, then $\lambda_n \to \lambda$ as $n \to \infty$ in $\mathcal{S}'$ if and only if $\lambda_n \to \lambda$ as $n \to \infty$ in $\mathscr{H}_p$ for some $p \in \mathbb{R}$.

**Lemma A.2** | *The two families of semi-norms $\|\cdot\|_{\mathscr{H}_p}$ in* (A.7) *for $p \in \mathbb{N}_0$ and $\|\cdot\|_m^*$ in* (A.2) *for $m \in \mathbb{N}_0$ on $\mathcal{S}(\mathbb{R}^d)$ are equivalent. That is, for each $m$, there are a $p$ and a universal constant $c_m$ such that*

$$\|\phi\|_m^* \le c_m \|\phi\|_{\mathscr{H}_p}\,,$$

*and for each $p$, there are an $m$ and a universal constant $c_p$ such that*

$$\|\phi\|_{\mathscr{H}_p} \le c_p \|\phi\|_m^*\,,$$

*both for all $\phi \in \mathcal{S}(\mathbb{R}^d)$. For each $m$, the number $p$ and the constant $c_m$ depend on $d$, and similarly for each $p$, so do $m$ and $c_p$.*

***Proof*** See e.g. Holley and Stroock [30, Appendix, (A.14)–(A.20)]. □

In view of the definition of $\|\cdot\|_m^*$ in (A.2), the above lemma shows in particular that functions in $\mathscr{H}_p$ enjoy certain differentiability properties as well as polynomial decay up to a certain order. An immediate consequence of this is in particular that $\mathscr{H}_p(\mathbb{R}^d) \hookrightarrow \mathsf{Lip}(\mathbb{R}^d)$, the space of Lipschitz functions on $\mathbb{R}^d$, for some $p > 0$ large enough depending on $d$.

## A.2 Probability theory

**Spaces of probability measures.** A special role is our study is taken by positive measures, and in particular probability measures. We thus take a moment to clarify how they relate to our exposition so far. For a topological space X, we write $\mathbf{M}_1^+(\mathrm{X})$ for the set of probability measures on $\mathcal{B}(X)$. If $\mathrm{X} = \mathbb{R}^d$, then we oftentimes write $\mathbf{M}_1^+ := \mathbf{M}_1^+(\mathbb{R}^d)$.

The *narrow topology* $\tau_{\mathrm{wk}^*}$ is induced by duality with continuous bounded functions $\mathbf{C}_\mathrm{b}(\mathrm{X})$. Specifically, if $(\lambda_m)_{m\in\mathbb{N}}$ is a sequence in $\mathbf{M}_1^+(\mathrm{X})$, then $\lambda_m \to \lambda$ as $m \to \infty$ relative to $\tau_{\mathrm{wk}^*}$ if $\lambda_m[f] \to \lambda[f]$ for all $f \in \mathbf{C}_\mathrm{b}(\mathrm{X})$. We write $\mathcal{P}_{\mathrm{wk}^*}(\mathrm{X}) := (\mathbf{M}_1^+(\mathrm{X}), \tau_{\mathrm{wk}^*})$ for the space $\mathbf{M}_1^+(\mathrm{X})$ endowed with the narrow topology $\tau_{\mathrm{wk}^*}(\mathrm{X})$. If X is Polish, then so is $\mathcal{P}_{\mathrm{wk}^*}(\mathrm{X})$; see, e.g. Aliprantis and Border [2, Thm. 15.15].

Since $\mathbf{M}_1^+ \subseteq \mathcal{S}'$, we can consider the trace topology of $\mathcal{S}'$ on $\mathbf{M}_1^+$. This topology is, however, different from narrow convergence. To discuss this difference, denote by $\mathbf{M}_{\le 1}^+ = \mathbf{M}_{\le 1}^+(\mathbb{R}^d)$ the space of subprobability measures on $\mathcal{B}(\mathbb{R}^d)$, i.e. positive measures of no more than unit mass. Evidently $\mathbf{M}_1^+ \subseteq \mathbf{M}_{\le 1}^+ \subseteq \mathcal{S}'$. The set $\mathbf{M}_1^+$ is *not closed* in $\mathcal{P}_{\le 1,\mathcal{S}'}$; see e.g. [12, Rem. A.2]. The set $\mathbf{M}_{\le 1}^+$ equipped with the trace topology of $\mathcal{S}'$ is denoted by $\mathcal{P}_{\le 1,\mathcal{S}'} := (\mathbf{M}_{\le 1}^+, \tau_{\mathcal{S}'})$. Specifically, if $(\lambda_m)_{m\in\mathbb{N}}$ is a sequence in $\mathbf{M}_{\le 1}^+$, then $\lambda_m \to \lambda$ as $m \to \infty$ in $\mathcal{P}_{\le 1,\mathcal{S}'}$ if for all $\phi \in \mathcal{S}$, we have that $\lambda_m[\phi] \to \lambda[\phi]$.

**Lemma A.3** | *If $(\nu_m)_{m\in\mathbb{N}}$ is a sequence in $\mathbf{M}_1^+$, then:*
*1) If $\nu_m \to \nu$ in $\mathcal{P}_{\mathrm{wk}^*}$, then $\nu \in \mathbf{M}_1^+$ and $\nu_m \to \nu$ in $\mathcal{P}_{\mathcal{S}'}$.*
*2) If $\nu_m \to \nu$ in $\mathcal{P}_{\le 1,\mathcal{S}'}$ and $\nu[\mathbb{R}^d] = 1$, then $\nu_m \to \nu$ in $\mathcal{P}_{\mathrm{wk}^*}$.*

***Proof*** Part 1) follows from the fact that $\mathcal{S} \subseteq \mathbf{C}_b$ and that the notions of weak and strong convergence in $\mathcal{S}'$ coincide for sequences; see Huang and Yan [32, Thm. 3.12]. Part 2) follows from classical properties of narrow convergence; see Klenke [39, Thm. 13.16]. □

# B OMITTED PROOFS

I thank Martin Schweizer for generalizing and streamlining an initial version of this proof, and for adding the reference to the work of Stricker [60].

***Proof of Proposition 2.6*** The proof is subdivided into several steps.

***Step 1*** Assumptions 1.5 and 1.2 give us the weak solution $\boldsymbol{X}^n$ of (1.1), (1.2) on $(\Omega^n_{\mathrm{par}}, \mathcal{G}^n, \mathbb{P}^n_{x_0})$ via Lemma 1.3. From Lemma 1.4, we know that $\mu^n$ is $\mathcal{S}'$-valued with continuous trajectories. Moreover, for each $\phi \in \mathcal{S}$, the calculations using Itô's formula carried out in (2.15)–(2.18) give us the representation

$$\mu_t^n[\phi] = \delta_{x_0}[\phi] + \int_0^t \mu_s^n[\mathscr{L}_s(\mu_s^n)\phi]\,\mathrm{d}s + N_t^{n,\phi} \qquad t \in [0,T]\,, \tag{B.1}$$

where $N^{n,\phi}$ is the continuous $(\mathbb{G}^n, \mathbb{P}^n_{x_0})$-local martingale given, in view of (2.19), (2.20) and (B.1), (2.18), by

$$\begin{aligned} N_t^{n,\phi} &:= \frac{1}{n}\sum_{i=1}^n (M^\phi)_t^{i,n} + \frac{1}{n}\sum_{i=1}^n (\bar{M}^\phi)_t^{i,n} \\ &= \frac{1}{n}\sum_{i=1}^n \int_0^t \mathscr{U}_s(\mu_s^n)\phi(X_t^{i,n})\,\mathrm{d}B_s^{i,n} + \int_0^t \mu_s^n[\mathscr{R}_s(\mu_s^n)\phi]\,\mathrm{d}Z_s^n\,. \end{aligned} \tag{B.2}$$

The real-valued process $\mu^n[\phi]$ is thus for each $\phi \in \mathcal{S}$ a continuous $\mathbb{R}$-valued semimartingale, and the regularization result from Pérez-Abreu [49, Thm. 1] (see also part 1) of Lemma 2.2) shows that $\mu^n$ is a continuous $\mathcal{S}'$-valued semimartingale relative to $(\mathbb{G}^n, \mathbb{P}^n_{x_0})$ in the sense of Definition 2.1.

***Step 2*** Our next aim is to obtain a representation in terms of the canonical process $\Lambda$ on $(\Omega, \mathcal{F})$, and this is more subtle. First, $\Lambda$ has under the pushed measure $\mathbf{P}^n_{x_0}$ from (1.5) the same law as $\mu^n$ under $\mathbb{P}^n_{x_0}$. However, the filtration $\mathbb{F}$ on $(\Omega, \mathcal{F})$ is generated by $\Lambda$ (and made right-continuous), whereas the filtration $\mathbb{G}^n$ on $(\Omega^n_{\mathrm{par}}, \mathcal{G}^n)$ is generated by $(X^n, Z^n)$ (and made right-continuous). To obtain a match between $\Lambda$ and $\mu^n$ as processes, let $\tilde{\mathbb{G}}^n := (\tilde{\mathcal{G}}^n_t)_{t\in[0,T]}$ with $\tilde{\mathcal{G}}^n_t := \bigcap_{0<\varepsilon<T-t} \sigma(\mu^n_s\,;\, 0 \le s < t+\varepsilon)$ for $t \in [0,T)$ and $\tilde{\mathcal{G}}^n_T := \sigma(\mu^n_s\,;\, 0 \le s \le T)$ be the right-continuous version of the raw filtration generated by $\mu^n$, viewed as process with values in $(\mathcal{S}', \mathcal{B}(\mathcal{S}'))$. Then (1.1), (1.2) and (2.4) imply $\tilde{\mathbb{G}}^n \subseteq \mathbb{G}^n$, and rewriting (B.1) as $N_t^{n,\phi} = \mu_t^n[\phi] - \delta_{x_0}[\phi] - \int_0^t \mu_s^n[\mathscr{L}_s(\mu_s^n)\phi]\,\mathrm{d}s$, $t \in [0,T]$, shows that $N^{n,\phi}$ is $\tilde{\mathbb{G}}^n$-adapted and hence also a $\tilde{\mathbb{G}}^n$-semimartingale by Stricker [60, Thm. 3.1]. But $N^{n,\phi}$ is continuous; so it is even a special $\tilde{\mathbb{G}}^n$-semimartingale by [60, Thm. 2.6]. Moreover, the canonical decomposition of $\mu^n$ relative to $(\tilde{\mathbb{G}}^n, \mathbb{P}^n_{x_0})$ is the same as relative to $(\mathbb{G}^n, \mathbb{P}^n_{x_0})$. Since $\mathrm{Law}_{\mathbb{P}^n_{x_0}}(\mu^n) = \mathrm{Law}_{\mathbf{P}^n_{x_0}}(\Lambda)$ by (1.5), and the filtrations $\mathbb{F}$ for $\Lambda$ and $\tilde{\mathbb{G}}^n$ for $\mu^n$ match by construction, we see that for each fixed $\phi \in \mathcal{S}$, the process $\Lambda[\phi]$ satisfies in analogy to (B.1) the dynamics

$$\Lambda_t[\phi] = \delta_{x_0}[\phi] + \int_0^t \Lambda_s[\mathscr{L}_s(\Lambda_s)\phi]\,\mathrm{d}s + M_t^\phi\,, \qquad t \in [0,T]\,, \tag{B.3}$$

with $\mathrm{Law}_{\mathbf{P}^n_{x_0}}(M^\phi) = \mathrm{Law}_{\mathbb{P}^n_{x_0}}(N^{n,\phi})$ so that $M^\phi$ is a continuous local $(\mathbb{F}, \mathbf{P}^n_{x_0})$-martingale. Note that we can also view (B.3) as the definition of $M^\phi$.

***Step 3*** For the canonical decomposition (B.3) of $\Lambda[\phi]$ for fixed $\phi \in \mathcal{S}$, we note that $\Lambda[\phi]$ is bounded uniformly in $(t,\omega) \in [0,T] \times \Omega$ because $\phi$ is bounded and $\Lambda$ takes values in $\mathbf{M}_1^+$ $\mathbf{P}^n_{x_0}$-a.s. Moreover, the finite-variation term $A^\phi$ in (B.3) is absolutely continuous with drift

$$|\Lambda_s[\mathscr{L}_s(\Lambda_s)\phi]| \le \int_{\mathbb{R}^d} |\mathscr{L}_s(\Lambda_s)\phi(x)|\,\Lambda_s(\mathrm{d}x)\,.$$

But by the definition (2.8) of $\mathscr{L}_s$ and because $b$, $\sigma$, $\bar{\sigma}$ are bounded by Assumptions 1.5 and 1.2,

$$|\mathscr{L}_s(\Lambda_s)\phi(x)| \le \Big(\|b\|_\infty + \frac{1}{2}(\|\sigma\|_\infty^2 + \|\bar{\sigma}\|_\infty^2)\Big)\|\phi\|_2^*\,,$$

where $\|\phi\|_2^* = \max_{|\alpha|\le 2} \sup_{x\in\mathbb{R}^d} |(1+|x|^2)\partial^\alpha \phi(x)| < \infty$ was defined in (A.2). Again using that $\Lambda$ takes values in $\mathbf{M}_1^+$ $\mathbf{P}^n_{x_0}$-a.s., we conclude that $\mathbf{P}^n_{x_0}$-a.s.,

$$\int_0^T |\mathrm{d}A_s^\phi| = \int_0^T |\Lambda_s[\mathscr{L}_s(\Lambda_s)\phi(x)]|\,\mathrm{d}s \le T\Big(\|b\|_\infty + \frac{1}{2}(\|\sigma\|_\infty^2 + \|\bar\sigma\|_\infty^2)\Big)\|\phi\|_2^* \tag{B.4}$$

so that the variation of $A^\phi$ is bounded uniformly in $(t,\omega)\in[0,T]\times\Omega$. In consequence,

$$\begin{aligned}\sup_{t\in[0,T]} |M_t^\phi| &= \sup_{t\in[0,T]} |\Lambda_t[\phi] - \delta_{x_0} - A_t^\phi| \\ &\le 2\|\phi\|_\infty + \int_0^T |\mathrm{d}A_s^\phi| \\ &\le \bigg(2 + T\Big(\|b\|_\infty + \frac{1}{2}(\|\sigma\|_\infty^2 + \|\bar\sigma\|_\infty^2)\Big)\bigg)\|\phi\|_2^* \end{aligned}\tag{B.5}$$

shows that the local martingale $M^\phi$ is in fact for each $\phi\in\mathcal{S}$ a bounded continuous $(\mathbb{F},\mathbf{P}^n_{x_0})$-martingale.

***Step 4*** By Steps 1 and 2, $\Lambda$ is an $\mathcal{S}'$-valued $(\mathbb{F},\mathbf{P}^n_{x_0})$-semimartingale; so we can rewrite (B.3) as

$$\Lambda_t = \delta_{x_0} + A_t + M_t\,, \qquad t\in[0,T]\,,$$

with the $\mathcal{S}'$-valued processes $A$, $M$ defined using (2.26) by

$$\begin{aligned} A_t[\phi] &\coloneqq A_t^\phi = \int_0^t \Lambda_s[\mathscr{L}_s(\Lambda_s)\phi]\,\mathrm{d}s = \int_0^t \mathscr{A}_s(\Lambda_s)[\phi]\,\mathrm{d}s\,, \\ M_t[\phi] &\coloneqq M_t^\phi\,. \end{aligned}\tag{B.6}$$

Then the boundedness results in Step 2 allow us to use part 2) of Lemma 2.2 and conclude that $M\in\mathscr{M}_2^c(\mathbb{F},\mathbf{P}^n_{x_0};\mathcal{S}')$ and $\Lambda, A, M\in\mathbf{C}([0,T];\mathscr{H}_{-p})$ $\mathbf{P}^n_{x_0}$-a.s. for any $p\ge p_0$ for some $p_0>0$. At this point, $p_0$ can still depend on $n$, but we show in Step 6 below that this dependence can be eliminated. Note also that using the standard linearity property of the Bochner integral, e.g. [33, Eqn. (1.2)], we can write

$$A_t[\phi] = \int_0^t \mathscr{A}_s(\Lambda_s)[\phi]\,\mathrm{d}s = \Big(\int_0^t \mathscr{A}_s(\Lambda_s)\,\mathrm{d}s\Big)[\phi]\,,$$

where the first integral is a Lebesgue integral and the second a Bochner integral. This establishes (2.30).

***Step 5*** To establish the representation (2.30) of $\langle\langle M\rangle\rangle$, we can exploit earlier computations. Due to $\mathrm{Law}_{\mathbb{P}^n_{x_0}}(\mu^n) = \mathrm{Law}_{\mathbf{P}^n_{x_0}}(\Lambda)$ and the correspondence between the filtrations $\mathbb{F}$ and $\tilde{\mathbb{G}}^n$ from Step 2, we can obtain $\langle M^\phi, M^\psi\rangle$ under $\mathbf{P}^n_{x_0}$ by computing $\langle N^{n,\phi}, N^{n,\psi}\rangle$ under $\mathbb{P}^n_{x_0}$ and then replacing $\mu^n$ by $\Lambda$ in the resulting formulas. But by (B.1), $N^{n,\phi}$ is the local martingale in the canonical decomposition of the continuous $(\tilde{\mathbb{G}}^n,\mathbb{P}^n_{x_0})$-semimartingale

$\mu^n[\phi]$, and so (2.25) and (2.27), (2.28) yield

$$\begin{aligned}\langle N^{n,\phi}, N^{n,\psi}\rangle_t &= \langle \mu^n[\phi], \mu^n[\psi]\rangle_t \\ &= \int_0^t \Big(\mu^n_s[\mathscr{R}_s(\mu^n_s)\phi]\Big)\cdot\Big(\mu^n_s[\mathscr{R}_s(\mu^n_s)\psi]\Big)\,\mathrm{d}s \\ &\quad + \frac{1}{n}\int_0^t \mu^n_s\Big[\Big(\mathscr{U}_s(\mu^n_s)\phi\Big)\cdot\Big(\mathscr{U}_s(\mu^n_s)\psi\Big)\Big]\,\mathrm{d}s \\ &= \int_0^t \mathscr{Q}_s(\mu^n_s)[\phi,\psi]\,\mathrm{d}s + \frac{1}{n}\int_0^t \mathscr{C}_s(\mu^n_s)[\phi,\psi]\,\mathrm{d}s\,.\end{aligned}$$

Therefore (B.6) and standard linearity properties of the Bochner integral imply that for $\mathbb{F}$ and $\mathbf{P}^n_{x_0}$,

$$\begin{aligned}\langle M[\phi], M[\psi]\rangle_t &= \langle M^\phi, M^\psi\rangle_t \\ &= \int_0^t \mathscr{Q}_s(\Lambda_s)[\phi,\psi]\,\mathrm{d}s + \frac{1}{n}\int_0^t \mathscr{C}_s(\Lambda_s)[\phi,\psi]\,\mathrm{d}s \\ &= \Big(\int_0^t \mathscr{Q}_s(\Lambda_s)\,\mathrm{d}s + \frac{1}{n}\int_0^t \mathscr{C}_s(\Lambda_s)\,\mathrm{d}s\Big)[\phi,\psi]\,,\end{aligned}$$

and because

$$\langle\!\langle M\rangle\!\rangle[\phi,\psi] = \langle M[\phi], M[\psi]\rangle\,.$$

Since this uniquely characterizes the tensor quadratic variation, we obtain (2.30); see Section 2.1.

***Step 6*** It remains to show that the dependence in Step 4 of $p_0$ on $n$ can be removed, i.e. that we can find $p_0 > 0$ with $\Lambda \in \mathbf{C}([0,T];\,\mathscr{H}_{-p})$ $\mathbf{P}^n_{x_0}$-a.s. for all $n \in \mathbb{N}$ and also $M \in \mathscr{M}^{\mathrm{c}}_2(\mathbb{F}, \mathbf{P}^n_{x_0};\,\mathscr{H}_{-p})$ for all $n \in \mathbb{N}$, both for $p \geq p_0$. To this end, we first note that for any $n \in \mathbb{N}$, $\Lambda$ is by the construction of $\mathbf{P}^n_{x_0}$ in (1.5) $\mathbf{P}^n_{x_0}$-a.s. $\mathbf{M}^+_1$-valued. We can thus fix $m \in \mathbb{N}$ and use Lemma A.2 to obtain $r > 0$ such that

$$|\Lambda_t[\phi]| \leq \|\phi\|_\infty \leq \|\phi\|^*_m \leq c_m\|\phi\|_{\mathscr{H}_r}\,.$$

Via (A.8) and using the fact $\mathscr{H}_{-p}$ and $\mathscr{H}'_p$ are isometrically isomorphic by virtue of the map $J$ from Section 2.1, this yields

$$\sup_{t\in[0,T]}\|\Lambda_t\|_{\mathscr{H}_{-r}} = \sup_{t\in[0,T]}\|\Lambda_t\|_{\mathscr{H}'_r} = \sup_{t\in[0,T]}\sup\{|\Lambda_t[\phi]| \,:\, \|\phi\|_{\mathscr{H}_r} \leq 1\} \leq c_m$$

$\mathbf{P}^n_{x_0}$-a.s. for all $n \in \mathbb{N}$. By part 1) of Lemma 1.4 , the sequence $(\mathbf{P}^n_{x_0})_{n\in\mathbb{N}}$ is tight. Together with the above bound, this allows us to apply Kallianpur and Xiong [36, Thm. 2.5.2] and obtain $p_1 > 0$ such that

$$\Lambda \in \mathbf{C}([0,T];\,\mathscr{H}_{-p}) \subseteq \mathbf{C}([0,T];\,\mathscr{H}_{-p_1})$$

$\mathbf{P}^n_{x_0}$-a.s. for all $n \in \mathbb{N}$ and all $p \geq p_1$. This already yields the first statement in (2.31) as $p_1$ does not depend on $n$. Next, thanks to the estimate (B.5) for $M_t[\phi] = M^\phi_t$, a completely analogous argument for $M$ instead of $\Lambda$ shows that

$$M \in \mathbf{C}([0,T];\,\mathscr{H}_{-p}) \subseteq \mathbf{C}([0,T];\,\mathscr{H}_{-p_2})$$

$\mathbf{P}^n_{x_0}$-a.s. for all $n \in \mathbb{N}$ and all $p \geq p_2 > 0$. Finally, repeating this argument for $t = T$ shows that

$$\mathbf{E}^{\mathbf{P}^n_{x_0}}[\|M_T\|^2_{\mathscr{H}_{-p}}] < \infty \qquad \text{for all } n \in \mathbb{N} \text{ and all } p \geq p_2 \,.$$

As a matter of fact, we even have a bound uniformly in $\omega$. From Equation (2.1), this shows that $M \in \mathscr{M}^c_2(\mathbb{F}, \mathbf{P}^n_{x_0}; \mathscr{H}_{-p})$ for all $n \in \mathbb{N}$ and all $p \geq p_2$. Taking $p_0 := \max\{p_1, p_2\}$ thus completes the proof. □

***Proof of Lemma 4.3*** Measurability is established in e.g. Makarov and Podkorytov [45, Lem. 7.5.2]. For part 1), see e.g. Brézis [8, Prop. 4.20].

For parts 2)–4), we set $f^\delta(x) := \int_{\mathbb{R}^d} f(y) h_\delta(x - y)\, \mathrm{d}y$ for bounded and measurable $f : \mathbb{R}^d \to \mathbb{R}$ and observe the two following facts. First, defining $K'$ as in part 4), it can be verified directly from the definitions that $f^\delta(x) = (f \mathbb{1}_{K'})^\delta(x)$ for all $x \in K$; see for instance Makarov and Podkorytov [45, Lem. 7.5.4]. Second, the following version of Young's inequality holds (see e.g. Adams and Fournier [1, Thm. 2.24]): For all $r_1, r_2, r_3 \in [1, \infty]$ satisfying $1/r_1 + 1/r_2 + 1/r_3 = 2$ and all $w \in \mathbb{L}^{r_3}(\mathbb{R}^d, \mathrm{d}x)$, we have

$$\left| \int_{\mathbb{R}^d} f^\delta(x) w(x)\, \mathrm{d}x \right| \leq \|f\|_{\mathbb{L}^{r_1}(\mathbb{R}^d, \mathrm{d}x)} \|h\|_{\mathbb{L}^{r_2}(\mathbb{R}^d, \mathrm{d}x)} \|w\|_{\mathbb{L}^{r_3}(\mathbb{R}^d, \mathrm{d}x)} \,.$$

We can now prove part 2). Indeed, setting $f = b_t(\lambda)\mathbb{1}_{K'}$, which is bounded and measurable by Assumption 1.5, we have by [8, Thm. 4.22] that $\lim_{\delta \to 0} \|f - f^\delta\|_{\mathbb{L}^r(\mathbb{R}^d, \mathrm{d}x)} = 0$ whenever $r' < \infty$. Thus also $\lim_{\delta \to 0} \|f - f^\delta\|_{\mathbb{L}^r(K, \mathrm{d}x)} = 0$. But since $f = b_t(\lambda)\mathbb{1}_{K'}$, we have by the first fact observed above that

$$\left\| b_t(\lambda)\mathbb{1}_{K'} - \left(b_t(\lambda)\mathbb{1}_{K'}\right)^\delta \right\|_{\mathbb{L}^r(K, \mathrm{d}x)} = \|b_t(\lambda) - b_t(\lambda)^\delta\|_{\mathbb{L}^r(K, \mathrm{d}x)} \,.$$

So part 2) follows.

For part 3), note that Hölder's inequality and $\|h_\delta\|_{\mathbb{L}^1} = 1$ show that $|b^\delta_t(x, \lambda)| \leq \|b_t(\lambda)\|_{\mathbb{L}^\infty}$ for all $x \in \mathbb{R}^d$, so $\|b_t(\lambda)\|_\infty \leq \|b_t(\lambda)\|_{\mathbb{L}^\infty}$.

Finally, for part 4), we once more consider $f = b_t(\lambda)\mathbb{1}_{K'}$. Let $w$ be a function that is supported in $K$ with $\|w\|_{\mathbb{L}^{r'}(K, \mathrm{d}x)} = 1$. Young's inequality with $r_1 = r'/(r' - 1)$, $r_2 = 1$ and $r_3 = r'$ shows that

$$\left| \int_K f^\delta(x) w(x)\, \mathrm{d}x \right| = \left| \int_{\mathbb{R}^d} f^\delta(x) w(x)\, \mathrm{d}x \right| \leq \|f\|_{\mathbb{L}^r(\mathbb{R}^d, \mathrm{d}x)} \,.$$

But from the dual characterization of the $\mathbb{L}^r(K, \mathrm{d}x)$-norm, we get

$$\|f^\delta\|_{\mathbb{L}^r(K, \mathrm{d}x)} = \sup\left\{ \left| \int_K f^\delta(x) w(x)\, \mathrm{d}x \right| \,:\, \|w\|_{\mathbb{L}^{r'}(K, \mathrm{d}x)} = 1 \right\}.$$

Using the definition $f = b_t(\lambda)\mathbb{1}_{K'}$, the first fact we observed above gives

$$\|f^\delta\|_{\mathbb{L}^r(K, \mathrm{d}x)} = \|(\mathbb{1}_{K'} b_t(\lambda))^\delta\|_{\mathbb{L}^r(K, \mathrm{d}x)} = \|b^\delta\|_{\mathbb{L}^r(K, \mathrm{d}x)} \,,$$

and so

$$\|b^\delta\|_{\mathbb{L}^r(K, \mathrm{d}x)} = \|f^\delta\|_{\mathbb{L}^r(K, \mathrm{d}x)} \leq \|f\|_{\mathbb{L}^r(\mathbb{R}^d, \mathrm{d}x)} = \|b_t(\lambda)\|_{\mathbb{L}^r(K', \mathrm{d}x)} \,,$$

where the last equality uses once more the definition $f = b_t(\lambda)\mathbb{1}_{K'}$. This establishes part 4) and completes the proof. □